\documentclass{amsart}

  \def\pb{}

\usepackage{amsmath}
\usepackage{amsxtra}

\usepackage{enumitem}
\setlist[enumerate]{leftmargin=*}
  
\usepackage{amscd}
\usepackage{amsthm}
\usepackage{hyperref}

 \DeclareFontFamily{U}{rcjhbltx}{}
 \DeclareFontShape{U}{rcjhbltx}{m}{n}{<->rcjhbltx}{}
\DeclareSymbolFont{hebrewletters}{U}{rcjhbltx}{m}{n}
 
\DeclareMathSymbol{\shin}{\mathord}{hebrewletters}{152}

\def\Sh{\vartriangle}

\usepackage[T1]{fontenc}

\usepackage{amsfonts}
\usepackage{amssymb}
\usepackage{eucal}
\usepackage[all]{xypic}
\usepackage{color} 
\newtheorem{thm}[subsection]{Theorem}
\newtheorem{cor}[subsection]{Corollary}
\newtheorem{lem}[subsection]{Lemma}
\newtheorem{prop}[subsection]{Proposition}

\theoremstyle{definition}

\newtheorem{defn}[subsection]{Definition}

\newtheorem{rem}[subsection]{Remark}
\newtheorem{rems}[subsection]{Remarks}

\newtheorem{question}[subsection]{Question}

\newtheorem{example}[subsection]{Example}

\newtheorem{problem}[subsection]{Problem}

\usepackage{yhmath}
\DeclareSymbolFont{largesymbols}{OMX}{yhex}{m}{n}
\DeclareMathAccent{\widetilde}{\mathord}{largesymbols}{"65}

\newcommand{\thmref}[1]{Theorem~\ref{#1}}
\newcommand{\secref}[1]{\S\ref{#1}}
\newcommand{\lemref}[1]{Lemma~\ref{#1}}

\newcommand{\defref}[1]{Definition~\ref{#1}}
\newcommand{\propref}[1]{Proposition~\ref{#1}}
\newcommand{\corref}[1]{Corollary~\ref{#1}}

\newcommand{\remref}[1]{Remark~\ref{#1}}
\newcommand{\exref}[1]{Example~\ref{#1}}

\def\qref#1{Question~(\ref{#1})}

\newcommand{\nc}{\newcommand}
\nc{\renc}{\renewcommand}
\nc{\ssec}{\subsection}
\nc{\sssec}{\subsubsection} 
\nc\ol{\overline}
\nc\wt{\widetilde}
\nc\wh{\widehat}

 \def\bl{{{\bf \Lam}}}
\nc{\Aa}{{\mathbb{A}}}
\nc{\Bb}{{\mathbb{B}}}
 \nc{\Gg}{{\mathbb{G}}}  
\renc{\gg}{\mathcal{D}}
\renc{\L}{\mathcal{L}}
\nc{\Hh}{{\mathbb{H}}}
 \nc{\Nn}{{\mathbb{N}}}
\nc{\Pp}{{\mathbb{P}}}
\nc{\Rr}{{\mathbb{R}}}
\nc{\BV}{{\mathbb{V}}}
\nc{\BW}{{\mathbb{W}}}
\nc{\Zz}{{\mathbb{Z}}}
\nc{\Qq}{{\mathbb{Q}}}
\nc{\Ss}{{\mathbb{S}}}
\nc{\Cc}{{\mathbb{C}}}
\nc{\Ff}{{\mathbb{F}}}
\nc{\ff}{{\bar{f}}}
\nc{\EL}{{L_\infty}}
 \def\f{\mathsf{f}}
\nc{\CA}{{\mathcal{A}}}
\nc{\CB}{{\mathcal{B}}}
 
\nc{\CF}{{\mathcal{F}}}

  \def\wz{\mathsf{wz}}

\def\td{^{def}}%
 
\nc{\Las}{\mathsf{Las}}
\def\sl{\mathsf{L}}
\nc{\CG}{{\mathcal{G}}}

 \nc{\R}{{\mathsf{Def}}}

\nc{\CC}{{\mathcal{C}}}
\nc{\CM}{{\mathsf{M}}}

\nc{\CN}{{\mathcal{N}}}
\nc{\Oog}{{\mathbb{O}}}
\nc{\Oo}{{\mathcal{O}}}
 
\nc{\CQ}{{\mathcal{Q}}} 
\nc{\CS}{{\mathcal{S}}}
\nc{\CT}{{\mathcal{T}}}

\nc{\CV}{{\mathcal{V}}}
 
\nc{\CW}{{\mathcal{W}}}
\nc{\CZ}{{\mathcal{Z}}}

\def\CU{\mathcal{J}}
 
\nc{\modec}{{\overset{e.c.}\models}{}}
\nc{\subsetec}{{\underset{e.c.}\subset}{}}
\def\ta{\underset{\alpha}{\times}}

 \def\wo{\widetilde{\Omega}}

\nc{\bvo}{{\underset{\omega}\bigvee}}
\nc{\bwo}{{\underset{\omega}\bigwedge}}

\nc{\cM}{{\check{\mathcal M}}{}}
\nc{\csM}{{\check{\mathcal A}}{}}
\nc{\oM}{{\overset{\circ}{\mathcal M}}{}}
\nc{\obM}{{\overset{\circ}{\mathbf M}}{}}
\nc{\oCA}{{\overset{\circ}{\mathcal A}}{}}
\nc{\obA}{{\overset{\circ}{\mathbf A}}{}}
\nc{\ooM}{{\overset{\circ}{M}}{}}
\nc{\osM}{{\overset{\circ}{\mathsf M}}{}}
\nc{\vM}{{\overset{\bullet}{\mathcal M}}{}}
\nc{\nM}{{\underset{\bullet}{\mathcal M}}{}}
\nc{\oD}{{\overset{\circ}{\mathcal D}}{}}
\nc{\obD}{{\overset{\circ}{\mathbf D}}{}}
\nc{\oA}{{\overset{\circ}{\mathbb A}}{}}
\nc{\op}{{\overset{\bullet}{\mathbf p}}{}}
 
\nc{\oU}{{\overset{\bullet}{\mathcal U}}{}}
\nc{\oZ}{{\overset{\circ}{\mathcal Z}}{}}
\nc{\ofZ}{{\overset{\circ}{\mathfrak Z}}{}}
\nc{\oF}{{\overset{\circ}{\fF}}}

\nc{\fa}{{\mathfrak{a}}}
\nc{\fb}{{\mathfrak{b}}}
\nc{\fg}{{\mathfrak{g}}}
\nc{\fgt}{{\fg}_!}
\nc{\fgl}{{\mathfrak{gl}}}
\nc{\fh}{{\mathfrak{h}}}
\nc{\fj}{{\mathfrak{j}}}
\nc{\fm}{{\mathfrak{m}}}
\nc{\ft}{{\mathfrak{t}}}
\nc{\fn}{{\mathfrak{n}}}
\nc{\fu}{{\mathfrak{u}}}
\nc{\fp}{{\mathfrak{p}}}
\nc{\fr}{{\mathfrak{r}}}
\nc{\fs}{{\mathfrak{s}}}
\nc{\fsl}{{\mathfrak{sl}}} 

\nc{\fA}{{\mathfrak{A}}}
\nc{\fB}{{\mathfrak{B}}}
\nc{\fD}{{\mathfrak{D}}}
\nc{\fE}{{\mathfrak{E}}}
\nc{\fF}{{\mathfrak{F}}}
\nc{\fG}{{\mathfrak{G}}}
\nc{\fK}{{\mathfrak{K}}}
\nc{\fL}{{\mathfrak{L}}}
\nc{\fM}{{\mathfrak{M}}}
\nc{\fN}{{\mathfrak{N}}}
\nc{\fP}{{\mathfrak{P}}}
\nc{\fU}{{\mathfrak{U}}}
\nc{\fV}{{\mathfrak{V}}}
\nc{\fZ}{{\mathfrak{Z}}}

\nc{\bb}{{\mathbf{b}}}
\nc{\bc}{{\mathbf{c}}}
\nc{\bd}{\partial}
\nc{\be}{{\mathbf{e}}}
\nc{\bj}{{\mathbf{j}}}
\nc{\bn}{{\mathbf{n}}}
\nc{\bp}{{\mathbf{p}}}
\nc{\bq}{{\mathbf{q}}}
\nc{\bF}{{\mathbf{F}}}
\nc{\bu}{{\mathbf{u}}}
\nc{\bv}{{\mathbf{v}}}
\nc{\bx}{{\mathbf{x}}}
\nc{\bs}{{\mathbf{s}}}
\nc{\by}{{\bar{y}}}
\nc{\bw}{{\mathbf{w}}}
\nc{\bA}{{\mathbf{A}}}
\nc{\bK}{{\mathbf{K}}}
\nc{\bI}{{\mathbf{I}}}
\nc{\bB}{{\mathbf{B}}}
\nc{\bG}{{\mathbf{G}}}
\nc{\bC}{{\mathbf{C}}}
\nc{\bD}{{\mathbf{D}}}
\nc{\bP}{{\mathbf{P}}}
\nc{\bH}{{\mathbf{H}}}
\nc{\bM}{{\mathbf{M}}}
\nc{\bN}{{\mathbf{N}}}
\nc{\bV}{{\mathbf{V}}}
\nc{\bU}{{\mathbf{U}}}
\nc{\bL}{{\mathbf{L}}}
\nc{\bT}{{\mathbf{T}}}
\nc{\bW}{{\mathbf{W}}}
\nc{\bX}{{\mathbf{X}}}
\nc{\bY}{{\mathbf{Y}}}
\nc{\bZ}{{\mathbf{Z}}}
\nc{\bS}{{\mathbf{S}}}
\nc{\bSi}{{\bar{\Sigma}}}
\nc{\sA}{{\mathsf{A}}}
\nc{\sB}{{\mathsf{B}}}
\nc{\sC}{{\mathsf{C}}}
\nc{\sD}{{\mathsf{D}}}
\nc{\sF}{{\mathsf{F}}}
\nc{\sG}{{\mathsf{G}}}
\nc{\sK}{{\mathsf{K}}}
\nc{\sM}{{\mathsf{M}}}
\nc{\sO}{{\mathsf{O}}}
\nc{\sQ}{{\mathsf{Q}}}
\nc{\sP}{{\mathsf{P}}}
\nc{\sZ}{{\mathsf{Z}}}
\nc{\sfp}{{\mathsf{p}}}
\nc{\sr}{{\mathsf{r}}}
\nc{\sg}{{\mathsf{g}}}
\nc{\sff}{{\mathsf{f}}}
\nc{\sfb}{{\mathsf{b}}}
\nc{\sfc}{{\mathsf{c}}}
\nc{\sd}{{\ltimes}}

 \def\e{\epsilon}

  \nc{\vol}{{\mathop{\operatorname{\rm vol\,}}}}
  \nc{\pr}{{\mathop{\operatorname{\rm pr\,}}}}
\nc{\co}{{\mathop{\operatorname{\rm Core\,}}}}
\nc{\comm}{{\mathop{\operatorname{\rm Comm\,}}}}
  \nc{\gal}{{\mathop{\operatorname{\rm Gal\,}}}}
 
  \nc{\disc}{{\mathop{\operatorname{\rm disc}}}}
  \nc{\Sym}{{\mathop{\operatorname{\rm Sym}}}}
   \nc{\Aut}{{\mathop{\operatorname{\rm Aut}}}}
 \nc{\Spec}{{\mathop{\operatorname{\rm Spec}}}}
  \nc{\spec}{{\mathop{\operatorname{\rm Spec}}}}
\nc{\Ker}{{\mathop{\operatorname{\rm Ker}}}}
 \nc{\dom}{{\mathop{\operatorname{\rm dom}}}}
\nc{\End}{{\mathop{\operatorname{\rm End}}}}
 \nc{\Hom}{\operatorname{Hom}}
 \nc{\GL}{{\mathop{\operatorname{\rm GL}}}}
 \nc{\Id}{{\mathop{\operatorname{\rm Id}}}}
 \nc{\rk}{{\mathop{\operatorname{\rm rk}}}}
 \nc{\length}{{\mathop{\operatorname{\rm length}}}}
\nc{\supp}{{\mathop{\operatorname{\rm supp} \, }}}
\nc{\val}{{\rm val}}
 
\nc{\res}{{\mathop{\operatorname{\rm res}}}}

\def\Ind#1#2#3{{#1} {\downarrow}_{#3} {#2} }

\def\tensor{{\otimes}}
\def\meet{\cap}
\def\union{\cup}

\def\si{\sigma}
\def\g{\gamma}
\def\G{\Gamma}

\def\m{\smallsetminus}

\nc{\seq}[1]{\stackrel{#1}{\sim}}

\def\inv{^{-1}}
\def\claim#1{\smallskip {\noindent \bf Claim #1.\quad}}

\def\beq#1{\begin{equation} \label{ #1}}
\def\eeq{\end{equation}}
\def\normal{\trianglelefteq}

\def\prf{\begin{proof}}
\def\pv{\end{proof} }
 \def\eprf{\end{proof} }

\def\a{\alpha}

\def\bdr{\mathbf{\Delta_r}}

 \renc{\b}{{\beta}}

\def\la{\langle}
\def\ra{\rangle}

\def\Ind#1#2{#1\setbox0=\hbox{$#1x$}\kern\wd0\hbox to 0pt{\hss$#1\mid$\hss}
\lower.9\ht0\hbox to 0pt{\hss$#1\smile$\hss}\kern\wd0}

 \def\Lam{\Lambda}

 \def\cosp{\mathsf{CorSp}}
 \def\cogr{\mathsf{CoGr}}
 \def\appg{\mathsf{CApprGr}}
 \def\appsp{\mathsf{ApprSp}}

\def\lam{\lambda}

\def\U{\mathcal{U}}
  
 \def\Om{\Omega}
 \def\om{\omega}

 \def\ot{\leftarrow}

\def\wl{\widetilde{\Lam}}

\def\lcb{\bar{\lc}}

 \def\tl{\widetilde{\Lam}}
\def\lc{\mathsf H}
\def\al{\G}
\def\hG{ G}
\def\trellis{approximate lattice}

\setlist[itemize]{leftmargin=*}

\title{ Beyond the Lascar group} 
\author{Ehud Hrushovski}
\begin{document}

\begin{abstract}  We work   in a first-order setting where structures are spread out over a metric space,
with quantification allowed only over bounded subsets.  Assuming a doubling property for the metric space,
 we  define a canonical {\em core} $\CU$ associated to such a theory, a  locally compact structure that embeds into the type space over any model. 
The automorphism group  of $\CU$,  modulo certain infinitesimal automorphisms, is a locally compact group $\CG$.  The automorphism groups of models of the theory are related with $\CG$,   not in general via a homomorphism, but by a {\em quasi-homomorphism}, respecting multiplication up to a certain canonical compact error set.  This fundamental structure is applied to describe the nature of approximate subgroups.     Specifically we obtain a full classification of (properly) approximate lattices of $SL_n(\Rr)$ or $SL_n(\Qq_p)$.  
 
\end{abstract}

\maketitle
\setcounter{tocdepth}{1}
\tableofcontents

\section{Introduction}

  This paper consists of three successively more concrete parts,  logically forming a progression.  The first, purely model theoretic, constructs a canonical geometry  associated
  with any first-order theory.  This `core space'  was introduced in \cite{patterns};  here we recast it in a  `local' setting that  allows the core, and the attendant automorphism group, to be locally compact rather than compact.  Beyond that, the relation between 
 the core and the automorphism groups of a models of the theory is now brought out explicitly, in terms of a
 quasi-homomomorphism; the Lascar-Shelah neighbor relation plays the role of the error set. 
  The second part uses this construction  to give a structure theorem for approximate subgroups.  The `Lie model theorem' 
 proved in \cite{nqf} for near-subgroups, and used in \cite{bgt}, is extended in modified form to 
   general approximate subgroups;   a conjecture of Massicot-Wagner predicting that this should be possible is thus confirmed.  The third part specializes to {\em approximate lattices} in locally compact groups.  In the case of semisimple groups, 
  a full classification can be achieved, solving in particular a problem of Bj\"orklund and Hartnick \cite{BjH}.  
  
    The rest of the introduction follows this outline in more detail.   Type spaces -   the Stone spaces dual to the Boolean algebra of formulas     - occur for instance  in the study of countable categoricity.  
  These spaces contain all information about $T$,  but only as a convenient reformulation of the syntax;  they do not  bring out well the geometry of definable sets;  in particular,  they are inert and admit no automorphisms.     In Morley's proof of his theorem on uncountable categoricity, type spaces {\em over models} first played a leading role.  
Morley rank is immediately visible from their topology; and they carry a natural action of $Aut(M)$.  This does not, however, give a canonical space or group associated with $T$, since a choice of $M$ is required.   For some purposes, this is alleviated by uniqueness theorems for saturated models, but then the group
 becomes large and unwieldy,  falling outside any known structured family of groups.   
 
 Shelah \cite{shelah} recognized the algebraic imaginaries as a critical geometry associated with a first-order theory $T$.
The automorphism group is a profinite quotient of $Aut(M)$   (for $M$ a sufficiently homogeneous model of $T$),  and a 
perfect Galois duality holds between them.  For T=ACF one recovers Galois's duality in full.  This was for any theory; for 
stable $T$, Shelah showed how amalgamation of substructures is controlled by his group.  

 Lascar and Pillay \cite{lascar-pillay} defined the {\em bounded hyperimaginaries} as  a generalization
of the algebraic imaginaries; their automorphism group is a compact topological group.   (See also \cite{kim-pillay}; the term Kim-Pillay group, and space, comes from \cite{h-simple}.)   Once again this is valid for any first-order theory; for simple theories, it controls amalgamation.   This equally beautiful theory later became absorbed in continuous logic, and indeed in \cite{byu} is not even terminologically separated from Shelah's algebraic imaginaries.    

It was Lascar, more than a decade earlier, who  suggested that a  richer Galois group is needed in unstable theories.  He defined two candidates; the compact one  used by Kim-Pillay, and a 
   a bigger common quotient of the groups $Aut(M)$, called the general Lascar group.  The latter was much studied recently, but 
 was never found to be the automorphism group of  any meaningful geometry (see e.g. \cite{kms} for a related negative result.)
 Krupi\'nski, Pillay, and Rzepecki \cite{kpr}, \cite{kns} , working   group-theoretically, showed suggestively that   the general Lascar group was in turn a quotient of a compact group;    
 see  the introduction to \cite{patterns} for more detailed references.    
%
 
 In \cite{patterns},  
 a common subspace $\CU$ of all type spaces over models was identified,  generalizing the Shelah space and Kim-Pillay space for stable and simple theories, and having significant structure beyond the topology.   While canonical and closely related to definability in $T$, the automorphism groups $Aut(M)$ of models of $T$ do not naturally act on $\CU$.  We show here, however, that they do admit a useful {\em quasi-action}.   A {\em quasi-homomorphism} $\phi: G \to H:K$ 
between two groups $G,H$, relative to an `error set' $K \subset H$, is  a map $\phi$ satisfying 
$\phi(1)=1$ and 
$\phi(xy\inv) \in \phi(x)\phi(y\inv)K$ for all $x,y \in G$.  
This is of interest when $K$ is small in some sense.   
We find a quasi-homomorphism  $\varphi: Aut(M) \to Aut(\CU):\Sh$.   The error set $\Sh$ derives from a fundamental `neighbor' relation going back to \cite{lascar} (on $Aut(M)$
for highly saturated $M$) and \cite{shelah-simple} (for indiscernible sequences), and defined here on the core $\CU$.

%
 Given a quasi-homomorphism $\phi: G \to H:K$, one can of course factor out the group $\la K \ra$ generated by $K$
 so as to obtain an  ordinary homomorphism $G \to H/\la K \ra$.  In the case of the quasi-homomorphism $Aut(M) \to Aut(\CU):\Sh$ described above,   this would recover the general Lascar group. 
  
  The stance taken here is that factoring out by $\la \Sh \ra$ is the wrong thing to do.      An alternative way to forget 
 unwanted information while discarding much less of value is offered by category theory:  in place of changing the objects, add (iso)morphisms.   This is set out   in Appendix \ref{categories}.   In the paper itself we will use the categorical approach only in spirit.    
  We will however analyze mathematically the relevant morphisms. 
  
    Examples of quasi-homomorphisms include projective representations of a group $G$; here $H$ is a unitary group, the group
 of automorphisms of a Hilbert space, and $K$ is the one-dimensional center of $H$.    
Quasi-homorphisms with unitary targets were also studied notably in  \cite{turing}  (a 1938 paper first cited, according to  MathSciNet, in 2011), and in \cite{kazhdan},  with $K$ a set of elements of  small operator norm in both cases.      When $G$ is finite (Turing) or amenable (Kazhdan), they were shown to be close to ordinary homomorphisms.    

 The quasi-homomorphisms of interest to us will have  locally compact target groups, and normal, 
 compact error sets $K$.
The case $H = \Rr$ has been very well studied; a quasi-homomorphism into $\Rr$, with compact error set, is called a {\em quasimorphism}; it is significant in several branches of geometry.    (Readers unfamiliar with the notion may enjoy reading
the very short introduction \cite{kotschick} at this point).  


 Returning to first-order theories,  we can get a clear idea of the smallness of $\Sh$ if we allow a slight generalisation,
 to local logic.  This is introduced in \S 2.  In many-sorted logic, quantification is allowed only on a given sort, or finitely many sorts.   If one  thinks of a sort as compact  (say via the type space), then a many-sorted model is only locally compact.
  Local logic is a slight variation, likewise restricting quantification, but  refraining from destroying possible symmetries by naming specific bounded domains as sorts.  In the simplest case,
 one has a metric $m$ (a discrete metric is good enough for our purposes); and quantification is allowed over balls of bounded $m$-radius. 
 We will assume a `doubling' property - a $2$-ball  for the metric $m$ is covered by finitely many $1$-balls.  
   The construction of $\CU$ generalizes;  $Aut(\CU)$ is now locally compact rather than compact; it has a natural locally compact Hausdorff quotient, which is what we will actually work with.  The relation $\Sh$ on $\CU$ remains local:  
 $x \Sh y$ implies that their distance is at most $1$.   As a results, the set of automorphisms $\Sh^{Aut(\CU)}$   is a  normal,   compact  subset of $Aut(\CU)$.     In fact,    
 compactness alone does not fully capture the tightness of $\Sh$.    But along with conjugation invariance - arising from the 
 canonical nature of $\Sh$ - it is  sufficient for most of the group-theoretic applications. 
 
 At this point, the reader is encouraged to consult \thmref{summary1}, summarizing the main properties of the locally compact core $\CU$, and \thmref{qh-basic},
including the basic properties of the quasi-homomorphism $\varphi$.    The material is closely parallel to the non-local case of \cite{patterns}.    The exposition in sections 2 and 3 is a compromise, not repeating every lemma of \cite{patterns}, but attempting to be self-contained with respect to the main line and in particular the results that will be critical in the later sections.

In \secref{definablegroups}, we transpose to the case of definable groups.  In the local setting, we look at groups that are direct limits of definable sets; in the model theory literature they are known as piecewise-definable, locally definable, or strict ind-definable.
Thus $G=\union X_n$, where each $X_n$ is a definable set, and multiplication restricts to a definable map $X_n \times X_n \to X_{n+1}$ .  
The doubling condition amounts here to the assumption that $X_{n+1}$ is covered by finitely many translates of $X_n$
This makes each $X_n$ into an {\em approximate subgroup} of $G$, as defined in \cite{tao}.  
At least at first approach, approximate subgroups are best   understood up to commensurability,  i.e. without distinguishing $X$ and $X'$, if each can be covered by finitely many translates of the other. 

We will call an  approximate subgroup {\em {laminar}} if it contains a  sequence $(X_n)$ as above
indexed by {\em negative} integers.    A 
 commensurability class $\om$ is  {\em {laminar}} if it contains a laminar  approximate subgroup.   
Model-theoretically, $X$ is laminar iff it contains a subgroup cut out by an intersection of generic definable subsets.  
Quotienting such a a $\bigwedge$-definable subgroup out from a saturated version of $G$ leads to a homomorphism $\varphi: G \to \lc$
into a locally compact group $\lc$; and $\om$ can be recovered from $\varphi$ and $\lc$, since the  elements of $\om$ are intertwined
with pullbacks of compact open neighborhoods of $1$ in $\lc$.   
 We obtain something like a Bourgain system (\cite{green-sanders})  down to an 
arbitrarily {fine} scale.  In the non-laminar case, to extend the metaphor, eddies of a certain size prevent further descent.

We call $X$ {\em amenable} if there exists a translation invariant measure on $X^3$; more generally it is a {\em near subgroup} if
an appropriate invariant ideal on definable subsets of $X^3$ exists.   For   amenable approximate subgroups and and near subgroups it was proved
in \cite{nqf}, \cite{sanders}, \cite{massicot-wagner},  in various settings, that the commensurability
class of $X$ is {laminar}; in fact the sequence $X^n, n \geq 4$ can be continued downwards through the negative integers.

Here we prove, with no amenability assumptions, that {\em a map $\varphi: G \to \lc:K$ into a locally compact group still exists, with the same  property of intertwining {\em compact} with {\em definable}.}  Only now $\varphi$ is  a quasi-homomorphism,  a  compact, normal error set $K$ is allowed; 
and open neighborhoods of $K$ (rather than of $1$) are considered on the locally compact side.   See \thmref{grmain}. 

Model-theoretically, the characteristic definability property of $\varphi$ in the {laminar} case was {\em definable separation}:  
if $C,C'$ are disjoint compact subsets of $\lc$, then $\varphi\inv(C),\varphi \inv(C')$ are separated by a definable set.
If we give $\lc$ a metric, $C,C'$ are separated by some distance $\e>0$.
In the general setting, we still have the separation property but now assuming $CK,C'K$ are disjoint; they must  be seen to be distinct even if our resolution is limited to the radius of $K$.  
In the {laminar} case $\e>0$ is arbitrary, but in general a certain minimum separation is required.

Massicot and Wagner wrote, regarding approximate subgroups:   "We conjecture that even without the definable amenability assumption a suitable Lie model exists"  \cite{massicot-wagner}.     The term `Lie model' was coined earlier in \cite{bgt}, where it was the starting point for the classification of pseudo-finite  approximate groups.
The meaning of `suitable' was not further clarified, but 
I  view \thmref{ag2} as a confirmation of their intuition.   It asserts that {\em any commensurability class of approximate subgroups
arises by pullback of compact neighborhoods under an quasi-homomorphism into a Lie group $L$; moreover, the error set is contained
 a bounded subset of a closed, normal abelian subgroup $A \normal L$, $A \cong \Rr^N$, such that $L$-conjugation respects a Euclidean structure on $A$.}    (We refer to such actions as {\em rigid}.)

It was not immediately apparent 
 whether \thmref{ag2} has real content or merely replaces  one kind of approximateness (for subgroups) by 
another (for homomorphisms); the problems attacked in part 3, described below, were intended as a testing ground.  The  latter turned out to be a much tighter condition, thanks to the combined compactness and normality of the error set and the finite dimensionality of the target.  To bring this out directly we give two additional formulations within \S 5,
a group-theoretic and a model-theoretic one.   The group-theoretic version, \thmref{ag2gr}, points out that the quasi-homomorphism
$\varphi: G \to \lc:K$ of \thmref{ag2} can be  viewed as an (ordinary) homomorphism from an `almost' central   extension $\tilde{G}$  of $G$ by  $\Rr^n$, into the Lie group $\lc$.  Here  `almost'   means that the image in $End(\Rr^n)$ of conjugation by elements of $\tilde{G}$ ,
 while not trivial,  is contained in a compact subgroup.  
This resembles the situation with projective representations of Lie groups, seen as representations
of a central extension.   The cocycle coding the group extension will be a bounded one, forming one connection with the far-reaching results of bounded cohomology (e.g. \cite{gromov},  \cite{brooks},  \cite{burger-monod},\cite{monod-icm}.)

The second formulation is \thmref{ag2mt}.  The model-theoretic avatar of an approximate subgroup in the {laminar} case is an $\bigwedge$-definable group `of bounded index.'  Here we show that in general, an approximate subgroup is associated with
an $\bigwedge$-definable set, obtained as the `approximate kernel' of finitely many (homogeneous) quasimorphisms on an $\bigwedge$-definable group $\G$ of bounded index, i.e.  the   pullback of a compact neighborhood of the error set in $\Rr^n$.  This connects   again,  not in quite the same way, to bounded cohomology.

  For  instance,  
since groups of bounded exponent admit neither dense homomorphisms into connected Lie groups nor unbounded quasimorphism, 
 any approximate subgroup of a bounded exponent group is commensurable to an actual subgroup;  the conclusion is the same 
  as in \cite{nqf}, where this was used as an illustration of the theory, but the hypothesis includes no amenability assumption. 
  
The $\bigwedge$-definable group $\G$ is the kernel of a continuous homorphism into a locally compact group
   with respect to the logic topology; 
 the quasimorphism $\G \to \Rr^n$ is also continuous; so that 
 the bumpiness alluded to above is felt only at the scale of the error set. 
  
The (symmetrized) approximate kernel of a nontrivial quasimorphism  is always an approximate subgroup (\propref{converse}), but never a {laminar} one (\propref{notsame})

The third part of the paper, consisting of \S 6, \S 7 and Appendix A,  is concerned with the geometry of numbers; but this will not be obvious at the outset.    
Our working environment in \S 6 is an easily defined class of locally compact groups, that we call {\em abstractly semisimple}.  It has two defining requirements:   no abelian normal subgroups, no discrete conjugacy classes other than the trivial ones.   
These serve as ambient groups; within them, we consider discrete approximate subgroups.   (Perhaps it is worth recalling here, though I am aware of  no direct connection,   the fundamental role of discreteness in the theory of pseudo-finite approximate subgroups of \cite{bgt}.)  
We call a  commensurability class of discrete approximate subgroups {\em strictly dense} if the group generated by any element of the class is dense; this is at an opposite extreme from discrete subgroups;   under certain assumptions, a decomposition exists  into the two kinds.   Building on the results of \S 5,   we  begin by showing  
 that, for abstractly semisimple groups,    {\em a strictly dense approximate subgroup  is commensurably contained in a laminar one.}  
 See  \thmref{discrete1}.

Recall that a {\em lattice} is a discrete subgroup  of finite co-volume; it is called  {\em uniform}   if it is  co-compact.  In the case of the Euclidean space $\Rr^3$,  a lattice is the accepted mathematical model of a crystal, and their 
histories are inseparable.   Yves Meyer in  \cite{meyer} defined, in the setting of $\Rr^n$, {\em approximate lattices};
this later became the mathematical home for quasicrystals.   
 Meyer  proved
that they are {laminar} as approximate subgroups, with Lie model $\Rr^n$;    he did so in the equivalent, geometrically suggestive language of {\em model sets}. 
  On the other hand, lattices in non-commutative locally compact groups such as $SL_n(\Rr)$ play a major role in many parts of geometry, ergodic theory and number theory;   especially relevant to us are the arithmeticity 
   theorems of \cite{margulis}.    
  
  Bj\"orklund and Hartnick \cite{BjH} were first to investigate approximate versions.   They defined uniform approximate lattices,  and offered two ``tentative" definitions in the non-uniform case.   We will suggest another definition, defining 
   an {\em   approximate lattice in the sense of finite covolume} to be a discrete approximate subgroup $\Lam$ of $G$,   such that $\Lam C =G$ for some Borel set $C$ of finite measure.    As  in \cite{BjH},  $\Lam$ is said to be {\em uniform} if $C$ can be taken compact.   
   It seems to work;
the basic properties of approximate lattices, with this elementary definition, are developed in Appendix \ref{approxlattices}.   
After writing this, I ran into
Problem 6.4 of \cite{bats}, asking for the ``right" definition of an approximate lattice; Appendix \secref{approxlattices}, then, is our possibly naive proposal. 
In this paper,  approximate lattice  will always mean, by way of abbreviation, `approximate lattice in the sense of finite covolume'.   

Among the approximate lattices,  Bj\"orklund and Hartnick consider  the {\em Meyer sets}.  The definition (that we naturally extend from the uniform case) inolves not only $G$ but also an additional locally compact
$\lc$, that we will call the complementary group, and  a lattice $\G$ in $G \times \lc$.  A Meyer set $M$ is defined by the positive primitive formula of {\em local}  logic:  $(\exists h \in C)((x,h) \in \G)$, where $C$ is a compact neighborhood of $1$ in $H$. 
This description is equivalent to laminarity of $M$, with $\lc$ corresponding to the target of the homomorphism into 
a locally compact group.   
 
    I am grateful to Emmanuel Breuillard    for introducing me to this beautiful theory, and explaining this equivalence between the language of laminarity, or Lie models, and of Meyer model sets.     
       
Bj\"orklund and Hartnick posed the problem of  proving laminarity when possible (\cite{BjH}, Problem 1,  formulated there for 
uniform approximate lattices.)     Simon Machado, in  \cite{machado}, \cite{machado2}, did so in the nilpotent and solvable cases,
and for amenable groups.  In \secref{discreteapprox}, complementing this work, we solve Problem 1 for   semisimple groups $G$.   
To do so,  we will use a local logic where `local' implies 'discrete'.   This contrasts with the natural treatment of
say $SL_n(\Rr)$ in local logic, where ``local" means "compact" (see e.g. \cite{HPP}, section 7.)  
 Nevertheless the local compactness of $Aut(\CU)$ (as mediated by Theorems 4.11,5.10) 
   allows us to recover compactness, effectively 
 finding the complementary group and presenting the 
 approximate lattice as a Meyer set.   See  \corref{semisimple1}.  
 The duality between compactness and discreteness seen here is a key feature of this world.  
 
 At this point, Margulis arithmeticity becomes available.   See  {arithmodelset} for the definition of an {\em arithmetic approximate lattice}.  It is just like the definition of an arithmetic lattice in \cite{margulis}, except that 
 (as if in line with the recommendation of Hensel and Weil) 
  all places are treated as equals; no special provision is made for archimedean primes.  If we switch attention 
for a moment to rings, 
$\Zz$ is the simplest lattice (in $\Rr$), while the approximate ring $\{a/p^n \in \Zz[1/p]: |a| \leq p^n \}$ is the simplest 
 arithmetic approximate lattice (in $\Qq_p$).  In each case, the (approximate) lattice $R$ is the set of rationals whose norm in every completion, other than the ambient locally compact group, is at most $1$.    And in both cases, appropriately understood,
$G(R)$ gives rise to arithmetic lattices in $G(\Rr)$ or $G(\Qq_p)$ for semisimple groups $G$.  
   \thmref{arith1} shows that {\em  any  approximate lattice  in $G$ is commensurable with a product  of lattices and arithmetic approximate lattices  of direct factors of $G$. }     Note that the arithmeticity of the strictly approximate lattices is valid even for $n=2$; higher rank assumptions are not required.
%
   
  To illustrate the concrete character of this answer,  let us describe, up to commensurability, the approximate lattices  $\Lam$ in $SL_n(\Rr)$.
   They  depend  on a choice of a  number field $K \leq \Rr$ and an algebraic group $G \leq GL_N$ defined over $K$,
   and isomorphic to $SL_n$ over $\Rr$ (the case $G=SL_n$, $n=N$ is interesting enough).       Let $\Lam_K$ be
 the set of matrices in $G(K) \leq GL_N(K)$  of the form $I+ M$, where $I$ is the identity matrix and $M$ is a matrix whose nonzero entries are Pisot numbers $\alpha$ with $\Qq(\alpha)= K$.  Then $\Lam_K$ is an approximate lattice in $SL_n(\Rr)$.     
It follows from \thmref{arith1} that {\em any approximate lattice in $SL_n(\Rr)$ is commensurable to a conjugate of one of these, or to a lattice.}    
  
In positive characteristic, 
  strictly approximate lattice seem unnecessary, and in fact, do not exist (\corref{char-p}).  

 
The class of abstractly semisimple groups was introduced in \S 6 as simply  a convenient 
 home for  groups like $SL_n(K)$, with $K$ a local field, and their more demanding cousins; but  it eventually becomes natural to take interest in the class itself.     \thmref{arith2} 
shows that all strictly approximate lattices in abstractly semisimple groups are of arithmetic origin.  In particular, 
by  \remref{arith2b},    {\em an arbitrary abstractly semisimple group containing a strictly approximate lattice is virtually isogenous an adelic product of semisimple groups over local fields}, where
`virtually isogenous' allows homomorphisms with co-compact image and compact kernel.  
 The proof of this recognition theorem 
 uses Margulis arithmeticity again, and a method of pivoting over the complementary locally compact.
   
 Finally, we take a peek beyond soluble or semisimple groups.  We construct a simple example of a non-laminar \trellis.  Returning to bounded cohmology, we also give a criterion for an  approximate lattice $\Lam \leq G$, where $G$ is semisimple,   ensuring  laminarity of any approximate lattice $\Lam'$  lying above $\Lam$ in an algebraic group $G'$ with semisimple part $G$.  
  
 Leitfaden:   The appendices are, within the present paper, self-contained.
 Sections (2,3,4) form a logical sequence.  The main results of the later sections require, from this development, only parts (1,2) of  \thmref{grmain}. 
  The reader interested mostly in approximate lattices can read  Appendix \ref{approxsg}, the first few pages of \secref{approxsg}, including  \thmref{ag2}, and turn to \secref{discreteapprox}. 
  
    Questions of uneven scope and difficulty are scattered through the text, mostly at the end of the sections.    
 
   Thanks to Arturo Rodriguez-Fanlo    and Anand Pillay for their reading and comments, and to Yves de Cornulier for explaining
 the literature around \propref{lc1cc}.  
 

Immediately after writing this paper  I became aware that Simon Machado had been working on the structure of
approximate  lattices  in semisimple groups at the same time.   He
 obtained similar results including  arithmeticity and rigidity results for strong approximate lattices in higher rank Lie groups, via what appears to be a beautiful and completely different approach via   Zimmer cocycle rigidity;  see \cite{machado4}.    
 See also the introduction to \S A.  The basic definitions differ slightly and still need to  be compared in detail; see in particular the comments added to  Questions \ref{othergroups-q} and   \ref{bjhquestion}.      
 
 Machado also drew my attention, while I was lecturing on these results, to
 the paper \cite{bfs} of Bader, Furman and Sauer;   it is very close in both statement and proof to  \thmref{arith2} and \remref{arith2b}.
 I left the treatment as it was to  illustration of the naturality of the  approximate lattice approach,  and with future generalizations beyond the laminar framework in mind.

 While the results that approximate lattices in semisimple groups are close to lattices (in other groups) are very satisfactory,
 they also mean that the additional flexibility of approximate lattices must find use outside the semisimple setting,  or beyond the lattice framework.  
 In particular Kazhdan suggested investigating  discrete approximate subgroups in $GL_n$ over local fields, with respect to growth rates of  the number of points in increasing balls.  Can growth rates be achieved for approximate subgroups that are not possible for subgroups? 
 
 Moving to the approximate lattice setting allows investigating an  arithmetic lattice  in $GL_n(\Rr)$ within any $GL_n(\Qq_p)$;
 it suffices to extend to an adelic lattice, cut with neighborhoods of the identity in the other places and project to $GL_n(\Qq_p)$;
 little information is lost, and in particular the commensurability class of the original lattice can be reconstructed.  
 Taking ultraproducts, and noting that the degree $k$ of approximateness does not depend on $p$,  one can even obtain 
 such information within  linear groups over equicharacteristic Henselian fields $K$.   Can some of the theory of these lattices be seen more clearly when working in $GL_n(\Qq_p)$ or $GL_n(K)$?

 \section{Local logic} \label{local}

    
 We begin by reviewing the theory of existentially closed and saturated models.  This theory belongs to the model theory of the 1960's, and was already recalled in detail in \cite{patterns}.  But we need a  generalization here,  permitting a  localization of the quantifiers without disturbing a possible symmetry.  The extension is   slight in that upon adding a constant symbol $c$ in
 ${\gg}$, we return to the  previous positive-primitive setting, with balls of various radii around $c$ are viewed as sorts.
 As usual more care is needed to subtract a constant than to add one.
 
Throughout this paper,  unless  explicitly specified otherwise, 
 the word `definable' will mean:  by a formula of the given language, without parameters.

  We assume given a relational language  (with equality) ${\L}$, with a distinguished sort ${\gg}$, and a distinguished family $\{\mu\}$ of binary relations   on ${\gg}$, called {\em locality relations}.   For simplicity, we will think of  ${\gg}$ 
  as the only sort, and 
  assume we have only one locality relation $\mu_1$. 
      For any ${\L}$-structure $A$, we define a   metric $\mu$ on $A$ according to the Cayley graph of $\mu_1$:   
      \[\mu(x,y) \leq 1 \iff  x=y \vee \mu_1(x,y) \vee \mu_1(y,x) \]
    \[ \mu(x,y) \leq n+1 \iff   (\exists z)( \mu(x,z)\leq 1 \wedge \mu(z,y) \leq n ) \]

We write $\mu_n(x,y)$ for $\mu(x,y) \leq n$ (replacing $\mu_1$ by its symmetrization),
 and also $\mu_n(a):=\{y: \mu_n(a,y)\}$.

 By a pp formula,we mean one built from atomic formulas using conjunction and {\em local existential quantification}:  if $\phi(x,y,\cdots)$
 is a pp formula with a free  variable $x$ and another variable or constant $y$ (both of sort $\gg$), 
then  $(\exists x)(\mu_1(x,y) \wedge \phi)$ is also pp.     Thus iterated quantification ranges over a ball of bounded radius.   There are no pp sentences,  though they do appear as soon a a constant is named.

A {\em primitive-universal theory} $\CT$ includes negations of sentences of the form $(\exists x)(\bigwedge_{i=1}^m \psi_i)$,
where $x$ is a tuple of variables, and $\psi_i$ a tuple of atomic formulas, or equivalently a pp formula; thus a typical sentence
can be written $(\forall x) \neg \psi(x)$.   
A {\em primitive-existential theory}   includes sentences of the form $(\exists x)(\bigwedge_{i=1}^m \psi_i)$.
Given a primitive-universal theory $\CT$, let $\CT_{\exists} $ consist of all existential sentences $(\exists x) \psi(x)$ where $\neg (\exists x) \psi(x) \notin  \CT$.
 
 A {\em  weak model} of $\CT$ is a structure $A$ where the sentences of $\CT$ hold.  It is a {\em model} of $\CT$, in the sense of local logic, if it is a weak model where any two elements of ${\gg}(M)$  are at finite $\mu$-distance.   We will sometimes 
use the term {\em local model} for emphasis and clarity.

We say LJEP holds for $\CT$ if  for any two local models $A,B$  of $\CT$, there exists a local model $C$ of $\CT$
and homomorphisms $A \to C$ and $B \to C$.      In this case,    $\CT ^{\pm} := \CT \union \CT_{\exists}$ is consistent,
(and indeed any existentially closed model of $\CT$ is a model of $\CT^{\pm}$). 


In our application to the core, LJEP will hold   in a strong  form,   (LJEP$_k$):     If $(\exists x) \psi(x), (\exists y) \psi'(y) $ are in  $\CT_{\exists}$, where $\psi,\psi'$ are pp and $x,y$ are variables 
 (of sort $\gg$),    then  
 $(\exists x)(\exists y)(\mu_k(x,y)        \wedge \psi(x) \wedge \psi'(y) ) \in \CT_{\exists}$.
 
   In the main application, in fact, we will have LJEP$_2$; in case of transitivity, when ${\gg}$ has a single 1-type,  LJEP$_0$ holds. 

LJEP$_k$ implies LJEP:    If $A,B \models \CT$, pick $a \in {\gg}(A), b \in {\gg}(B)$;
  then the union of the diagrams of $A$ and of $B$, along with the sentence $\mu_k(a,b)$, is consistent.  Any weak model
  of this union, 
  restricted to the universe spanned by $A,B$, will be a model of $\CT$.
  
   \noindent{{\bf Remarks on LJEP}     One can define a distance on maximal pp types, namely $d(p,q) \leq n$ if 
  $\mu_k(x,y) \wedge p(x) \wedge q(y)$ is consistent with $\CT$.   Define an equivalence relation on types, $d(p,q)<\infty$.  Then
  LJEP holds iff 
   any two types are 
  equivalent.    Without the joint embedding property, the results of this section still hold, except for the uniqueness of the universal model in \propref{prop2.1}; 
  instead, there is one universal model realizing each equivalence class.

Local logic ultrapowers of a model $M$ of $\CT$ are formed by first taking the ordinary ultrapower, then restricting
to the elements at finite distance from some element of $M$.      We can similarly define ultraproducts of  models $M$ of $\CT_c$, i.e. models $M$ of $\CT$ with a distinguished element $a=c^M \in {\gg}(M)$.   

A model $A$  of $\CT$ is {\em existentially closed} (abbreviated e.c.) if ${\gg}(A) \neq \emptyset$, and for every homomorphism $f: A \to B$,  where $B \models \CT$, and 
any ${\L}_A$-{pp} sentence $\phi$, if $B_A \models \phi$ then $A_A \models \phi$.  Here ${\L}_A$ is ${\L}$
expanded by constants for the elements $a \in A$; they are interpreted as $a$ in $A_A$ and as $f(a)$
in $B_A$.   

 In particular, if $A$ is e.c. then $A$ is {\em simple}, i.e. every  homomorphism  $f: A \to B$ is an $\L$- embedding (an isomorphism onto the image.)  
Thus $A$ is e.c. iff it is simple, and existentially closed with respect to embeddings.  

It is clear from Zorn's lemma that any  model $A$ of $\CT$ admits a homomorphism into a simple one.    
The usual direct limit construction shows that any simple model $A$ of $\CT$ admits an embedding 
 into an existentially closed model $B$ of $\CT$, with $|B| \leq |{\L}_A| + \aleph_0$.     Write $M \modec \CT$ if $M$ is an e.c. model of $\CT$.    Also, for any local sentence $\psi$, write $\CT \modec \psi$ if 
$M \models \psi$ whenever $M \modec \CT$.
 
The class of e.c. models  admits amalgamation:   if $f_i: A \to B_i$, we may embed each $B_i$ in an ultrapower $A^*$ of $A$, then compose with a homomorphism to an e.c. model.  

An e.c.-model $\U \models \CT$ is {\em $\lambda$-saturated} if whenever  $A$ is an e.c. substructure of $\U$ and $h: A \to B$ an embedding with  $B \modec \CT$ with $|B| < \lambda$, there exists an embedding $g: B \to \U$ with $g \circ h = Id_A$.

For pp formulas $\phi(x),\psi(x)$ we write $\phi \perp \psi$ if $\CT \models (\forall x)\neg (\phi \wedge \psi)$.  
If $M \modec \CT$, for any pp formula $\phi(x)$ and any $a$ from $M$, we have $M \models \phi(a)$
or $M \models \psi(a)$ for some pp formula $\psi \perp \phi$.  This is proved just as in \cite{patterns} 2.3(1),
taking into account that pp formulas now use local quantifiers only.  

 Write $\phi \subsetec \psi$ if $\theta \perp \phi$ whenever
$\theta \perp \psi$; in e.c. models we have     $\phi \subsetec \psi$ iff $\phi(M) \subseteq \psi(M)$.  

We will be interested in a generalized finiteness condition on $\CT$:  

(ecB):  For some infinite cardinal $\kappa$, any e.c. model of $\CT$ has cardinality $\leq \kappa$.   

Here is a more concrete form of this  condition; we will name it positive primitive compactness.  Note the similarity, in complete metric setting, to the 'totally bounded' definition of local compactness;
here one thinks of $\neg \phi$ as an open neighborhood of the diagonal, with $\phi$ ensuring a distance bounded away from $0$, in any given bounded subset.

({ppC}):  Let $\phi(x,y)$ be a pp formula, implying $\mu_n(x,y)$ for some $n$, and with $\phi \perp \  =$; i.e. $\CT \models \neg (\exists x) \phi(x,x)$.     Then for some $k=k(\phi) \in \Nn$,
$\CT \models \neg (\exists x_0,\cdots,x_k) \bigwedge_{0 \leq i < j \leq k} \phi(x_i,x_j) $.  \label{ppC}

Walter Taylor isolated this condition (in the non-local setting) using the Erdo\"s-Rado theorem,  in connection with the Hanf number for e.c. models; i.e. he proved the equivalence
of (ecB) with (ppC).  
If ({ppC}) fails, it is  clear from compactness  that there exist arbitrarily large e.c. models of $\CT$; indeed they contain
arbitrarily long sequences $(a_i)$ with $\phi(a_i,a_j)$ for all $i< j$, precluding any possible equality $a_i = a_j$.    Conversely if 
 ({ppC}) holds, (ecB) follows  with $\kappa=2^{\L}$;  
 This is due (in the non-local setting) to Walter Taylor.
 we    prove this below by  constructing the universal model without  any mention of cardinals bigger than $2^{\L}$.

\begin{prop}\label{prop2.1}  Let $\CT$ be a local primitive-universal theory  enjoying ({ppC}). 
 Then:  \begin{enumerate}
\item  Any $M \models \CT$ admits a homomorphism to a unique $|\L|^+$-saturated e.c. model $\U$.
If $M \modec \CT$ then $\U$ is unique up to isomorphism over $M$.

  \item  $\U$ is maximal:   if $\U \leq B \modec \CT $ then $U=B$.   If $f: \U \to B \models \CT$ is a homomorphism, 
     there exists a homomorphism $r: B \to \U$ with $r \circ f = Id_\U$.  
 \item Any homomorphism on $\U$ is an embedding, and 
    any endomorphism $f: \U \to \U$ is an isomorphism. 
 \item  $\U$ is    homogeneous  for {pp} types.  
   \item   Any e.c. model $A$ of $\CT$ satisfies $|A| \leq 2^{|\L|}$.
\item  if $\CT$ has LJEP, then $\U$ is universal, i.e.   , i.e. any $M \models \CT$ admits a homomorphism $M \to \U$; and
$\CT$ is the unique universal e.c. model of $\CT$.
\end{enumerate}
\end{prop}

\prf   We may assume $M \modec \CT$, by finding a homomorphism from $M$ to an e.c. model.  
Existence of an $|\L|^+$-saturated e.c. model $\U$, embedding $M$, with $|\U| \leq 2^{|\L|+|M|}$, now
 follows from the usual construction of saturated models, using amalgamation repeatedly and taking direct limits.

\claim{} Let $\U$ be an $|\L|^+$-saturated e.c. model.   Then $\U$ is maximal. 
\prf   Suppose, for the sake of contradiction, that $B \modec \CT$, $\U \leq B$, and $b \in B \setminus U$.  By $|\L|^+$-saturation, one can find a sequence
$(a_i: i < |\L|^+)$ with $a_i \in \U$ such that $a_i$ has the same pp type as $b$ over $\{a_j: j<i\}$.  Since $a_i \in \U$
we have    $a_i \neq b$,
and so since $B$ is e.c. there exists a pp formula $R_i \perp \ =$ with $R_i(a_i,b)$ and such that $R_i$ implies $\mu_{n}$
for some $n=n_i$.   By refining the sequence we may
assume all $R_i$ are equal to some $R$, and all $n_i=n$.      Then $R(a_j,a_i) \wedge \mu_{2n}(a_i,a_j)$ holds for all $j<i$,   contradicting ({ppC}).
\eprf

Thus the first part of (2) holds.
Now uniqueness in (1) follows too: if $M$ admits a homomorphism into two such e.c. models $\U$ and $\U'$, by 
amalgamation over $M$ they embed over $M$ into a third model $\U''$; but then by the claim the embeddings $\U \to \U''$
and $\U' \to \U''$ are isomorphisms (over $M$.)

(3)  Any homomorphism from any $A \modec \CT$  into a model of $\CT$ is an embedding, 
by definition of e.c..   

 Let $f: {\U} \to {\U}$ be an endomorphism,
with image $U' \leq U$.   Since $\U$ is e.c., $f$ is an isomorphism $\U \to U'$;
hence $U'$ is also maximal.  Since $U' \leq \U$,   we have  $\U=U'$.
Thus $f$ is an automorphism of $\U$.

As ${\U}$ is universal, if $\U \leq B \models \CT$  there exists a homomorphism $f: B \to {\U}$;
 on ${\U}$ it induces an isomorphism $g$;  so $r=g \inv \circ f: B \to {\U}$ is as required in the second sentence of (2).
 
 (4)  
 Homogeneity for pp types:  let $C,C' \subset \U$ and let $g: C \to C'$ be a bijection preserving pp formulas.
 We may assume $C \neq \emptyset$, say $c \in C$.
 Using locality of $\U$, any existential sentence  $\theta$ true in $\U_C$ is implied by a stronger pp sentence, specifying
 also a distance of the witnesses from $c$.  Hence $\theta$, transposed via $g$, is also true in $\U_{C'}$.  
 Thus $g$ extends to an embedding $G$ of $\U$ in some elementary extension $\U^*$ of $\U$ (a weak model initially).
 Restricting to elements at finite distance from $g(c)$ will give a local model $V$ including all of $g(\U)$. 
 We also have the usual inclusion of $\U$ in $\U^*$, and again it is contained in $V$.  
  Let $r: V \to \U$ be a retraction.    Then $r \circ G: \U \to \U$ is an automorphism extending $g$.

  For (5), we may   assume that $A$ is $|\L|^+$-saturated.   By L\"owenheim-Skolem,  $A$ contains an  $|\L|^+$-saturated e.c. substructure $A'$ of  cardinality at most $2^{|\L|}$.  By the Claim we have $A=A'$.
 
 (6) Now assume LJEP.  Then any two $|L|^+$-saturated e.c. models  $\U,\U' \models T$ are isomorphic.  Indeed
 by LJEP  $\U,\U'$ both embed into some $\U''$, and then $\U \cong \U'' \cong \U'$ as in the proof of uniqueness in (1).  
 
 Any universal   $U \modec \CT$ is isomorphic to  $\U$:   since  $U$ is universal,
 there exists a homomorphism $f: U \to \U$, which must be an embedding; so we may assume $U \leq \U$.
 Then there exists a retraction $r: \U \to U$.  But endomorphisms of $\U$ are isomorphisms, so $\U \cong U$.
 
\eprf

In the absence of LJEP, we can define a distance on maximal pp types, namely $d(p,q) \leq n$ if 
  $\mu_k(x,y) \wedge p(x) \wedge q(y)$ is consistent with $\CT$.   Define an equivalence relation on types, $d(p,q)<\infty$.  
 Then one can generalize (6) to (6'):     there is one universal model realizing each equivalence class of types.   

 \ssec{The {pp} topology on ${\U}$.}    
 \label{pptop}
  Let us topologize the sorts of ${\U}$ (specifically, ${\gg}^m$ for $m \in \Nn$).
 We give each sort of $\U$ the coarsest topology such that every pp-definable set, with parameters, is closed.  
 Thus as a pre-basis for  ${\gg}^m$  we take  the complements of sets of the form  $\{x:  R(x,c) \}$, with  $x=(x_1,\ldots,x_m)$, $c=(c_1,\ldots,c_k) \in {\U}$, and  $R$  a {pp} relation.   While the existential quantifiers in $R$ must be bounded in line
 with the local definition of pp, there is no requirement that $R$ imply a finite distance between its variables.  In case every pp relation
 is equivalent to a conjunction of atomic formulas (as will hold for the pattern language, \thmref{summary1}), we can of course
 restrict attention to atomic $R$.
As we included equality in the language, singletons in ${\gg}^n$ are  pp-definable with  parameters and hence closed, so 
each sort of $\U$ is T1.

We will now make a further   assumption; essentially, it asserts that $\infty$ can be separated from any given point
of $\U$ by disjoint open sets.

(${Loc}$)  There exists a pp formula $\mu_{>n}(x,y)$ such that $\mu_{>n} \perp \mu_n$, and for some 
$N  >n$, ${\gg} \subsetec \mu_{>n} \union \mu_{N} $.

 \begin{lem}\label{LC2}    Let $a \in {\gg}(\U)$, $n \in \Nn$, and  let $B=\mu_n(a) \subset {\gg}$ be a $\mu$-ball.   Then $B^m \subset {\gg}^m$ is compact.
 Assuming (${Loc}$) holds, ${\gg}$ is locally compact; \footnote{We use the term to mean that every point has a compact neighborhood; it does not presuppose that ${\gg}$ is Hausdorff.} and every compact subset of $\gg(\U)$ is contained
 in some ball $\mu_n(a)$.  
 \end{lem}
 
\prf Consider a family $F_i$ of basic closed subsets of $B^m$  with the finite intersection property.
 $F_i$ is defined  by $R_i(x,c_i)$ with $R_i$ a finite disjunction of {pp} formulas.  In an elementary extension $U'$ of ${\U}$, 
 one can find $d=(d_1,\ldots,d_m)$ with $R_i(d,c_i)$ holding for all $i$, as well as $\mu_n(a,d_i)$ for each $i \leq m$. 
 We can initially allow $U'$ to be a non-local extension,  then discard all  elements at infinite distance from $a$; this will not effect  the truth value of local formulas for the remaining elements, and the $d_i$ remain in the domain.  
 By \propref{prop2.1} there exists $r: U' \to {\U}$,
 $r|{\U}=Id_{\U}$.    Then $R_i(r(d),c_i)$ holds for each $i$, so $r(d) \in \meet_i F_i$ and
 $\meet_i F_i \neq \emptyset$.
 
 If (${Loc}$) holds, then $a$ lies in the interior of $\mu_N(a)$, since it is contained in the complement of $\mu_{>n}(a)$ which in turn is contained in $\mu_N(a)$.   Hence $a$ has a compact neighborhood, and so ${\gg}$ is locally compact.   If $Y$ is a compact
 subset of $\gg(\U)$, as it is covered by neighborhoods $\mu_N(a')$, it must be  covered by finitely many such neighborhoods,
 $\mu_N(a_1),\cdots, \mu_N(a_m)$;  
 and as the  `centers' $a_i$ are at finite distances from each other, it follows that $Y$ is contained in a large enough ball.
 \eprf 
     
\begin{lem} \label{LClem1}  Let $A \models \CT$, $a \in A$.   Assume:
\begin{enumerate}
\item  Each ball $\mu_n(a)$ is compact in the pp-topology.  
\item  For some $n$, if $\phi(x)$ is a pp-formula in one variable then either $\CT \models \neg \exists x \phi$, or $\phi$ is realized in $\mu_k(a)$.
\end{enumerate}  Then
$A$ is universal, i.e. every model of $\CT$ admits a homomorphism to $A$; and $\CT$ is pp-compact (ppC).  
 \end{lem}
 
 \prf     To prove universality it suffices to show that any $B \modec \CT$ admits a homomorphism into $A$.    Let $b_0 \in B$.
 Then any finite number of pp formulas true of $b_0$ are realized in $\mu_k(a)$; the sets of realizations are closed; hence by (1), they are realized
 simultaneously, i.e. some $a_0 \in \mu_k(a)$ realizes the pp-type of $b_0$.      Let   $f$ be a map from $C$ into $A$ preserving pp-formulas, with maximal
 possible domain $C$ containing $b_0 $.      To prove universality we have to show that  $C=B$.  Otherwise, let $b \in B \m C$.  Let $p$ be the  pp-type of $b$ over $C$.   
Then $p$ ensures $\mu_n(x,b_0)$ for some $n$.  By compactness of $\mu_{n+k}(a)$, as in the first step, $f(p)$ is realized in $A$.  This contradicts the maximality of $C$ and
proves that $C=B$, and so $A$ is universal.

   Let us now prove {ppC}.  Let $\phi(x,y)$ be a local pp formula contradicting equality,  as in the definition of {ppC}; $\phi(x,y)$ implies $\mu_m(x,y)$.    We may take $\phi$ symmetric
   (replacing  it by $\phi(x,y) \vee \phi(y,x)$ if necessary.)  
Suppose there is no finite bound   on the size $\phi$-cliques.   Then by compactness there exists a   weak model $B$of  $\CT$ and an infinite $\phi$-clique in $B$; it must
be contained in a $\mu_m$-ball, and so we can  take $B$ to be a local model.  Then $B$ embeds into $A$, so in fact $A$ contains an infinite $\phi$-clique within some ball
$\mu_n(a)$.  Extend it to a maximal $\phi$-clique $I \subset \mu_n(a)$.  Then $\mu_n(c)$ is covered by the open sets $\neg \phi(c,x)$, as $c$ ranges over $I$.  Hence $\mu_n(a)$
is covered by a finite  number of sets $\neg \phi(c,x)$, where now $c \in I_0$ for some finite $I_0 \subset I$.  If $d \in I \setminus I_0$ we have however $\phi(c,d)$ for each
$c \in I_0$.  This contradiction proves that a bound   exists.  
\eprf 
 
 It follows from universality that  $\CT$ has LJEP; of course, condition (2) really implies LJEP$_k$.

 \ssec{The pp topology on  $\Aut({\U})$.}
 Let $G= \Aut({\U})$.   Give $G$ the   topology whose pre-basic open set have the form
   \[\{g: \neg R (ga_1,\ldots,ga_n,b_1,\ldots,b_m) \}\] 
where $a_1,\ldots,a_n,b_1,\ldots,b_m \in \U$ and $R$ is pp.      The action $G \times \U \to \U$ is then continuous in the first variable; for fixed $g$ it is also continuous in the second variable, since $g$ is an automorphism.  
It is also easy to see that multiplication $G \times G \to G$ is continuous in each variable, and that inversion $g \mapsto g \inv$
is continuous.

As $\U$ is T1,  it follows from left-continuity of the action that $G$ is T1.

Let $G(a,n) = \{g \in G:  \mu_n(g(a),a)\}$.

\begin{lem}\label{LC3}    { $G(a,n)$     is compact.   Assuming (${Loc}$), $G$ is a locally compact space.}\end{lem}

\prf  Fix $a$ and $n$.  Let $u$ be an ultrafilter on a set $I$, and let $g_i \in G(a,n)$, $i \in I$;
  we need to find a limit point $g$ of $(g_i)_i$ along $u$.  Let $\U^*$ be the ultrapower of $\U$ along $u$,
  and let $g_*: \U \to \U^*$ be the ultraproduct of the maps $g_i: \U \to \U^*$.  Since $\mu_n(g_i(a),a)$, 
  we have $\mu_n(g_*(a),a)$ and hence $\mu_{n+k}(g_*(c),a)$ whenever $\mu_k(c,a)$.     Thus $g$ is defined on $\U$ into
  the local ultrapower.  Let $j$ be the diagonal embedding
  $\U \to \U^*$, ultrapower of $Id: \U \to \U$.  
      As $\U^* \models  \CT^{\pm}$,
    \propref{prop2.1} provides    a homomorphism $r: \U^* \to \U$ with $r \circ j = Id_{\U}$.  Let $g=r \circ g_* $.
 Then $g \in End(\U)=\Aut(\U)$.  If $R(g_ia_1,\ldots,g_ia_n,b_1,\ldots,b_m)$ holds for $u$-almost all $i \in I$,
 then $\U^* \models R(g_*a_1,\ldots,g_*a_n,jb_1,\ldots,jb_m)$ so $\U \models R(ga_1,\ldots,ga_n,b_1,\ldots,b_m)$.
Hence   $g$ is indeed a limit point of $(g_i)_i$ along $u$. 

  Assuming (${Loc}$), $G$ is locally compact:  since translations are continuous, it suffices to show
  that the identity element $1 \in G$ has a compact neighborhood.  Indeed  as in \lemref{LC2}, the compact set $G(a,N)$ contains 
  an open set including $1$, namely the complement of $\{g: \mu_{>n}(a,g(a))\}$. 
\eprf

Let $X$ be a $G$-invariant subset of ${\gg}$.   Define the infinitesimal elements of $G$, in its action on $X$,   to be 
\[  {\fg_X}= \{g \in G: \hbox{for every nonempty open } U \subset X, \    \       gU \meet U \neq \emptyset  \}  \]

Fix $X$ to be either ${\gg}(\U)$ or a complete 0-definable pp type $P \subset {\gg}(\U)$; and 
  let $\fg=\fg_X, \CG = G/\fg, \pi: G \to \CG$ the quotient map.

  By Lemma C.1 of \cite{patterns}, $\fg$ is a closed normal subgroup of $G$, and $\CG = G/\fg$ is Hausdorff.  Note that $G/\fg$
  remains locally compact, with continuous inversion and right and left multiplication.  
It follows from Ellis' joint continuity theorem \cite{ellis} 
that $G/\fg$ is a locally compact Hausdorff topological group, acting continuously 
 on $\U/\fg$.    (A  direct proof in our setting is given in  \cite{patterns} 3.27.)   

We have shown the main part of:  

\begin{prop} \label{2.5}   Assume (ecB) and (${Loc}$).     Let $\U$ be the universal e.c. model,  $G=Aut(\U)$,$\fg=\fg_{\U}$, $\CG=G/\fg$.   
 Then  $\CG$   is a locally compact topological group.      With $\pi: G \to \CG$ the quotient map, we have:

\begin{itemize}
\item  
If $Y \subset \CG$ is precompact, then  $\pi \inv(Y)$ is bounded.     
\item 
If $W \subset G$ is bounded, then $\pi W$ is precompact.  \end{itemize}
 \end{prop}
 Here $W \subset G$ is said to be bounded if    for some/any $a \in \gg(\U)$, for some $N \in \Nn$,
 $Wa \subset \mu_N(a)$; see \defref{groupdefs}  \ref{bounded}.  

\prf    
 We saw that $\CG$ is a locally compact topological group, and that  $\fg$ is a closed normal subgroup of $G$.   
 Fix $a \in \gg(\U)$.   
By \lemref{LC2},  $a$ has  a compact   neighborhood  $W=\mu_n(a)$.  Then by definition of $\fg$,
if $g \in \fg$ then $gW \meet W \neq \emptyset $;       so $\mu_n(a,gb)$ holds for some $b$ with $\mu_n(a,b)$;
hence $\mu_n(ga,gb)$, and $\mu_{2n}(a,ga)$.  This shows that $\fg \subset G(a,2n)$, which is compact by \lemref{LC3}.
Thus $\fg$ is compact.  
It follows easily\footnote{[Let $\pi: G \to \CG$ be the quotient map.  Let $Y \subset G$ be closed with $\pi(Y)$ compact.  Then $Y$ is compact:
Let $y_i \in Y, i \in I$, and $u$ an ultrafilter on $I$.  The $\lim_{i \to u} \pi(y_i) = \pi(b) $ for some $b$.
Since $b\fg$ is compact, it has a compact neighborhood $U=  U \fg$; and we have $y_i \in U$ for almost all $i$;
so $y_i$ approaches a limit along $u$.]} that if $Y \subset \CG$ is compact, then so is the pullback to $G$.   Thus
for the first item 
 it suffices to show that if $Y \subset G$ is compact, then $Y \subset G(a,N)$ for some $N$.  Indeed as $Y$ is compact,
$Y a$ is compact.   By \lemref{LC2}, 
 $Y a \subset \mu_N(a)$ for some $N$;
hence by definition $Y \subset G(a,N)$.

  If $W \subset G$ is bounded, $W \subset G(a,N)$ for some $N$;  we saw that $G(a,N)$ is compact; hence so is $\pi(G(a,N))$.
So  $\pi(W) \subset \pi(G(a,N))$ is precompact.
  
 \eprf

 \begin{prop} \label{autlc}
   Assume (ecB) and (${Loc}$); also that ${\L}$ is countable.  Let $P$ be a  complete pp type  in ${\gg}(\U)$, 
    $\CG(P)=G/\fg_{P}$.   Then $\CG(P) $ 
 is a second countable, locally compact topological group.   If $P'$ is another complete pp type  in ${\gg}(\U)$, 
 then $\CG(P'),\CG(P)$ differ only by a compact isogeny.
 \end{prop}
  
 \prf  We have the natural homomorphism $G \to \CG(P)$, by restriction and then going to the quotient.  It is clearly continuous, when 
 $\CG(P)$ is endowed with the topology described above, with respect to the action of $G$ on $P$.
  The kernel is compact, indeed with $n$ as in the proof of \propref{2.5},  for any $a \in P$, the kernel   is contained
 in $G(a,n)$.    It follows that $\CG(P)$ is locally compact, and (regardless of $P$) differs
 by a compact isogeny from $\CG$.  When ${\L}$  is countable, the proof of \cite{patterns} 3.24  shows that $\CG(P)$ is second countable.  (using $\si$-compactness
 via  $G = \union_{n \in \Nn} G(a,n)$   in place of compactness; note that Lemma C.4 remains valid when $\ft_1$ is 
 $\si$-compact, and $X \times Y$ is $\ft_2$-compact for $\ft_1$-compact $X,Y$, as it the case here.)  
  \eprf

   \pb
  
 \section{The core space}

Let $T$ be a complete first-order theory, in a language $L$ with a distinguished sort ${D}$ and locality relation on ${D}$, that we take to be generated by a single reflexive, symmetric binary relation ${m_1}$ on ${D}$.    
 Let ${d_m}(x,y)$ be the Cayley-graph distance formed from ${m_1}$, as in the discussion of local logic.  A model of $T$ will be said to be {\em ${m}$-local}
 if any two elements of ${D}$ are at finite ${d_m}$-  distance.      
  We treat  $T$ in local logic.  
   This   means, first of all, that we 
  restrict attention to ${m}$-local models $M$ of $T$,  and in particular  to the part of the  type space consisting of types realized in elementary extensions that are themesleves ${m}$-local.    If a weak model $N'$ is an elementary extension of $M$ 
  in the usual sense, and $N$ is obtained by simply removing any points of ${D}(N')$ at infinite distance from points of ${D}(M)$,
  then by fiat, $N$ is a (local) elementary extension of $T$ in the local logic sense.

We will use only local formulas, i.e. formulas built up using  local quantifers only; with the one exception that we  assume $T$ is complete
 with respect to sentences $(\forall x) \phi$ and $(\exists x)\phi$, where $\phi$ is local.  
\footnote{ We can assume or arrange by Morleyzation that $T$ admits local  quantifier-elimination. i.e. every local formula
 is equivalent to a quantifier-free one.    In this case $T$ is determined by its universal part; the models of $T$ are e.c. models of $T_{\forall}$.     In \cite{patterns} it was useful to allow more general universal theories, restricting attention to e.c. models; this could be done here too.}.      
  Note that the 'local' restriction on   types $p(x) \in S_{D}(M)$ amounts simply to saying that   some formula $m_n(x,b)$ be represented in $p$.      Following convention in enriched logics, we denote the usual quantifier-free type
     space by $S^{weak} (M)$; and the subspace of local types by $S(M)$.

 For simplicity we
     will concentrate on the space of local $1$-types in ${D}$, $S_{D}$.   
      $S_{{D}}(M)$ is an open subset of $S_{D}^{weak}(M)$
     with the usual topology;   hence it is  locally compact.  
            
\ssec{The pattern language}   Let us recall the construction in \cite{patterns} of the  pattern language ${\L}$, and 
 the  ${\L}$-structure on the type spaces
  $S_{x}(M)$, $M \models T$.   Here we 
  concentrate on a single sort ${D}$ of $T$,  and will be interested only in one sort  of ${\L}$, to be denoted $\gg$,  corresponding to  local  $1$-types on ${D}$.  Let $x,x_i$ be variables 
  ranging over ${D}$.  
Let $y$ be an additional tuple of variables (the {\em parameter variables}.)  
  Let $t=(\phi_1,\ldots,\phi_n;\alpha)$ be
an $n$-tuple of formulas $\phi_i(x_i,y)$, and let  $\alpha(y)$ be a formula.
 To each such $t$ we associate a   relation symbol $\R_t$ of ${{\L}}$, taking variables $(\xi_1,\ldots,\xi_n)$ ranging over $\gg$.
For any $M \models T$, 
 we define an ${{\L}}$-structure with universe $S=S_D(M)$.   
The interpretation of $\R_t$ is:  
 \[\R_{t} ^S =\{(p_1,\ldots,p_n) \in S:   (\forall a \in \alpha(M)) \bigvee_{i \leq n} (\phi_i(x,a) \in p_i) \} \]

Thus when $n=1$, $\R_{\phi;\alpha}$ indicates a partial definition scheme for $\phi$, thus the notation $\R$.
(In \cite{patterns}, $Om_{\neg \phi; \alpha}$ was used instead, where $Om$ stands for {\em omitting} $\neg \phi$.
These are equivalent in the presence of negation  in $L$, but even then $\R$ is more natural.)
 For typographical reasons, we will sometimes write $\R[t]$ for $\R_t$.
  
   \ssec{Locality structure on $\gg$}

  \begin{defn}  A model $M$ of $T$ is {\em orbit-bounded} (in the sort ${D}$) if for some $k \in \Nn$, 
  for all $a,b \in {{D}}(M)$, there exists $\alpha \in Aut(M)$ such that $a,\alpha(b)$ are at distance $\leq k$.  \end{defn}

We will make two assumptions about $T$ as a local theory.
 
 (OB)   $T$ has an orbit-bounded model $M$.
 
 (DP)  For each $r \in \Nn$, for some $m=m(r) \in \Nn$, every  $r+1$-ball is a   union  of at most $m$ $r$-balls.   
 
  \medskip
 
 The second assumption, (DP), is the {\em doubling} property for the metric $d_m$.     Assuming 
 any two points are connected by a finite chain of points at distance $1$ (i.e. the locality is generated by a single relation),
   (DP) follows from the same statement for $r=1$.
 If it holds for one ball    in an orbit-bounded model, then 
 the uniform version will be valid too. 
 
    (DP) is needed in order  to define the locality structure for the pattern language,   by a pp formula.
    It follows from doubling that every $n$-ball is a union of finitely many $1$-balls; hence
 a type $p$ of ${D}$ over a model $M$ is local iff it represents the formula $m(x,y)$, i.e. $m(x,c) \in p$ for some $c \in {D}(M)$.
( In other words we will not have $ \R[{\neg m(x,y);D(y) }]$.)
  
  The terminology orbit-bounded is short for `the space of orbits of $Aut(M)$ on $M$ is bounded, under the induced metric'.
In our applications, (OB) will hold in every local model of $T$.  In general, if $T$ has one orbit-bounded model, then
every $\aleph_0$-homogneous model will be orbit-bounded.
   
   \begin{defn} \label{mu} We define the locality generating relation $\mu_1(\xi_1,\xi_2)$ for ${\L}$ to be 
   \[\R[{\neg m_1(x_1,y_1), \neg m_1(x_2,y_2); \neg m_5(y_1,y_2)}]\]  
   
   Towards (${Loc}$), we also define $\mu_{>1}(\xi_1,\xi_2)$ by the formula:
   
 \[  \R[{\neg m_1(x_1,y_1) , \neg m_1(x_2,y_2) ; m_{5}(y_1,y_2)    }  ]    \] 
   \end{defn}

 (We used $5$ above, where $3$ would do, in order to ensure below that close Lascar neighbors have $\mu$-distance at most $1$;
 it is not an essential point.)
   
For the space of types over a model $M$, if $p_1,p_2 \in S(M)$ are local types, and $S(M) \models \mu_1(p_1,p_2)$,
by (DP) there exist $b_1,b_2 \in D(M)$
with $m_1(x_i,b_i) \in p_i$; now $\mu_1(p_1,p_2)$ guarantees that $m_5(b_1,b_2)$ holds, so that any two realizations
of $p_1,p_2$ will be at $m$-distance at most $7$.  

Conversely if realizations $b_i$ of $p_i$ exist with $m_1(b_1,b_2)$, 
then   $\mu_1(p_1,p_2)$ holds (even with $3$ in place of $5$): we have to show that if $M \models \neg m_5(c_1,c_2)$
then $\neg m_1(x_1,c_1) \in p_1$ or $\neg m_2(x_1,c_2) \in p_2$.  Otherwise, $m_1(x_1,c_i) \in p_i$ so $m_1(b_i,c_i)$.
Then (initially in an elementary extension of $M$) the chain $c_1,b_2,b_2,c_2$ shows that $m_3(c_1,c_2)$ holds,
a contradiction.  

\begin{lem}\label{OB2}   Assume $M$ is orbit-bounded; then 
    $S_{D}(M)$ is  orbit-bounded too,  even just with respect to the action of $Aut(M)$.           \end{lem}
\prf  Let $p, q \in S_{{D}}(M)$.    
Since $p(x)$ is local,   for some $n \in \Nn$, $p(x)$ contains a formula defining an $m$-ball of radius $n$.    By the doubling property (DP) it follows that $p(x)$ contains a formula defining an $m$-ball of radius $1$;  i.e.
  $m_1(x,c) \in p$ for some $c \in M$.  Similarly $m_1(x,d) \in q$ for some $d \in M$.  Let $\si \in Aut(M)$ be such that
$m(d,\si(c)) \leq k$.   Then $\si$ induces an automorphism of $S(M)$,  and $m_1(x,\si(c)) \in \si(p)$. 
 Let $d=d_1,\cdots,d_k=\si(c)$ be a chain of elements of $D(M)$, with $m_1(d_i,d_{i+1})$.   (Say $k$ is even.)
Also write $d_i$ for the algebraic type $x=d_i$.  We have $\mu_1(q,d_2), \mu_1(d_2,d_4), \cdots, \mu_1(d_{k-2},d_k),
\mu_1(d_k,\si(p))$.   Thus $\mu_{k/2+1} (\si(p),q)$ holds in $S(M)$.    \eprf

    \ssec{The  pattern theory  of $T$}

 We define ${\CT}$ be the set of  local, primitive universal ${\L}$-sentences true in   $S(M)$ for some $M \models T$.  
 Defined this way, it is not clear that ${\CT}$ is consistent; but this,  along with  LJEP,  follows from \lemref{orbitbounded2}.

 \begin{lem} \label{orbitbounded2}   Let $M$ be an orbit-bounded model of $T$.  Then $S(M) \models {\CT}$.  \end{lem}

 \prf   
 Let $N \models T$.  If $S(M) \models (\exists x_1)\cdots(\exists x_n) \Theta$, where $\Theta$  is a finite conjunction of basic relations of ${\L}$, we  have to show that $S(N) \models (\exists x_1)\cdots(\exists x_n) \Theta$.   For simplicity, let us assume
 $\Theta$ is a basic relation $\R[{\phi_1,\ldots,\phi_n;\alpha}]$.   
  
  Let $p_1,\ldots,p_n \in S(M)$  be such that $S(M) \models \Theta(p_1,\ldots,p_n)$.  
   Since $M$ and $S(M)$ are local, each $p_i$ is at some finite $m$-distance $d$
 from some $a \in M$, i.e. $m_d(a,x_i) \in p_i$.   Thus the following partial type over $M$ is consistent:
 \[ \{  \bigvee_{i=1}^n \phi_i(x_i,c):   c \in \alpha(M)  \} \union \{ m_d(a,x_i): i =1,\ldots,n \} \]

 Hence,  for any $\si \in Aut(M)$,  the same partial type is consistent if we replace $a$ by $\si(a)$.   By the orbit boundedness
 of $M$, it follows that for {\em any} $a' \in D(M)$, the type below is consistent too:
  \[ \{   \bigvee_{i=1}^n \phi_i(x_i,c):   c \in \alpha(M)  \union \{ m_{d+k}(a',x_i): i =1,\ldots,n \} \]
  Thus for any $m$, we have in $M$:
\[ (\forall u  \in {D})(\forall c_1,\ldots,c_m)( \bigwedge_j\alpha(c_j) \to   (\exists x_1,\ldots,x_n) \bigwedge_i m_{d+k}(u,x_i) \wedge
 \bigwedge_j (  \bigvee_{i=1}^n \phi_i(x_i,c_j))  )  \] %
  This applies equally to $N$ (using the completeness of $T$), and so retracing out steps, we find that for any $b' \in N$ there exist 
  types $q_1,\ldots,q_n$ with $\mu_{d+k}(b',x_i) \in q_i$ and $S(N) \models \R_t(q_1,\ldots,q_n)$.    The $q_i$ are local types,
  and show that $S(N) \models (\exists x_1)\cdots(\exists x_n) \Theta$, as promised.

  \eprf

 Regarding (${Loc}$),  clearly $\mu_1,\mu_{>1}$ are incompatible in any local type space. 
 On the other hand, if $p_1,p_2 \in S(M)$, then either $\mu_{>1}$ holds or $m_i(x_i,b_i) \in p_i$ for some
 $b_1,b_2 \in \gg(M)$ with $m_5(b_1,b_2)$.    Then as in the proof of \lemref{OB2} we have $\mu_4(p_1,p_2)$.
 Thus $S(M) \models (\forall \xi)(\mu_{>1}(\xi) \vee \mu_{4}(\xi))$, so (say by \lemref{embedinmodel}) 
 ${\CT} \modec  (\forall \xi)(\mu_{>1}(\xi) \vee \mu_{4}(\xi))$.   

  \begin{lem}\label{embedinmodel} Let  $M \models T$ be orbit-bounded.  $A \models {\CT}$.  Then there exists a homomorphism $f: A \to S(M)$.   In particular any
 e.c. model of ${\CT}$ is isomorphic to a substructure of $S(M)$.   
  ${\CT}$ has the properties   (ecB)   and (${Loc}$).
  \end{lem}
\prf  We have to check that 
 the usual (locally compact) type space topology of $S_D(M)$ satisfies the hyptheses of \lemref{LClem1}.  (2) there was already
 proved in \lemref{OB2}.

  The interpretation of $R=\R[\phi_1,\ldots,\phi_n; \alpha]$ is closed, since if $p=(p_1,\ldots,p_n) \notin R$,
then $\neg \phi_i(x_i,c) \in p_i$ for some $c \in \alpha(M)$; the set of all tuples $q=(q_1,\ldots,q_n)$  with $\neg \phi_i(x_i,c) \in q_i$
is an open neighborhood of $p$ in the product topology, disjoint from $R$.  

As for compactness of $\mu_n(p)$, where $p \in S_D(M)$:   if $q \in \mu_n(p)$, there exists a chain
$p=p_0,\cdots, p_n=q$ with $\mu_1(p_i,p_{i+1})$.  Pick $a_i \in D(M)$ with $m_1(x,a_i) \in p_i$.  Then
$m_5(a_i,a_{i+1})$ , so $m_{1+5n}(x,a_0) \in q$.  The types $q$ with  $m_{1+5n}(x,a_0) \in q$ form 
a compact set  $W$ in the usual topology on the weak type space, which is entirely contained in the space of local types;
hence $W$ is compact.  Now $\mu_n(p)$ is a subset of $W$  defined by a pp relation, hence is closed in $W$ and thus
  it is compact too.   

 \eprf
 
\propref{prop2.1} and \propref{2.5} are thus valid for ${\CT}$.   We call $\U$ the core space of $T$.   To summarize:

 Let $T$ be a local theory in the sort $D$, complete  (with respect to local formulas),
 and satisfying (DP) and (OB).
We have defined the {\em pattern language ${\L}$},
 a local primitive universal theory in the sort $\gg $, valid in $S_D(M)$
when $M \models T$.  

\begin{defn}
A {\em core} for (type spaces of models of) $T$ is an  ${\L}$-structure $J$ such that:
\begin{itemize}
\item For any orbit-bounded $M \models T$, there exists an ${\L}$- embedding $j: J \to S(M)$.
\item  For any such $j$, there exists a homomorphism $r: S(M) \to J$
such that $r \circ j = Id_J$.
\end{itemize}
\end{defn}

\begin{thm}   \label{summary1}     Let $T$ be a complete local theory in the sort $D$.
  Assume $T$ has an orbit-bounded model, where every $2$-ball is contained in finitely many $1$-balls.  Then:
\begin{enumerate}
\item  A core   exists and is unique up to isomorphism.  We denote it $\CU=\co(T)$.
 \item Reflection:   let $\Theta$ be a 
   universal sentence $(\forall p_1,\ldots,p_k) R$, where $R$ is an arbitrary Boolean combination of local pp  formulas in the pattern
 language $\mathcal{L}$ .   Let $M$ be an orbit-bounded model of $T$.
 If $S(M) \models \Theta$ then   $\CU \models \Theta$.
\item  Quantifier elimination:  
$\mathcal{L}$-atomic   types in $\CU$ are orbits of $Aut(\CU)$.   Every pp formula is equivalent 
(in any $S(M)$ as well as in $\CU$) 
to a certain conjunction of atomic formulas.

\item  $Aut(\CU)$   has a natural  locally compact topology:    \[ W(R;a,b) :=   \{g:  \R_t (ga_1,\ldots,ga_n,b_1,\ldots,b_m) \}\] 
 is a basic closed set.   The canonical Hausdorff quotient  $\CG(T) : = Aut(\CU)/\fg
 $ is a locally compact topological group.
 
 \end{enumerate}
 \end{thm}

We have just proved (1) and (4).   By definition of a core, there exists an ${\L}$- embedding $j: \CU \to S(M)$,
as well as  $r: S(M) \to \CU$ with $r \circ j = Id_{\CU}$. It follows that the truth value of any pp $\psi$ formula is preserved by $j$;
i.e.  $\CU \models \psi(a)$ if and only if $S(M) \models \psi(j(a))$.  Hence for an arbitrary Boolean combinations $\Psi$ of pp formulas,
if  $S(M) \models (\forall x)\Psi$,  then $S(M) \models \Psi(j(a))$ so $\CU \models \Psi(a)$, this for any $a$ from $\CU$.   
This proves (2).  For (3)   see \cite{patterns},  3.12 (1).

\ssec{The local Lascar-Shelah relation} \label{LSh}

We continue to assume $T$ has an orbit-bounded model $M$, with bound $k$,  and that doubling (DP) holds; and we concentrate
on a single sort $D$.

A {\em definable family of finite local partitions} is a pair of formulas $(B,E)$, with $B$ a   non-empty 0-definable subset of $D^n$ for some $n \geq 1$, and 
$E(x,y; z_1,\ldots,z_n)  \subset D^2 \times B$,  such that for some $f \in \Nn$, with $f$ strictly above the radius $k$ 
of $M/Aut(M)$  we have:
\begin{enumerate}
\item  $E$ implies $m_{f}(x,z_1)$ and $m_{f}(y,z_1)$. 
\item   
 $B$ implies $m_f(z_1,z_i)$ for each $i \leq n$.
\item  For any $b=(b_1,\ldots,b_n) \in B$, $E_b(x,y):= E(x,y;b)$ defines an equivalence relation on $m_{f}(b_1)$, with at most $f$ classes.
\end{enumerate}

 Two elements $a,a' \in D$ are  {\em  Lascar-Shelah neighbors}, denoted $a {\Sh} a'$, if for each $n \in \Nn$, every non-empty 0-definable set $B \subset D^n$, and  every definable family $E \subset D^2 \times B$ of finite local partitions, 
 \[M \models (\exists z_1,\cdots,z_n)  E(a,a',z_1,\ldots,z_n) .\] 
 If this holds for a single $E$, we will say that they are  Lascar-Shelah neighbors with respect to $E$.
 The notation is intended to suggest an isosceles triangle, with $a$ and $b$ at the distinct endpoints of the equal sides.
 We will nevertheless sometimes write $\Sh(x,y)$ for $x \Sh y$.

   We call two (local) types $p_1,p_2$ over a model $N$ {\em Lascar neighbors} if  
there exist $a_i \models p_i$ such that $a_1,a_2$ are Lascar neighbors, i.e. if $p_1(x_1) \union p_2(x_2) \union {\Sh}(x_1,x_2)$
is consistent.   

Call $p_1,p_2$  {\em close neighbors} if  for every definable family of finite local partitions
 $(E_b)_{b \in B} $, for some $b \in B(N)$ and $d' \in D(N)$, $x_i E_b d' \in p_i$.  This implies that $p_1 \union p_2 \models {\Sh}(x_1,x_2)$, and
in particular $p_1,p_2$ are neighbors.    

 If $p_i,p_{i+1}$  are Lascar neighbors for $i=1,\ldots,n-1$, we will say that $p_1,p_n$ are at Lascar distance at most $n$, and write $\Sh^n(p_1,p_n)$.


 \begin{lem} \label{shrem} Let  $M$ be an orbit-bounded model of $T$,  $a,a' \in D(M)$.  Then $a \Sh a'$ iff for some   elementary extension $M'$ of $M$,
 and elementary submodel $N \prec M'$, $N \cong M$, we have $tp(a/N)=tp(a'/N)$.   (Here by convention $M'$ is local, and $tp$ refers
 to the quantifier-free type.)   
 
 Also, $a{\Sh} a'$ implies $m( x,y) \leq 2$.
 \end{lem}
 
 \prf  Assume $tp(a/N)=tp(a'/N)$.  We have $D(N) \neq \emptyset$; pick $d \in D(N)$.   Since $M'$ is local, $M' \models m_l(a,d), m_l(a',d)$  for some $l$.  By (DP),
 $m_l(d)$ is a finite union of balls $m_1(c_i)$, $c_1,\ldots,c_r \in D(N)$.  So $m_1(a,c)$ holds for some $c=c_j \in D(N)$.
 Since   $tp(a/N)=tp(a'/N)$, we have $m_1(a',c)$ too.  This proves that $m( x,y) \leq 2$.
 
 Now let $(B,E)$ be a definable family of  finite local partitions.   
 
 As $B \neq \emptyset$, and using  (OB), 
 there exists $b=(b_1,\ldots,b_n) \in B(N)$ with $m(c,b_1) \leq k$.   Thus $m_{k+1} (a,b_1) $ and $m_{k+1}(a',b_1)$
 hold.  Since $E_b$ defines a finite partition of $m_{f}(b_1)$ and $f \geq k+1$, and using  that $N$ is a model and 
 $tp(a/N)=tp(a'/N)$, we have $E(a,a',b')$.   This proves that $a \Sh a'$.  
 
 Along the way we saw that $m(a,a') \leq 2$.  
 
 Conversely, assume $a \Sh a'$.  Let $y_c: c \in M$ be variables, and let $\Delta$ be the diagram of $M$,
 i.e. the set of $L$-local  formulas in the variables $y_c$ that become true in $M$ under the assignment $y_c \mapsto c$.  
 Fix one of the variables $y_c$, and call it $y_1$. 
 Let $\Delta'$ be the set of all local formulas $\phi(x,y) \iff \phi(x',y)$, where the $y$ are some of the $y_c$, including $y_1$, and $\phi$ is a  formula of $L$.
 Since $a  \Sh a'$, it is easy to see that $T \union \Delta \union \Delta' \union \{m_{k+1}(a,y_1) ,m_{k+1}(a',y_1)\}$ is consistent.  Hence this set of formulas is realized
 in some weak model, and restricting to    elements at finite distance from $a$ and $a'$, it is realized in some model.
 This leads to $N \cong M$ with $N,M \prec M'$ for some $M' \succ M$.
 \eprf
  
\ssec{The neighbor relation in the pattern language}  \label{neighbor2}

 \begin{lem}   \label{LSdef}   
 \begin{enumerate}  
 \item (${\L}$-$\bigwedge$-definability of ${\Sh}$.)  There exist atomic formulas $\psi_m$ of ${\L}$ such that for any $N \models T$, and any $p,q \in S(N)$, $p{\Sh} q $ iff for each $m \in \Nn$,  $S(N) \models \psi_m(p,q)$.  
 \item  The close neighbor relation on $S(M)$  is an intersection of $\leq |L|$ open sets in the pp topology.   
\item Let $\rho: S(N) \to S(N)$ be a retraction, i.e. a homomorphism with $\rho^2=\rho$.   Let $p \in S(N)$.
 Then $p,\rho(p)$ are close Lascar neighbors. 
\item  Assume $T$ is orbit-bounded.   Let $J$ be a copy of $\CU$ in $S(M)$, and assume $p,q \in J$ cannot be separated by disjoint open subsets of $J$.
 Then $p,q$ are at Lascar distance at most $2$.
 \end{enumerate}
  \end{lem}
 
 \prf 
 

(1)   Let $p,q \in S(N)$.   If $p {\Sh} q$ then for every  definable  family $(B,E)$ of  finite local partitions such that
 \[T \models \neg (\exists x,x',u,y) (\phi(x,u) \wedge \phi'(x',u) \wedge \alpha(u) \wedge B(y) \wedge  E(x,x',y)) \]
we have \[S(N) \models \R [\neg \phi, \neg \phi'; \alpha](p,q)\]
Conversely,  if $(p,q) \notin \Sh$, then $p(x) \union q(x')$ (a set of formulas of $L(M)$, including the diagram of $M$) must be inconsistent with the following set of formulas:
\[   \{(\exists y) (B(y) \wedge  E(x,x',y) ):   (B,E) \} \]
where $(B,E)$ ranges over all definable families of finite local partitions, $y=(y_1,\ldots,y_k)$.   Note that these formulas imply a bounded distance
between $x,x'$ and each $y_i$, and between them and the elements of $M$ (the latter because $p,q$ are local types.)
By compactness, for some tuple $c$ from $N$, $\phi(x,c) \in p, \phi(x',c) \in p'$, $\a \in tp_{N}(c)$,   and for some finite set $(B_j,E_j)$ of definable families of finite local partitions,
\[T \models \neg (\exists x,x',u,y_1,\ldots,y_r) ( \phi(x,u) \wedge \phi(x',u) \wedge \a(u) \wedge \bigwedge_j B_j(y_j) \wedge E_j (x,x',y_j) )\] 
Let $y=(y_1,\ldots,y_r)$, $B(y) = \bigwedge_j B_j(y_j)$, $E(x,x',y) = \bigwedge_j E_j(x,x',y_j)$.
Then 
\[T \models \neg (\exists x,x',u,y) ( \phi(x,u) \wedge \phi(x',u) \wedge \a(u) \wedge  \bigwedge B(y) \wedge E (x,x',y) )\] 
yet $c$ shows that $S(N)$ does not omit $\phi,\phi', \alpha$.  Thus
 $\R [\neg \phi, \neg\phi'; \alpha](p,q)$ holds for each $\phi,\phi',\alpha$ as above iff $p{\Sh}q$.

 (2)  Let $q=\rho(p)$.  Let $\bar{E}=(E_b)_{b \in B} $ be a    definable  family of local finite partitions.
 Consider the atomic formula $\psi = \psi_{\bar{E}}$, 
 \[\psi(\xi_1,\xi_2) = \R[\neg x_1E_y u, \neg x_2 E_y u; B(y) \wedge D(u) ] \]
   It is clear that   $p,q \in S(M)$ are close Lascar neighbors with respect  to $\bar{E}$
   iff $S(M) \models \neg \psi(p,q)$.   And $p,q \in S(M)$   close Lascar neighbors iff they are 
   close Lascar neighbors with respect  to each $\bar{E}$.
 
 (3) Since a realization of $q$
 must be $E_b$-equivalent to some $d' \in N$ (as there are only finitely many classes, all represented in $N$),
$S(N) \models \neg \psi(q,q)$.  Now if we had  $S(N) \models   \psi(p,q)$, applying the  
 homomorphism $\rho$ would give $\psi(q,q)$, which we have just ruled out; so $S(N) \models \neg \psi(p,q)$.
 This means that  for some $b \in B(N)$ and $d' \in D(N)$, $x_1 E_b d' \in p$ and $x_2 E_b d' \in q$.  As this holds for all definable families of finite partitions, $p,\rho(p)$ are close neighbors. 
 
 (4) By \lemref{embedinmodel},
 $\CT$ has Loc.  
  Let $N$ be a bounded open set containing both $p$ and $q$.   
   By compactness, it suffices to show that for each  definable  family $\bar{E}$ of local finite partitions,
  there exists $r \in cl(N) $ that is a Lascar neighbor to $p$ and to $q$ with respect to $\bar{E}$.  In fact we will show
  that $p,q$ have a common close Lascar neighbor in $N$, with respect to $\bar{E}$.  By (2), the close neighbor
  relation with respect to $\bar{E}$  is pp-open.  Let $N'$ the $\bar{E}$-close neighbors of  $p$ within $N$;
  then $N'$ is open and $p \in N'$.  Similarly let $N''$ be the set of $\bar{E}$-close neighbors of $p$ within $N$, an open neighborhood of $q$.  Since $p,q$ cannot be separated we have $N' \meet N'' \neq \emptyset$.  Any 
  $r \in N' \meet N'' $ is a common  close Lascar neighbor of $p,q$ in $N$.
 
 \eprf

    Note that the $\psi_m$ define the same relation on $S(N)$ and on any retract, notably any copy $J$ of $\CU$.
 We will refer to the relation they define on any model of $\CT$ as the Lascar neighbor relation, and denote it by $\Sh$.
 Later we will use the same symbol on $\CU$ and on $Aut(\CU)$ and $Aut(\CU)/\fg$;  when the need arises to clarify, we will use a supersript, such as $\Sh^{\CU}$.

\ssec{Automorphisms of models and pattern automorphisms} \label{LSh2}

In this section,  $M,N$ will denote orbit-bounded (local) models of $T$.  
We say that a subset $W$ of $Aut(M)$ is {\em bounded} if for some $a \in D(M)$, $Wa :=\{w(a): w \in W \}$ is
contained in a $d_m$-ball of finite radius.  Equivalently, for {\em all} $a \in D(M)$, $Wa$ has finite radius. (\defref{groupdefs} \ref{bounded}.)

By \lemref{LSdef} (1), $\Sh$ is defined by a conjunction $\Psi$ of atomic sentences in any $S(N)$, $N$ an orbit-bounded model of $T$.
We define $\Sh$ by the same formula on $\CU$.   Since  $\CU$ embeds into $S(N)$, this does not depend on the choice of $\Psi$.

On $G=Aut(\CU)$, we let  
\[\Sh= \Sh^{Aut(\CU)}  = \{g \in Aut(\CU):  (\forall a \in \CU)(a,g(a)) \in \Sh^{\CU} \} \]
Let $\Sh^{\CG} $ be the image of  $ \Sh^{Aut(\CU)}$ in $\CG$.

$\Sh=\Sh^{\CG}$ is a normal\footnote{I.e. conjugation-invariant.  See \defref{groupdefs} for group-theoretic terminology.},
symmetric, compact subset of $\CG$.  This follows from the same assertion in $Aut(\CU)$.  In $Aut(\CU)$ we saw that 
the $2$-ball is compact and contains $\Sh$ as a closed subset.

 \begin{thm}  \label{qh-basic}  Assume (DP).  Let $M$ be an orbit-bounded model of $T$.
  $\iota: \CU \to S(M)$ an embedding,
 and $\rho: S(M) \to \CU$ a retraction.  
 Then there exists a canonical  quasi-homomorphism ${\varphi}:Aut(M) \to (\CG(T): {\Sh})$.   We have:
 \begin{enumerate}
 \item  When $C \subset \CG$ is precompact, ${\varphi} \inv(C)$ is bounded in $Aut(M)$.
 \item  When $W \subset Aut(M)$ is bounded, ${\varphi}(W)$ is precompact in $\CG$.
 \end{enumerate}
 \end{thm}

  \prf   
  Let $J=\iota(\CU)$,  and identify $J$ with $\CU$; so that $\iota$ is the inclusion
  map $J \to S(M)$, and $\rho|J = Id_J$.    
     The definitions of ${\Sh}$ on $J$ and on $S(M)$ are compatible, via $\iota$ and $\rho$.
 
 For $g \in Aut(M)$, define
  \[  {\phi}(g) = \rho \circ g \circ \iota \in Aut(\CU)   \]
  
 And let ${\varphi} = \pi \circ \phi$, where $\pi:  Aut(\CU) \to \CG$ is the quotient map.    Since $\pi$
 is a continuous and proper homomorphism, it suffices to prove the quasi-homomorphism property and (1,2) for $\phi$.   
 
  We have  ${\phi}(Id_M)=Id_J$ since $\rho \circ \iota = Id_J$.     We have to show, for $g,h \in Aut(M)$, 
  that ${\phi}(gh) \in {\phi}(g) {\phi}(h) {\Sh}$.  For any $q \in S(M)$, we have $  \rho(q) {\Sh} q$ in $S(M)$ by \lemref{LSdef}(3).

  Thus 
  \[  \rho h(q)   {\Sh} h(q)   \] 
  As $ \rho g: S(M) \to \CU$ is an ${\L}$-homomorphism,
 $\rho g  \rho h(q)  {\Sh} \rho g h(q)$; in particular, for $p \in \CU$,  $\rho g  \rho h \iota (p)  {\Sh} \rho g h \iota(p)$;
  so 
  \[  \rho g \rho h \iota  {\Sh}^{Aut(\CU)} \rho g h  \iota \]
  i.e. (using $\iota \rho = \rho$) ${\phi}(g){\phi}(h) {\Sh} {\phi}(gh)$.  Since $End(J) = Aut(J)$, we can write ${\phi}(gh)={\phi}(g){\phi}(h) k$
  for some $k \in Aut(J)$, and we have ${\phi}(g){\phi}(h) k {\Sh} {\phi}(g){\phi}(h)$.   Applying the ${\L}$-endomorphism 
  $({\phi}(g){\phi}(h)) \inv$ we see that $k \in {\Sh}^{Aut(J)}$.  
  
  This shows that ${\phi}:Aut(M) \to Aut(\CU):{\Sh}$ is a quasi-homomorphism.
Composing with the Hausdorffization map $\pi:Aut(\CU) \to \CG(T)$, we see that $\phi:\CG(T): {\Sh}$ 
is also a quasi-homomorphism.

(1)  In view of \propref{2.5}, it suffices to prove now that if $C \subset G=Aut(\CU)$  is bounded with respect to the action on $(\CU,\mu)$, then $W:=\phi \inv(C)$ is bounded with respect to the action on $D(M)$.    Identifying $\CU$
with a subset $J$ of $S(M)$, the retraction $\rho$ shows every point of $S(M)$ is at distance at most $1$
from a point of $J$;   moreover
since $\rho$ changes distances by at most $1$, it is clear that $W$ is bounded with respect
to the action on $S(M)$.  On the other hand, letting $i(a)$ denote the algebraic  type $x=a$, 
we have an $Aut(M)$-invariant embedding $i:M \to S(M)$, and again we saw
using (DP)  that
any element of $S(M)$ is at bounded distance from some element of $i(M)$.   So if $a \in M$, then 
$W i(a)$ has finite radius with respect to $\mu$.  But immediately following 
 \defref{mu}, we saw that for $c,d \in i(M)$, the $m$-distance $m(c,d)$ is at most $7 \mu(i(c),i(d))$.    Thus $W a$ 
 has finite radius with respect to $m$, and so $W$ is bounded   with respect to the action on $D(M)$.

(2)   Similarly, if $W \subset Aut(M)$ is bounded with respect to the action on$(M,m)$, then it is bounded with respect to the induced action on types (to check boundedness we may choose a "center" $a$ lying in $M$.)  It follows as in (1) that the set $\rho \circ W$
has bounded action on $J$; so the image in $\CG$ is precompact by \propref{2.5}.  
\eprf

\begin{rem} \label{error-improve} If we introduce new sorts $\gg_n$ corresponding to $D^n$, and define $\Sh_!^{Aut(\CU)}$   to be the set of automorphisms $\g$ of $\CU$ (on all these sorts) with
 $  \g(a) \Sh a$ for any $a \in \gg_n$, the proof still shows that  ${\phi}:Aut(M) \to (Aut(\CU): {\Sh_!})$.  
 As we will have no immediate use  for this apparent improvement, we will stick with $\Sh$. \end{rem}

\begin{rem} \label{same} In the various constructions of \thmref{qh-basic},  
 we used the {\em fact} that the notion of locality is generated by a definable relation $m$, and we used the  notion of finite distance 
in defining the notion of model and of type space.   But $m$ itself was not further used.  Thus  replacing $m$ by another relation
 $m'$, where $m_1'(x,y)$ implies $m_k(x,y)$ for some $k$, and vice versa,  will change nothing:   the class of models,       
the definition of the local type space, of  the langauge $\L$, of $\CU$, the notion of an embedding and a retraction between $\CU$
and $S(M)$,  and the maps $\phi$ and $\varphi$, all remain the same.           
\end{rem}

 \ssec{Choices}  \label{choices}
 The homomorphism of \thmref{qh-basic} depends not only on $M$, but also on $\iota$ and $\rho$.  
  Letting $\Sigma=S(M)$,
the relevant spaces are $Hom(\CU,\Sigma)$ (for choosing $\iota$), and for $\iota \in Hom(\CU,\Sigma)$, $J=\iota(\CU)$,  the  space of retractions $Hom_J(\Sigma,J) $.     Their study is of considerable interest; in part, I think they are relevant to
Gromov's questions about symmetry breaking in Ramsey theory.    Here we will only mention a  
degree of canonicity that is nevertheless present, despite the choices.   This will not be  relevant to the main applications
to approximate subgroups and lattices.

 As discussed in \cite{patterns}, the choice of $\iota$ amounts to expanding $M$ to a model of a certain universal theory $T_{\td}$.
  $T_{\td}$ is no longer a local theory; it has a relation for each element of $\CU$, and the action of $\Aut(\CU)$
 passes from the models to the langauge itself. 
 
  \begin{rem} \label{choices2}    $\CU$ is itself independent of the choices, and the
 quasi-homomorphism $Aut(M) \to \CU$ (for a saturated $M$) can be made  as unique as   the   saturated model itself:  
in an appropriate model of set theory, we fix a cardinal $\kappa \geq L_{\td}$ with $2^\kappa=\kappa^+$, and let $\CM_{\td}$ be
an existentially closed model of $T_{\td}$, realizing every qf type over an existentially closed substructure of cardinality
$\leq \kappa$.    It was shown in  \cite{patterns} (see 3.17) that  $T_{{\td}}$ is an irreducible universal theory, i.e. the joint embedding property holds; this continues to apply sortwise; so $\CM_{\td}$ is unique up to isomorphism; 
and the reduct $\CM$ to $L$ is a $\kappa$-saturated model of $T$. 
 
 As for $\rho$, the arbitrariness is mitigated by \propref{qh2} (1);  in the   spirit
   of Appendix \ref{categories},  it implies that $\phi$
 gives  a well-defined homomorphism $Aut(M) \to Aut_{\appsp} (\CU_h,\Sh_2)$. 
 \end{rem}
 
 Let us formulate a saturation hypothesis on a model $M$ of $T_{\td}$:   
 
 (SH)  \begin{enumerate} 
 \item $M$ embeds  all models of $T_{\td}$ of size $\leq |\CU|$; 
 \item $M$ is   homogeneous for $L$-elementary embeddings $M_1 \to M$ with  $M_1 \models T$  with $|M_1| \leq  |\CU|$.  (I.e. for fixed $M_1$, all such embeddings are $Aut(M)$-conjugate.)  (In particular, $M$ is orbit-bounded.)
\end{enumerate}

Note (SH) holds if either $M$ is saturated as a model of $T_{\td}$, or else if it 
 is ec and  saturated for qf types of $L_{\td}$.  
 
 \begin{prop}  \label{qh2}  In   \thmref{qh-basic}, we further have: \begin{enumerate} 
   \item If $f,f': Aut(M) \to Aut(\CU)$ are obtained from two choices of retractions $\rho$,$\rho'$, then $f \Sh_2 f'$.
          \item Take $M$ to satisfy (SH).   
            Then 
              $\rho$ can be chosen so that the induced quasi-homomorphism
      $\phi:Aut(M) \to Aut(\CU)$ (and   $\varphi:Aut(M) \to  \CG(T) $ ) are surjective.   
           
    \end{enumerate}
   \end{prop}
 
 \prf  
 
 (1) Let $\alpha \in Aut(M)$.  
 We have  $  \rho \circ \alpha   {\Sh} \alpha$ by \lemref{LSdef} (3).
 and $\rho'  \circ \alpha   {\Sh} \alpha$, so
 $\rho \circ \alpha \Sh_2 \rho' \circ \alpha$ and
 $ \rho \circ \alpha \circ \iota  \Sh_2 \rho' \circ \alpha \circ \iota$.

 (3)  The $L_{\td}$-structure on $M$ can be viewed as a homomorphism $\iota: \CU \to S(M)$.  
 Let $M_0$ be a small   elementary submodel of $M$
 as an $L_{\td}$-structure, and let $r_0: S(M) \to S(M_0)$ be 
 the restriction map.  Then $\iota_0:=r_0 \circ \iota: \CU \to S(M_0)$
 admits, by the retraction property of $\CU$, a one-sided inverse
  $\rho_0: S(M_0) \to \CU$, with $\rho_0   \iota_0=Id_{\CU}$.   
 Let $\rho: S(M) \to \CU$ be the composition  $\rho_0 r_0 $.
 Let $f(\alpha) = \rho \circ \alpha \circ \iota$.  
 
 \claim{} $f$ is surjective. 
 
  To see this, identify $\CU$ with $J=\iota(\CU)$, so that $\iota$ becomes the inclusion map of $J$ in $S(M)$,
  and $\iota_0 = r_0|J$.
  Let $\beta \in Aut(J)$.  We can now view $\beta$ also as a homomorphism $J \to S(M)$.  
  According to \cite{patterns} A.11(3), there exists $\si \in Aut(M)$ with $r_0 \circ \si |J =r_0 \circ \beta$.  
  So $f(\sigma) =  \rho   \si   \iota
  = \rho_0 r_0   \si  \iota = \rho_0 r_0 \beta \iota =\rho_0 \iota_0 \beta \iota  = \beta$.   Since $\beta \in Aut(J)$
  was arbitrary, $f$ is onto.

 \eprf

  \ssec{The pullback of $\Sh$.}

Define  the binary relation $\Sh_k$ to be the $k$-fold composition of  ${\Sh}$ with itself; it is defined by a pp formula (the existential quantifiers are local since
 $\Sh(x,y)$ implies $m_2(x,y)$,   \lemref{shrem}.)  It is easy to check that the same expression
  defines
 the $k$-fold iteration of $\Sh$ on $\CU$ (or see \cite{patterns} 3.12.)  
 
For any $\mathcal{A},\mathcal{B} \modec {\CT}$, we define a relation $\Sh_k$ on 
$Hom(\mathcal{A},\mathcal{B})$ by
\[h_1 \Sh_k h_2 \iff (\forall a \in \mathcal{A} ) (h_1(a),h_2(a)) \in \Sh_k^{\mathcal{B}} \]
For an element of $Aut(\mathcal{A})$, we write $\Sh_k(h)$ for $h \Sh_k 1$, where $1=Id_{\mathcal{A}}$.
For $h_1,h_2 \in Aut(\mathcal{A})$, we have $h_1 \Sh_k h_2$ iff  $(h_1(a),h_2(a)) \in \Sh_k^{\mathcal{A}}$
for all $a$; by definability of $\Sh_k$, this is iff   $(a,h_1 \inv h_2(a)) \in  \Sh_k^{\mathcal{A}}$ for all $a$ iff $\Sh_k(h_1 \inv h_2)$. 
 It follows that  $\Sh_k$ is a normal, symmetric subset  of   $Aut(\mathcal{A})$.  
 
 The above notation applies to $Aut(\CU)$ and to $Aut(S(M)) = Aut(M)$ (via its action on $S(M)$.)
 
 Finally, we define $\Sh_k^{\CG(T)} $ to be the image of $\Sh_k^{Aut(J)}$ under the quotient homomorphism
$Aut(J) \to \CG(T)$.     

Note that in the groups $Aut(J)$ or $\CG(T)$ we have  $\Sh_k \Sh_{k'} \subset \Sh_{k+k'}$, but (presumably) equality need not hold in general. 

\begin{prop}\label{qh3}  In \thmref{qh-basic}, one can add:
\begin{enumerate} 
 \item  $\varphi(\Sh_n^{Aut(M)}) \subset \Sh_n^{\CG} $
\item   If $\phi(g_1)\phi(g_2)\inv \in \Sh_n^{Aut(\CU)}$, then $g_1 g_2 \inv \in \Sh_{n+4}$.    
\item  $\phi^{-1} \Sh_n^{Aut(\CU)} \subset \Sh_{n+3}^{Aut(M)}$
 \item  $ \varphi^{-1} \Sh_n^{\CG(T)} \subset \Sh_{n+5}^{Aut(M)}$. 
\item     If $\varphi(g_1)\varphi(g_2)\inv \in \Sh_n^{\CG}$, then $g_1 g_2 \inv \in \Sh_{n+6}$. 
  
\end{enumerate}
\end{prop}
    \prf
    Identify $\CU$ with $J \leq S(M)$, and let $\rho: S(M) \to J$ be a retraction.
 
(1)  If $g \in \Sh^{S(M)}$, so that $g(p) \Sh p$ for all $p \in S(M)$, then for $p  \in  J$ we have 
  ${\phi}(g)(p) = \rho g(p) L_n \rho(p) = p$.  Thus  ${\phi}(\Sh) \subset \Sh$.  Similarly, $\phi(\Sh_n) \subset \Sh_n$.

(2)  Let  $g_1,g_2 \in Aut(M)$ and assume $\phi(g_1) \Sh_n \phi(g_2)$ in $Aut({\CU)}$.  Let $p \in S(M)$.   
By \lemref{LSdef}(3),  $p \Sh \rho(p)$.   Thus (using the fact that $g_1,g_2,\rho$ are homomorphisms),
 
 \[ g_1(p) \Sh \rho g_1(p) \Sh \rho g_1 (\rho(p))   \Sh_n \rho g_2 (\rho(p))  \Sh \rho g_2(p) \Sh g_2(p) \]
 
So $g_1(p) \Sh_{n+4} g_2(p)$, and since $p$ was arbitrary, $g_1 g_2 \inv  \in  \Sh_{n+4} $.  

 (3)  The penultimate step in the displayed formula in  (2) can be skipped when $g_2=Id$, since $\rho \rho p = \rho p$.

(4)  
Now suppose $\varphi(g) \in \Sh_n^{\CG}$; i.e. $\phi(g) = \pi(h)$ for some $h \in \Sh_n^{Aut(J)}$.
 If $g \in \ker \pi$ then $g(p),p$ are inseparable by open sets.  By \lemref{LSdef} (4), we have $(p,g(p)) \in \Sh_2$ so 
 $\ker \pi \subset \Sh_2$.
So ${\phi}(g) h \inv \in \ker \pi \subset \Sh_2$.
Hence ${\phi}(g)  \in \Sh_{n+2}$, and by (3), $g \in \Sh_{(n+2)+3}$.

(5) Similarly, if $\varphi(g_1)\varphi(g_2)\inv \in \Sh_n^{\CG}$, then by (2) we have $g_1 g_2 \inv \in \Sh_{(n+2)+4}$.  
 
    \eprf

The example below is adapted from the basic component of the  example of \cite{clpz} of a non-compact Lascar group; local logic permits a simpler presentation.
It is also in essence the same as Example 3 of a quasimorphism in \cite{kotschick}, where one can find an illuminating discussion of its provenance in geometry.

\begin{example}   \label{milnor}
Let $T=Th_{loc;\mu_1}(\Rr,<,S)$.   Here $S(x)=x+1$,   
  \[\mu_1(x,y) \equiv x \leq y \leq S(x) \vee y \leq x \leq S(y),\]
and  the subscript loc designates that we take the theory in local logic, i.e. using   quantifiers with bounded range only.

 The space $\Sigma$ of (local) $1$-types over $M=(\Rr,<,S)$ can be identified with $\{-,0,+\} \times \Rr$;
 here $(0,\alpha)$ denotes the algebraic type $x=\alpha$; $(-,\alpha)$ denotes the type just below $\alpha$, 
 and dually $(+,\alpha)$.  The two non-local 0-definable types, at $\infty$ and $-\infty$, are of course  not part of the space $\Sigma$.    
 Let $J =   \{-,0,+\} \times \Zz$, a ${\L}$-substructure of $\Sigma$.    Pick a  cut $\gamma \in (0,1]$     
  and define 
 a ${\L}$-retraction $r:\Sigma \to J$ as follows:  for any $n \in \Zz$, $\delta \in [0,1)$  and $* \in \{-,0,+\}$,  
 let $r(*, n+\delta)=(*,n)$ if  $\delta \leq \gamma$ and $r(*, n+\delta)=(*,n+1)$ if 
 if $\gamma < \delta $.      It is easy to see that no further collapse is possible, and $J$ is the universal e.c. model of $\CT$.
 We have $Aut(J)=\Zz$.

$Aut(M)$ is  the much  bigger group  of order preserving maps $\Rr \to \Rr$ commuting with $S$.  
The quasimorphism on $Aut(M)$, associated with the above retraction, is the translation number\footnote{I take this term from \cite{kotschick}.  Emmanuel Breuillard pointed out to me that these ideas go back to Poincar\'e's {\em rotation number}.
A homeomorphism $\rho$ of the oriented circle $\Rr/\Zz = M/S^{\Zz}$ lifts uniquely to an element $r$ of $Aut(M)$ fixing a given point $m_0$;
the rotation number of $\rho$ is the translation number of $r$.}.

  Restricted to the group of 
translations $\Rr$ as a subgroup of $Aut(M)$, we obtain a choice of 'nearest integer' map, which is of course
a quasi-homomorphism $\Rr \to (\Zz;\{-1,0,1\})$.  

In connection with the discussion in \secref{choices}, we note also that in this example, we have
    $Hom(\CU,\Sigma)=\Rr$ and $Hom_J(\Sigma,J) $ is the space 
  of cuts in $(0,1)$. 

\end{example}


 \begin{question}  Is it possible to develop an analogous theory   uniform in the doubling constant?    Assume every $2$-ball
 is uniformly a union of a $Y$-parameterized family of $1$-balls, where the definable set $Y$ is in some sense
 `small' compared to $D$; at a minimum, there should exist types over $M$ that properly increase $M$ without increasing
 $Y$, so that finite distance still implies distance $1$ from some point.    According to a theorem of Shelah \cite{shelah-2},  such types exist in pseudo-finite theories provided  $|{D}|$ is not polynomially bounded by $|Y|$.  Thus in theory, this could 
 give a possible route to  various questions of polynomial bounds.    Developing a theory on this basis will not be easy, but  even   much stronger assumptions of relative smallness would be of interest.   The theory of locally compact groups  would need to be generalized, however, in any case, with the compact/ noncompact dichotomy replaced by a model-theoretic analogue.\end{question}

 \pb
 \section{Definable groups}    \label{definablegroups}

Recall Tao's  definition of   approximate subgroups.

\begin{defn}[\cite{tao}] \label{ag-tao} Let $G$ be a group.  A symmetric subset $\Lam \subset G$ is an  {\em approximate subgroup} 
if $\Lam \sim \Lam^2$.
\end{defn}

\begin{thm} \label{grmain}  Let $G$ be a group,  generated by an approximate subgroup $\Lam$.  Then 
there exists a  second countable 
locally compact topological group $\lc$,  
  a compact normal subset $\Sh \subset \lc$ 
    and a quasi-homomorphism 
    \[f: G \to \lc:\Sh\]  
  such that: 
    \begin{enumerate} 
\item For $C \subset \lc$ compact,  $f \inv(C)$ is contained in some $\Lambda^i$.  
\item For each $i$ there exists a compact $C \subset \lc$ with $\Lam^i \subset f \inv(C)$.   
\item Specifically, $f \inv(\Sh) \subset \Lam^{12}$.
\item  Let $X,X'$ be compact subsets of $\lc$, with $\Sh^2X \meet \Sh^2 X' = \emptyset$.  Then there exist  disjoint definable subsets $D,D'$ of $\Lam^k$ for some $k$ (with parameters in $\Lam$) such that $f \inv(X) \subset D$, and $f \inv(X') \subset  D'$.  
\end{enumerate}
 
  \end{thm}
  
  The purely group-theoretic properties (1,2) will be the main engine  of the group-theoretic applications, including \thmref{ag2} and \thmref{discrete1}; the second-countability will not be needed.
  
  Note that $\lc$ is compact iff $\Lam$ generates $G$ in a finite number $s$ of   steps.  Even in this case, when $s$ is large,
  we obtain a nontrivial statement due to (3); in particular, $f: G \to \lc:  f(\Lam^{12})$.  See \corref{slowly2}.

    (4) should be compared to the familiar form of continuity or continuous-logic definability of a map into a compact space,
  asserting that the preimages of two disjoint compact sets can be separated by a definable set.  Here this is  asserted for
  two compacts at some discrete distance from each other.

      We will begin with a proof of (1,2); 
the proofs of (3,4) are postponed till we discuss definability.%

  \prf  We define a local structure  $M$.  The universe $D(M)$ of $M$, as a set,  is a copy of $G$. 
    The   relations on $D$ include the binary relation $ y x \inv \in \Lam$, denoted $m_1$: 
  and the quaternary relation $yx \inv  = z w \inv \in \Lam$.   The locality generator is $m_1$.     There are no further 
  relations in the structure $M$; thus the action of $G$ on $D$ from the right,
  $(g,d) \mapsto d g \inv$, is an automorphism of $M$.   This gives a homomorphism $\delta: G \to Aut(M)$. 
  Since $G$ is transitive on $D$, it is clear that $M$ is orbit-bounded.    The balls can be identified
  as $m_1(a)=\Lam a$, and $m_2(a) = \Lam^2 a$.   As $\Lam$ is  an approximate subgroup,
  $\Lam^2 \subset \union_{i=1}^k \Lam c_i$, so $m_2(a) = \union_{i=1}^k  m_1(c_ia)$; this gives the doubling property (DP).
  Thus  \thmref{summary1} applies:  $Th(M)$  admits a core $\CU$, and  $Aut(\CU)$   has a natural  locally compact topology;
  the canonical Hausdorff quotient  $\CG(T) : = Aut(\CU)/\fg_D$ is a locally compact topological group;  we also
  defined a canonical compact,  normal, symmetric subset $\Sh \subset \CG(T)$.    There exists an embedding
  $\iota: \CU \to S(M)$, and a retraction $\rho: S(M) \to \CU$.
According to 
  \thmref{qh-basic}, 
  there exists a    quasi-homomorphism ${\varphi}:Aut(M) \to (\CG(T): {\Sh})$  such that
  compact subsets of $\CG$ pull back to bounded subsets of $Aut(M)$, and bounded subsets of $Aut(M)$ push forward to
  precompact subsets of $\CG$.   
  
  Let $f = \varphi \circ \delta$.   If $C \subset \CG$ is compact, then $\delta(f \inv(C))$ is bounded in $Aut(M)$, being a subset of $\varphi \inv(C)$.  
Picking $d_0 \in D$  to be the identity element of $G$,  we have $\delta(f \inv(C)) (d_0) = d_0 f \inv(C) = f \inv(C)$; so 
$f\inv (C)$ is 
  bounded in $D(M)$.  This means that $f \inv(C) \subset \Lam^i$ for some $i$.    Conversely, $\delta(\Lam^i)$ is by definition bounded, so $f(\Lam^i) = \varphi(\delta(\Lam^i)) $ is precompact, hence $C= cl(f(\Lam^i))$ is compact and $\Lam^i \subset f \inv(C)$.

To obtain second countability, we modify the end of the construction slightly.   Recall $\varphi$ was the composition
$\pi \circ \phi$, with $\pi: Aut(\CU) \to \CG$ the quotient map, and $\phi: Aut(M) \to (\CG(T): {\Sh})$ as constructed in \ref{qh-basic}.
Let $P$ be the complete atomic type of $\rho(d_0)$.
 Let $\fg=\fg_{P}$ be the infinitesimal subgroup with respect to the action on $P$.    By \propref{autlc}, $\fg_P$ is a closed, compact subgroup of $Aut(\CU)$, and  $\CG(P) $ 
 is a second countable, locally compact topological group.  Let $\pi_P: Aut(\CU) \to \fg_P$ be the projection.
 We let $\varphi_P = \pi_P \circ \phi$, and again $f=\varphi_P \circ \delta$.  We still have the above properties, and now
 according to \propref{autlc},  $\CG(P) $ 
 is also second countable.
 
  \eprf

 We will also mention a slightly more general version.     The proof is identical once we have a  generalization of the local logic of
  \secref{local}, where we do not assume that the locality relation is generated by a single relation; we will not repeat details but  will describe the changes needed in the construction.

    \begin{prop}\label{grmain2} 
     Let $G$ be a group, $\om$ a   family of 
     commensurable approximate subgroups of $G$.    Assume:
    if  $\Lam,\Lam' \in \om$,  then  $\Lam \union \Lam' \in \om$, and $\Lam^2 \in \om$; and $\union \om = G$.
  Then there exists a  locally compact topological group $\lc$,  
  a compact normal subset $\Sh \subset \lc$ 
    and a quasi-homomorphism 
\[{f}: G  \to \lc:\Sh\]
  such that: 
    \begin{enumerate} 
\item For $C \subset \lc$ compact,  and $\Lam \in \om$,    $f \inv(C) $ is covered by a finite union of right translates of $\Lam$. 
\item For   $\Lam \in \om$ there exists a compact $C \subset \lc$ with $\Lam \subset f \inv(C)$.   
  \item If $\Lam_0 \in \om$ and $X \subset \Lam \in \om$, we have $f \inv f(X)  \subset X \Lam_0^{12} $.
 If $f(g_1) \in f(g_2) \Sh_k $ then $\delta(g_1g_2 \inv) \in \Sh_{12+2k}$.
  \item  Let $X,Y$ be compact subsets of $\lc$, with $\Sh^2X \meet \Sh^2 Y = \emptyset$.  Then there exists a definable subset $D$ of  some $\Lam \in \om$ (with parameters in $\Lam$) such that $f \inv(X) \subset D$, and $f \inv(Y) \subset \Lam \setminus D$.  
  \item If $\om$ is countable, we may take $\lc$ to be second countable.
 \end{enumerate}
  \end{prop}

\prf      

Define a structure $M$ as in \thmref{grmain}, with universe $D=G$;  for {\em each} $\Lam \in \om$, we posit
a  binary relation $m_\Lam: \  y x \inv \in \Lam$; and a quaternary relation $yx \inv  = z w \inv \in \Lam$.   A subset of $D$ is viewed
as local if it is contained in a ball $m_\Lam(a)=\{b: m_\Lam(a,b)\}$; or  equivalently, in the union of  finitely many such balls.
As before, only existential quantifiers of the form $(\exists x)m_\Lam(x,y) \wedge \cdots$ are allowed in the formation of pp formulas,
but now all $\Lam \in \om$ are allowed.  
 If we pick one $\Lam_0 \in \om$, 
it remains true that any ball $m_\Lam(a)$ is covered by finitely many balls $m_{\Lam_0} (b_i)$.     A weak model is local if
 for any $a,b \in D(M)$, $m_\Lam(a,b)$ holds for some $\Lam \in \om$.  A type $p(x)$  is local if $m_\Lam(x,y)$ is represented
 in $p$ for some - or equivalently, by DP, for all - $\Lam \in \om$.  The construction of $\CU,Aut(\CU), \fg,f$ and proof of their properties is identical.  

\eprf

\begin{rem} \label{4.3} Let $\om'$ be any family of commensurable approximate subgroups of $G$,  
     such that if $\Lam \in \om$ and $g \in G$ then $g \inv \Lam g$ is commensurable
     with $\Lam$.  Then $\om'$ may be completed to  a class $\om$ satisfying the hypotheses of \propref{grmain2},
     whose union is $G$.
     \end{rem}
     \prf  
If   $\Lam \in \om$, then any two-sided translate $a \Lam b$ is contained in a finite union of right translates of $\Lam$:
we have
$a \Lam a \inv  \subset \union_i \Lam c_i$
for some $c_1,\ldots,c_l$; and $a \Lam b = (a \Lam a \inv) ab \subset \union_i \Lam c_iab$.   In particular, each
$c_iab\Lam$ is contained in a finite union $\union \Lam d_{ij}$, so $a \Lam b \Lam  c  \subset \union _{i,j} \Lam^2 d_{ij}$,
and it follows that $a \Lam b \Lam c$ is contained in a finite union of right translates of $\Lam$.  From here it is 
clear that closing $\om$ under the operations $\Lam,\Lam' \mapsto \Lam \union \Lam'$, $\Lam \mapsto \Lam^2 $,
and $\Lam \mapsto (a \Lam \union \{1\} \union \Lam a)$
will not  damage the hypotheses. \eprf

\ssec{Definability}  We will now place our approximate subgroups in a definable setting more explicitly.      A complete first-order theory $T$
is given with predicates $\Lam_i$ ($i \in I$, where $I$ is a partially ordered index set. )      If $i\leq j$, a definable injective map
$\Lam_i \to \Lam_j$ is given, forming a direct limit system; the direct limit will be denoted $G$.  We assume that for each $i,j \in I$ for some $k \in I$, a definable map  $\Lam_i  \times \Lam_j \to \Lam_k$ is given; inducing a group multiplication on $G$.  
  In this situation $G$ is called a piecewise definable group of the theory $T$;  or a strict Ind-object in Grothendieck's category of   Ind-definable sets.  
We identify $\Lam_i$ with a subset of $G$, and assume further that any two $\Lam_i$ are commensurable in $G$;
in particular each $\Lam_i$ is an approximate subgroup. 

 We say that a subset of $G$ is
{\em locally $\bigwedge$-definable} if the intersection with each $\Lam_i$ is $\bigwedge$-definable.  

Given a theory $T$ with a definable commensurable family of approximate subgroups, as above, 
we introduce a new local sort $D$, a torsor for $G$.   For $M \models T$, $D(M)$ is a copy of $G(M)$.  The  new relations introduced
are   binary relations $m_i(x,y): \ \ \equiv \ \ ( y x \inv \in \Lam_i)$ on $D$, and the relation $yx\inv = u \in \Lam_i$ on $D^2 \times \Lam_i$.  
We   define a locality structure on $D$, using the $m_i$.   Restricting attention to 
the induced (local) structure on $D$, we form $\CU$ (in the corresponding sort $\gg$.)      Denote by $L^+$ 
the langauge obtained from $L$ by adding the local sort $D$ and these relations; $M^+$ the result
of adding the torsor $D$ to the structure as above; and $T^+:= Th(M^+)$.


Next recall the universal theory $T^+_{\td}$ of \cite{patterns} (see 3.17).  A model $N^+_{\td}$ is nothing more   than a model 
$N^+$ of $T^+$, along with a homomorphism $\iota: \CU \to S_D(N^+)$.    The language $L^+_{\td}$ is obtained from $L^+$ by adding
a new relation $(d_px)\phi$
for each $p \in \CU$ and each $\phi(x;y_1,\ldots,y_k) \in L^+$, to be interpreted as a subset of $D^k$, namely  as
$\{(b_1,\ldots,b_k) \in D^k:  \phi(x,b_1,\ldots,b_k) \in \iota(p)\}$.  Fix some $p_0 \in \CU$;  the choice will not matter, up
  to bi-interpretation.  We  have in $L^+_{\td}$ a
unary relation $m_k^{p_0}$ in $L^+$ such that $m_k^{p_0}(b)$ holds iff $m_k(x,b) \in p_0$.   We will only consider models of
$T^+_{\td}$ that omit the partial type:  $\neg m_k^{p_0}(x,b), k=1,2,\cdots$.  This  amounts to treating  the sort $D$ of 
 $T^+_{\td}$ (a single local sort in $T^+$)  as the union of infinitely many `sorts' defined by $m_k^{p_0}$,  with inclusion maps among them.   Any model of
  $T^+_{\td}$ (in this sense), restricted to an $L$-structure,  is a (local) model of $T^+$.  
  
  Since some model of $(T^+_{\td})^{\pm}$
  is a model of $T$, if $T^+_{\td} \models \neg (\exists x,y,z)(\alpha(x,y) \wedge \gamma(y,z))$ with $\alpha \in L^+$,
  and $T^+ \models  (\forall y)(\exists x)(\alpha(x,y))$, then   
$T^+_{\td} \models \neg (\exists y,z)( \gamma(y,z))$.)   Thus $T^+_{\td} \modec T^+$.  
%
%

In \thmref{grmain} and \propref{grmain2}, we used the action of $G$ on $D$ by right translation.  This action depends on a choice
of an element $d_0$ of $D$, used to identify $G$ with $D$ via $g \mapsto gd_0$.    The homomorphism $\delta:G \to Aut(M^+)$
is thus also dependent on this choice, and may be better denoted $\delta_{d_0}$.   However when $X$ is a normal subset of $Aut(M^+)$, $\delta \inv(X)$ is well-defined and does not depend on the choice of $d_0$.  In particular, 
$\delta \inv(\Sh_k)$
 is well-defined.   (Where $\Sh_k \subset Aut(M^+)$ is defined above \lemref{qh3}.)  Let us compare $\delta \inv(\Sh_k)$
 to $\Sh_k$, as computed directly in $T$ (i.e. in $G$ and not in $D$.)

Consider 
\[P:= \{ab \inv:  a,b \in G(M), a \Sh b\} \] 
Here $\Sh$ is  computed in  $\aleph_0$-saturated models of $T$.

While there is no multiplication on $D$, we can define $ab \inv \in G$ for $a,b \in D$, so that $(ab \inv) b=a$.  Then
$P=\{a b \inv:   a,b \in D(M), a \Sh b\} $.  (If $a,b \in D$ and $tp(a/N)=tp(b/N)$, then picking $d \in D(N)$ we have $ab\inv = (a d\inv) (b d \inv) \inv$
     and $a d \inv \Sh^M b d \inv$.)   From this it follows easily that  $\delta \inv (\Sh_k^{Aut(M^+)}) \subset P^k$.  
     
     On the other hand, by locality and the doubling property, if $a,b \in D$ and $a \Sh^D b$ then $ab\inv \in \Lam^2$.  (Since $tp(a/N)=tp(b/N)$
     implies that $a,b$ lie in a common $1$-ball.)  
     
     By \propref{qh3} (2), we obtain:  

 \begin{prop}\label{qh3b}  In   \ref{grmain} and \ref{grmain2},   one can add:
If $f(g_1) f(g_2) \inv \in \Sh_k^{\CG}$ then  $\delta(g_1g_2 \inv) \in \Sh_{k+6}$ so
 $g_1 g_2 \inv \in P^{k+6} \subset \Lambda^{2(k+6)}$.  Likewise, if $f(g) \in \Sh_k^{\CG}$ then $\delta(g) \in  \Lam^{2(k+5)}$.
\end{prop}
 
 In particular, this proves  \thmref{grmain} (3).
 
 \begin{rem} All conjugates of $P$ are contained in $P^3$.    Hence $P \subset Q \subset P^3$ for some normal
 $\bigwedge$-definable set $Q$.  \end{rem}
 \prf Let $c \in {P}$, and $g \in G$.  We have $c=c_1 \inv c_2$ with $tp(c_1/N)=tp(c_2/N)$.  
 Take $N$ to be $\aleph_1$-saturated, and $N_0$ a countable elementary submodel; then $tp(g/N_0) = tp(m/N_0)$
 for some $ m \in N$; so  $ m \inv g \in {P}$.  
   Now $tp(m \inv c_1 m/M ) = tp(m \inv c_2 m/M)$,  so
 $m \inv c m \in {P}$.   Thus $g  \inv c g  = l \inv ( m \inv c m) l \in P^3$.     Now let $Q$ be the intersection of all conjugates of $P^3$.\eprf

\begin{prop} \label{grmain3} \label{4.8} 
Let $G$ be a piecewise definable group in a theory $T$, limit of definable approximate subgroups.    
Let $M^+_{\td} \models T^+_{\td}$, let $M^+$  be the reduct to $L^+$, and $\iota: \CU \to S(M^+)$  the embedding corresponding
to $M^+_{\td}$.  Further let   $\rho: S(M^+) \to \CU$ be a retraction, and let $d_0$ be an element of $D(M^+)$,  
Let $\delta: G \to Aut(S(M)^+)$ the map induced by right multiplication
if $(D,d_0) $ is identified with $(G,1)$:   $\delta(g)(hd_0) = h g \inv d_0$.   
      Finally let 
$\phi: Aut(S(M)^+) \to Aut(\CU)$ be defined by $\g \mapsto \rho \circ \gamma \circ \iota$,    $\varphi$ be the composition  of $\phi$
with $Aut(\CU) \to \CG$, and $f= \varphi \circ \delta: G(M) \to \CG$.   
Then:
\begin{enumerate} 
\item    Let $W$ be a compact subset of $\CG$.  Then there exists a   $\bigwedge$-definable 
   $W^* \subset G$  in $L_{\td}^+(d_0)$ with 
 $f \inv(W) \subset W^* \subset f \inv(W \Sh)$.   $W^*$ is defined by universal formulas of $L^+_{\td}(d_0)$.     
  
 \item 
 Let $W_1,W_2$ be compact subsets of $\CG$ with $\Sh W_1 \meet  \Sh W_2 = \emptyset$.  Then  
 there exist  definable sets $V_1,V_2 \subset G$ in $L_{\td}^+(d_0)$ with $f \inv (W_i) \subset V_i$ and 
 $V_1   \meet  V_2= \emptyset$.

\end{enumerate}
\end{prop}

\prf   
{\em During the present proof, the word {\em definable} and its attendants will refer to $L^+_{\td}(d_0)$.}

We show first that we may assume some saturation of $M^+_{\td}$; this will be needed for (2).  
  Let $((M^+_{\td})^*,d)$ be an elementary extension of   $(M^+_{\td},d)$.    Define $\rho^* = \rho \circ r$, 
where $r: S(M^+) \to S(M)$ is the restriction map.    We obtain $\iota^*,\delta^*,\phi^*, \varphi^*$;   and it is easy to verify that
the various maps commute, and we have:  $\phi \circ \delta = \phi^* \circ \delta^* |M$, and $f = f^* |M$.  Thus 
by passing to an elementary extension 
we may assume $M^+_{\td}$ is $\aleph_0$- saturated, and $M^+$ is $\aleph_0$- homogeneous.

We begin with the proof of (1).     

  Identify $\CU$ with  $\iota(\CU)$.    Let ${\f}$ denote the composition  
  \[ G(M) \to_{\delta}  Aut(M^+) \to  Aut(S_D(M^+)) \to_{\phi} Aut(\CU) \]
so that $f = \pi \circ {\f}$, with $\pi:  Aut(\CU) \to  \CG$.
 We will write $\hat{g}$ for $\delta(g)$ or for the image of $g$ in $Aut(S_D(M^+)) $ or $Aut(\CU)$.
The homomorphism $Aut(\CU) \to \CG$ is continuous and proper; so it suffices to show that for each closed bounded subset $V$ of $Aut(\CU)$, there exists a $\bigwedge$-definable set $V^*$ with
\[ {\f} \inv(V) \subset V^* \subset {\f} \inv(V \Sh)  \ \ \ \ \    \]
Assuming the displayed statement  holds for some family of  sets $V_i$, closed under finite intersections,  it 
will also be true for  $\meet_i V_i$,   with 
$(\meet_i V_i)^* = \meet_i V_i^*$; this, since by compactness $\meet_i V_i \Sh = (\meet_i V_i) \Sh$.  
  Hence, recalling that basic closed sets are closed under finite intersections, 
   it sufices to take $V$ to be a basic closed set in $Aut(\CU)$, i.e. 
\[ V = \{\g \in Aut(\CU):  R(a,\g(a)) \} \]
Here (using \thmref{summary1} (3)) $R$ can be taken to be a basic relation of $\L$ ensuring that $\g$ lies in a bounded subset of $G$, when $a$ is fixed; 
 $a$ is a tuple of elements of $\CU$, and $\g((a_j: j \in J)) = (\g(a_j): j \in J)$.  While $R$ will  make use of only finitely many of the $a_j$ and $\g(a_j)$, it will be convenient to allow $a$ to be an infinite tuple, indeed a tuple enumerating all of $\CU$.  Let $p$ be the atomic type of $a$ in $\CU$.  For two tuples $(a_j: j \in J), (b_j: j \in J)$, write $\Sh ^1(a,b)$ if $\Sh(a_j,b_j)$ holds for each $j$.  
Define a subset $V^*$ of $G$ thus:
\[ V^* = \{ g \in G:  S(M^+) \models  (\exists \xi) (p(\xi) \wedge (\xi \Sh \hat{g}(a)) \wedge R(a,\xi) )\}  \]

Let $g \in G(M)$, $\g={\f}(g)$ the image of $g$ in $Aut(\CU)$. We have $\g= (\rho \circ \hat{g}) | \CU$.

 If $\g \in V$,   then   letting $\xi = \g(a)$,
we see (using the definition of $V$ and \lemref{LSdef}(3))   that $g \in V^*$.   Thus   ${\f} \inv(V) \subset V^* $.

In the opposite directionm assume   $g \in V^*$.   Then
 \[ \CU  \models  (\exists \xi) (p(\xi) \wedge (\xi \Sh \rho\hat{g}(a)) \wedge R(a,\xi) ) .\]
Let $b$ be a witness for $\xi$.  Since $p(b)$ holds,  $a \mapsto b$ defines an $\L$-homomorphism $\CU \to \CU$; 
by \propref{prop2.1}, it is an  automorphism of $\CU$, call it $\g'$.  Since $\g ' (a) = b$ and  $b \Sh \g (a) = \g (\g') \inv (b)$, we have $\g    (\g') \inv   (b) \Sh b $ and so $b \Sh \g' \g \inv (b)$.  
As $b$ enumerates $\CU$,  $\g '   \g  \inv \in \Sh^{Aut(\CU)}$.
But $R(a,b)$ shows that $\g' \in V$.  So $\g \in V \Sh$.  
 
Note that $V^*$ is bounded.  It remains only to see that $V^*$ is (locally) $\bigwedge$-definable.   Using local compactness of $S(M)$ we see that $V^*$ is the intersection of all approximations to it, obtained by replacing $\Sh$ by a definable approximation inside the matrix of 
 definition of $V^*$;   these approximations can be chosen so as to still ensure that the set of $g$ they define is bounded.
 
   Hence 
  $V^*$ is a conjunction of
 pp relations concerning $(a, \hat{g}(a))$.   By \thmref{summary1} (3),  $V^*$ can also be written as an (infinite) conjunction of
 atomic formulas $\R_t$; so  it suffices to prove that each such formula is  locally  definable.
 To simplify notation, consider the atomic formula $\R_t(p,\hat{g}(q))$, 
 $t=(\phi_1(x_1,u),\phi_2(x_2,u))$, with $\phi_1,\phi_2$ formulas of $T^+$.  Here  $p,q \in \CU$ are fixed and $g$ is the variable.  We have 
 $\R_t(p,\hat{g}(q))$ iff there is no $c$ with $\neg \phi_1(x,c) \in p$ and $\neg \phi_2(x,c) \in \hat{g}(q)$.  
 Now $\phi_2(x,c) \in \hat{g}(q)$ iff $\phi_2(x,\hat{g \inv}(c)) \in q$.  
 Thus 
 \[\R_t(p,\delta(g) (q)) \iff \neg (\exists c) 
 (d_{p} x)\neg \phi_1(x,c) \wedge (d_q x)\neg \phi_2(x,\hat{g \inv}(c)) \] 
The right hand side  is   definable in  $L^+_{\td}(d)$.  
 
(2) Let $W_1^*,W_2^*$ be as in (1).   Then $W_i^*$ is   $\bigwedge$-definable (in particular contained in a definable
approximation to $G$), and 
$W_1^* \meet W_2^* \subset f\inv(W_1 \Sh) \meet f \inv(W_2 \Sh) = \emptyset$.   If  $W_1^*,W_2^*$ were not  separated by 
   definable sets $V_i \subset G$ (in $L^+_{\td}(d_0)$), by compactness and since every type over $\emptyset$ is
   realized , we would have $W_1^* \meet W_2^* \neq \emptyset$, a contradiction.  

\eprf

\begin{rem} \label{beth} 
    The Beth-Craig-Robinson definability theory is valid for $\bigwedge$-definable relations.  Here is the version we need:
 Let  $M$  be an $\aleph_0$-homogeneous $L$-structure.  Let $M'$ be an expansion to an $L'$-structure, in which every type is realized.  Let $X,Y$ be  $\bigwedge$-definable subsets in $M'$.  Assume $X \meet \si Y = \emptyset$ for every automorphism
 $\si$ of $M$.  Then $X,Y$ are separated by a 0-definable set of $M$.   [Proof:  Let $SX = \{tp(a; M'):  a \in X\}$, and similarly $SY$.
 Then $SX,SY$ are compact sets in the space $S'$ of $L'$-types (over $\emptyset$.)  Let $\pi: S' \to S$ be the projection to the space of $L$-types.  Then $\pi SX, \pi SY$ must also be compact sets, and they are disjoint since if $tp_L(a)=tp_L(b)$ then there exists
 an $L$-automorphism $\si$ with $\si(a)=b$.  Hence there exists a 0-definable set $D$ of $L$ separating $\pi SX, \pi S Y$ and thus also $X,Y$.]   
 
  It follows in particular that if  $X$ is $Aut(M)$-invariant, then $X$ is $\bigwedge$-definable in $L$.
 \end{rem}

Definability in the following corollary (and beyond) refers again to $L$, not to $L_{\td}$.

\begin{cor} \label{G00} In \propref{grmain3}, and hence also in \ref{grmain} and \ref{grmain2},  one can add that 
$f \inv( cl(\la \Sh \ra))$ is    locally 
$\bigwedge$-definable in $L$.   In fact it is locally an intersection of 0-definable sets commensurable with $\Lam$,
and the smallest subgroup with this property.
 \end{cor} 

\prf  As in \propref{grmain3}, we may assume $M$ is $\beth_2^+$-saturated  in $L_{\td}$, and $\aleph_0$-homogeneous in $L$. 
We will show that $f \inv( cl(\la \Sh \ra))$ is $Aut(M)$-invariant.   

We can identify $Aut(M)$ with a subgroup of $Aut(M^+)$, namely the subgroup fixing $d_0$.  Under this identification,
we have, for $\si \in Aut(M)$ and $g \in G(M)$,
\[  \delta(\si(g)) = \si \inv \delta(g) \si \]

We have  $f= \varphi \circ \delta: G(M) \to \CG: \Sh$; where $\varphi: Aut(M) \to \CG: \Sh$ and $\delta: G \to Aut(M)$.  

 Let $\si \in Aut(M)$,  $g \in G(M)$, and let $Y$ be a normal subset of $\CG$.
Note first   that 
 $\si \inv( \varphi \inv(Y) ) \si  \subset \varphi \inv (Y \Sh^3)$.  Composing with $\delta$, it follows that
 $\si  (f \inv(Y)) \subset f \inv(Y \Sh^3)$.  In case $Y = Y \Sh$, this shows that $f \inv(Y)$ is $Aut(M)$-invariant.

If in addition $Y$ is closed, by   \propref{grmain3}, $f \inv(Y)=Y^*$  is also   locally $\bigwedge$-definable
in $L_{\td} $;  i.e.  $\Lam\meet f \inv(Y)$ is  $\bigwedge$ -  definable in $L_{\td}$, for any $\Lam \in \om$.    By \remref{beth}, 
(using  the $Aut(M)$-invariance), $f \inv(Y)$ is then locally $\bigwedge$-definable in $L$.

In particular, this applies to $f \inv(cl(\la \Sh \ra))$.   
Since $\Lam / f \inv(cl(\la \Sh \ra))$  has   cardinality $\leq 2^{\aleph_0}$,  it is clear that  any definable subset of $\Lam$ containing $f \inv(cl(\la \Sh \ra))$ is commensurable with $\Lam$. 
 
In the other direction, let $H$ be a  locally $\bigwedge$--definable subgroup of $G$ of bounded index.  Then $H$ induces
a  locally $\bigwedge$--definable equivalence relation $E$ on the $G$-torsor $D$,  namely $ y x \inv \in H$; and thus an 
equivalence relation (also denoted $E$) on the type space $S_D(M)$.    It is easy to see that $E$ is   an intersection of
$\L$-definable relations, so that $E$ is closed in the pp topology too; and we have $\rho(p) E \rho(q)$ iff $p E q$. The elements of $Aut(J)$ preserving each $E$-class thus form a closed subgroup of $Aut(J)$.   Identifying
$\CU$ with the image $J$ of $\rho$,  we see that the image of $H$ in $Aut(\CU)$ contains a closed subgroup that in turn
contains $\Sh^{Aut(\CU)}$; so the image contains $cl(\la \Sh \ra)$.  This clearly remains true modulo $\fg$ too.
\eprf

 See \cite{arturo} for a systematic treatment of piecewise hyper-definable groups.      $G^{00}_0$ is defined there to
 be the  smallest locally - $\bigwedge$-definable subgroup of $G$ of bounded index.  In this notation, we have
$f \inv(cl(\la \Sh \ra))= G^{00}_0$.

  \ssec{ Proof of \thmref{grmain} (4)} 
  
   Unlike \propref{G00}, the conclusion here allows for parameters.  Let $M_0$ be an elementary submodel 
  of the given structure, with  
  $|M_0| \leq |\Om|$.  Let $M$ be a highly saturated model of $T^+_{\td}$, with $\aleph_0$-homogeneous reduct
  to $L(M_0)$.  We have the map $\varphi:  Aut(M) \to \CG$.   By \propref{qh3} (1), $\varphi(\Sh^{Aut(M)}) \subset \Sh^{\CG}$.
  Now $Aut(M/M_0) \subset \Sh^{Aut(M)}$, by definition of $\Sh$.  
  Thus if $\si \in Aut(M/M_0)$ and $a \in \varphi \inv(X) $ then $\si \inv a \si \in \varphi \inv(X \Sh^2)$.   Let 
  $X^*,Y^*$ be as in \propref{grmain3}.  Then $X^* \subset \varphi \inv(X \Sh)$ so $\si \inv (X^*) \si \subset \varphi \inv(X \Sh^3)$.
  Since $\Sh Y \meet \Sh ^3 X = \emptyset$,
  we have $\si \inv (X^*) \si \meet Y^* = \emptyset$.   Thus by Beth,  as formulated above, $X^*,Y^*$ are separated by $L(M_0)$-definable sets.     
   
  It follows that the pullback of closed, $\la \Sh \ra$-invariant sets is $\bigwedge$-definable with parameters.  Thus:
  
  \begin{cor}  \label{G00b}  The induced homomorphism $G/f \inv( cl(\la \Sh \ra)) \to \lc/ cl(\la \Sh \ra)$ is   
  continuous, i.e. the inverse image of any compact subset is $\bigwedge$-definable with parameters.
    \end{cor}

\propref{grmain2} applies, in particular, to the top commensurability class of a group $G$, namely the  generic definable subsets.
 Here the results of \cite{patterns} apply directly,   there is no need for locality, 
and $\lc$ is compact; nevertheless, the distance between $\Sh$ and $\lc$ reflects the level of approximateness.
For instance we have, from   \thmref{grmain} (3):

\begin{cor}\label{slowly2}  Let $G$ be a group,  let $X$ be a    symmetric definable 
subset commensurable to $G$.  Let $T^+$ 
  be the associated theory of a $G$-torsor $D$, 
$\CU$ the core space of $T^+$ on the torsor sort $D$, and let $C=Aut(\CU)/\fg$ be the canonical
 compact Hausdorff group associated to $(G,\cdot,X)$.  If $C=1$ then $X^{12} =G$. 
  \end{cor}
 
 It suffices of course that the image of $G$ in $C$ be trivial.   The corollary gives, in principle, a 
  compact model for
slowly converging approximate subgroups $\Lam$, i.e. where $\Lam^m=G$ but the smallest such $m$ is large. 
Whether it can be improved to a Lie model depends, in part, on an effectivization of the theory of normal conjugacy classes
in locally compact groups.  
 
     \begin{question}  How far can the 12 in \thmref{grmain}  be improved?   \end{question}
 
     \begin{question} Describe the image of $f$ in \thmref{grmain}.    Is it equal to $Aut((\CU_D)_h / (\CU_G)_h)$?   \end{question}

 \pb

\section{Approximate subgroups} \label{approxsg}

 Let $G$ be a group.  
 For $\Lam \subset G$, define $\Lam \inv = \{a \inv: a \in \Lam\}$,  $\Lam^2 = \{ab: a,b \in \Lam\}$.
 $\widetilde{\Lam} = \union \Lam^n $ will denote the group generated by $\Lam$.  $\Lam$ is 
 {\em symmetric} if $1 \in \Lam = \Lam \inv$.  
 
    If $Y \subset G$ is covered by finitely many right translates
of $X$,  $Y  \subset  \union_{j=1}^n X b_j$, we will say that $X$  {\em  commensurably covers} $Y$.  
 Two subsets $X,Y \subset G$ are {\em commensurable}, $X \sim Y$,  if each is covered by finitely many translates
of the other:  $X \subset \union_{i=1}^m Y a_i$, $Y \subset  \union_{j=1}^n X b_j$.

Let $\Lam$ be an approximate subgroup of a group $G$ (\defref{ag-tao}.)
  The {\em commensurator} $\comm(\Lam)$ of an approximate subgroup   $\Lam \subset G$ is the group of elements $g \in G$ with $g\Lam g \inv$ commensurable to $\Lam$.    The  {\em  commensurability class $\om$ of $\Lam$} is the 
  set of approximate subgroups commensurable to $\Lam$.  The 
commensurator of $\Lam$ depends only on the commensurability class, and can be denoted $\comm(\om)$ or $\la \om \ra_{max}$.
it is the pullback of the subgroup of $Aut(G)$ fixing $\om$, under the natural map $G \to Inn(G) \leq Aut(G)$.  

 \begin{lem}\label{commlem}  Let   $\Lam $ be an approximate subgroup of $G$, $\om$ the commensurability class.  Then 
 $\la \om \ra_{max} = \union \om$;  $\comm(\Lam)$ is the  union of all approximate subgroups commensurable to $\Lam$.
   \end{lem}
 \prf   Note first that $\Lam \subset \comm(\Lam)$:  if $g \in \Lam $ then $g \inv \Lam g \subset (\Lam)^3$, so 
 $g \inv \Lam g$ is commensurably covered by $\Lam$, and conjugating back by $g$ we obtain the converse.  
 If $\Lam' \sim \Lam$ then $\comm(\Lam')=\comm(\Lam)$, so $\Lam' \subset \comm(\Lam)$.  Conversely, 
 if $g \in \comm(\Lam)$, then  $g\Lam \union \Lam \union \Lam g\inv$ is an approximate subgroup in $\om$;
 see the proof of  \remref{4.3}.
 \eprf

  Let $\lc$ be a locally compact topological group. Let $U$ be a compact neighborhood of $1$.  Then $U$ is an approximate subgroup.    

The compact, symmetric neighborhoods of $1$ form a unique commensurability class of approximate subgroups; we call it the {\em fundamental commensurability class} of $\lc$.      If $f: G \to \lc$ is a homomorphism, 
 the commensurability class in $G$  of the pullback of any compact open neighborhood of the identity in $\lc$ is 
 denoted $\chi(f)$, and again called the fundamental commensurability class of $f$.  

If $H$ is a subgroup of $G$, any commensurability class $\om'$ of approximate subgroups of $H$ extends to a commensurability class $\om$ of approximate subgroups of $G$.  We say in this case that $\om$ {\em belongs to }  $H$.   (This terminology does not preclude belonging to still smaller subgroups.; in fact if $\om$ belongs to $H$, then it belongs to any finite index subgroup of $H$.  In many cases of interest, there does exist a unique smallest commensurability class of subgroups to which $\om$ belongs; see 
\propref{minimalclass}.)

\begin{lem}\label{commlem2}  
Let $H \leq G$.  If  $X,Y \subset H$ are commensurable as subsets of $G$,   they are commensurable in $H$.   \end{lem}
\prf  Say $X \subset \union Ya_i$; then
either $Ya_i$ can be deleted, or $a_i \in Y \inv X \subset H$.   thus $X$ is also contained in a union of $H$-translates of $Y$; and dually.   \eprf

\ssec{Laminar commensurability classes}
A commensurability class  $\om=\chi(f)$ arising from a homomorphism $f$ on a subgroup of $G$ into a locally compact group
$\lc$ will be said to be   {\em {laminar}}.  It is shown in \cite{nqf} that  in this case, $\lc$ can be taken to be a connected  Lie group
with no compact normal subgroups; this Lie group is then uniquely determined by the class $\om$.   
This uses Yamabe's theorem, that guarantees that any locally compact group has an open subgroup 
 compactly isogenous to a Lie group.   This subgroup however is neither canonical nor normal, and making use of it exacts a   technical toll later on,  before one  regains a measure of normality.

In the opposite direction, it is sometimes desirable to extend the domain of $f$ as much as possible.     It follows
from \lemref{commlem} that the domain of any $f$ with $\chi(f)=\om$ is contained in the commensurator group $\om$;
and indeed it is possible to find $f$ with this domain.  

\begin{lem}\label{commsmooth}  If $\om$ is a {laminar} commensurability class of approximate subgroups, there exists 
a homomorphism $F$ with domain $\comm(\om)$  into a locally compact group, with $\chi(F)=\om$.  \end{lem}

\prf This can be seen via \propref{grmain2}.  More directly,  view $\comm(\om)$ as an  ind-definable group, 
where the inductive system is the family of all approximate subgroups in $\om$; by \lemref{commlem},
their union of $\comm(\om)$.   Let $G*$ be a saturated elementary extension, $\lam^*=\lam(G^*)$,
and 
$\om^* = \union_{\lam \in \om} \lam^*$.    Let $I$ be the set of $\lam_1 \in \om$ such that there exist
$ \lam_2 \cdots \in \om$ with $\lam_{k+1}^2 \subset \lam_k$,  $\lam_k \inv = \lam_k$.   If $\lam,\lam' \in \om$,
then we saw that $\lam^2 \meet (\lam')^2 \in \om$.  Hence if $\lam_1,\lam_1' \in I$ and $\lam_k,\lam'_k$
are as in the definition of $I$, since  $\lam_{k+1}^2 \meet (\lam'_{k+1})^2 \in \om$, we have
$\lam_k \meet \lam'_k \in \om$ for each $k$, and so $\lam_1  \meet \lam_1' \in I$.  It is clear that
$I$ is conjugation-invariant.  Let $N$ be the 
  $\bigwedge$-definable subgroup of $\om^*$, formed by the intersection of all $\lam \in I$.   Then  $N$ is a normal subgroup of $\om^*$, $N \subset \lam_1$ (for any choice of $\lam_1 \in I$), and $\om^*/N$ is locally
  compact in the logic topology.  Let $F$ be the natural map $\om^* \to \om^*/N$, restricted to $\om$.      
  \eprf

 For example, the approximate subgroup $\{a/2^n: n \in \Nn, a \in \Zz, -2^n \leq a \leq 2^n \}$ of $\Qq$ has commensurability class $\chi(f)$, where $f: \Zz[1/2] \to \Rr$; and also $\chi(g)$, with $g:  \Qq \to \Rr \times \Pi_{p \neq 2} \Qq_p$ the natural inclusion.

We note immediately a  uniqueness statement for the Lie group explanation of a given commensurability class of approximate subgroups, with a given subgroup as domain.     

 If $N$ is a compact normal subgroup of $\lc$, $\pi: \lc \to \lc/N$ the quotient map, then the composed map $\pi \circ f:G \to \lc/H$ 
 determines the same commensurability class of approximate subgroups as $f$.  
 This is however
 the only ambiguity, i.e. $f: G \to \lc$ is unique up to compact isogeny:

\begin{prop}\label{uniquelc}  Let $h_i: G \to H_i$ be homomorphisms with dense image into locally compact groups, and assume
$\chi(h_1)=\chi(h_2)$.    Then there exist 
compact normal subgroups $N_i$ of $H_i$ and an isomorphism 
$\phi:H_1/N_1 \to H_2/N_2$ , with ${p}_2 h_2 = \phi \circ {p}_1 h_1$ (where ${p}_i: H_i \to H_i/N_i$ is the natural projection.) \end{prop}

\prf     Let $F_0 \leq H_1 \times H_2$ be the image of $G$ under the map $(h_1,h_2): \Lam \to H_1 \times H_2$.
  Let $F$ be the closure of $F_0$.    Then $F$ is a closed subgroup  of $H_1 \times H_2$.   
  Let $\pi_i: F \to H_i$ be the projection,   let $N_1 = \{x \in H_1: (x,1) \in F \}$, $N_2=\{y \in H_2: (1,y) \in F \}$.

\claim{1} $N_i$  is compact; and indeed $\pi_i$ is proper.
  
  \prf  Let $U_i$ be a precompact open neighborhood of $1$ in $H_i$.   Then $h_2 \inv(U_2)$ is commensurable with $\Lam$,
  and hence with $h_1 \inv(U_1)$; so $h_2 \inv(U_2)$ is covered by finitely many translates of  $h_1 \inv(U_1)$; hence
  $h_1  h_2 \inv(U_2)$ is covered by finitely many translates of  $U_1$, so it is precompact.

  So     $N_1 \subset \pi_1(F \meet (H_1 \times U_2)) \subset
\pi_1  cl(F_0 \meet  (H_1 \times U_2)) \subset  cl(h_1  h_2 \inv(U_2))  $ and so $N_1$ is precompact.  Being closed it is compact, and likewise $N_2$. 
Along the way we showed properness:  $\pi_2 \inv( U_2) \subset  cl(h_1  h_2 \inv(U_2)) \times U_2$ so it is pre-compact.
  
  \eprf
  
  \claim{2}    $F$ projects onto $H_1$ and onto $H_2$.   
  \prf  Let $U$ be an open subset of $H_2$, with compact closure $\bU$.  By Claim 1, $\pi_2$ is proper, 
  so $\pi_2 \inv(\bU)$ is compact.   Thus $\pi_2(\pi_2 \inv(\bU))$ is compact.    By the density assumption, 
  $\pi_2(\pi_2 \inv(U))$ is dense in $U$.  Thus $\pi_2(\pi_2 \inv(\bU))$ is closed, and contains a dense subset of $U$,
  hence it contains $U$.  It follows that $F$ projects onto $H_2$, and similarly $H_1$.
  \eprf

Now $N_i \normal H_i$, and $F$ induces a group isomorphism $\phi: H_1/N_1 \to H_2/N_2$.
By the same argument as in   Claim 1, $\phi (X)$ is compact iff $X$ is compact.  
   Let $V_i$ be a compact neighborhood of $1$ in $H_i$, and $W_1 = V_1 \union \phi \inv(V_2)$.   Then $W_1$ is compact, as is $W_2 = \phi(W_1) = \phi(V_1) \union V_2$. 
Thus  $\phi$ restricted to $W_1$ is a homeomorphism near $1$, which suffices.

\eprf

There is   a  one-sided  version of \propref{uniquelc}:    
\begin{rem}  If $h: G \to H$ and $h': G \to H'$ are homomorphisms with dense image into locally compact groups, 
and $\chi(h)$ lies below $\chi(h')$ in the sense that there exist $\Lam \in \chi(h), \Lam' \in \chi(h')$ with $\Lam \leq \Lam'$,
then there exists a compact normal subgroup $N'$ of $H'$ and a homomorphism $\phi: H \to H'/N'$ such that
$\phi \circ h = p' \circ h'$.     \rm 

This  could be proved in the same way, but also follows from the statement of \propref{uniquelc}:   define $h_2: G \to H \times H'$
by $h_2(x)=(h(x),h'(x))$, and let $H_2$ be the closure of the image of $h_2$.  Then $\chi(h_2)=\chi(h)$ and \propref{uniquelc} applies, and proves the remark.  \end{rem}

Given a commensurability class   $\theta$ of approximate subgroups, if it arises from a homomorphism  $f: G \to \lc$, there  is always a {\em universal} one with class $\theta$, mapping with compact kernel onto any other.
Often, there is also a minimal one, namely $f: G \to \lc$ such that $\lc$ has  no compact normal subgroups.
Both are unique up to a unique isomorphism under $G$.

 \ssec{Intersections}   \label{icc}  (See \cite{machado3}).
 Let $\om$ be a 
 commensurability class of approximate subgroups of a group $G$, and  let $H$ be a subgroup of $G$.  
 If $X ,Y \in \om$, then $Y \meet H$ is commensurably covered by $XX\inv  \meet H$.  To see this, say $Y \subset \union_{i=1}^m Xc_i$.  Then $Y \meet H \subset \union_{i=1}^m (Xc_i \meet H)$.  Now if $Xc_i \meet H \neq \emptyset$, 
 let $h_i \in Xc_i \meet H$.   Then $c_i \in X\inv h_i  $ so $Xc_i \subset X X \inv h_i$.   Hence the $H$-commensurability class of 
 $XX\inv \meet H$ does not depend on the choice of $X \in \om$.  In particular,
 as $(XX\inv \meet H)^2 \subset (XX\inv)^2 \meet H $ and $ (XX\inv)^2 \in \om$, $XX\inv \meet H$ is an approximate subgroup of $H$.     We denote this class by $\om \meet H$.   
 If $H$ is normal,
 we also have the class $\om /H$, represented by $X/H$ for any $X \in \om$.  
 Then $\om$ belongs to $H$ iff  $\om \meet H \sim \om$.    When $H$ is normal, this is equivalent to:  $\om/H \sim \{1\}/H$.

 One can similarly define an intersection of commensurability classes:  if $X,Y$ are approximate subgroups of $G$,
 then  so is $XX \inv \meet Y Y \inv$, and the commensurability class of $XX \inv \meet Y Y \inv$ depends only on the classes of $X$ and of $Y$.   In case $X,Y$ are commensurable, $X X \inv \meet Y Y \inv$ is commensurable to both.  This can be seen directly using a Rusza-style argument, starting from a maximal subset $\{a_1,\ldots,a_n\}$ of $X$ such that the sets $a_iY$ are disjoint.      A more conceptual derivation will soon be available:  we will see below (\thmref{grmain})  that any commensurability class of approximate subgroups can be written as $\chi(f)$, for some approximate homomorphism into a locally compact group; then it is easy to check that   $\chi(f \times g)  = \chi(f) \meet \chi(g)$, and $\chi(f \times f) = \chi(f)$.  

\ssec{Minimal commensurability class of supporting subgroups}
\propref{minimalclass} was obtained only for laminar approximate lattices in an earlier version of this paper, as a corollary of \propref{sensible}.   Thanks to Simon Machado for discussions regarding the general case below.

 \begin{lem}\label{minimalclass2lem}     Let $\om = \chi(g)$, $g: G \to \Rr^n:[-1,1]^n$ a homogeneous quasimorphism.   
  Let $S \leq S' $ be   subgroups of $G$ to which $\om$ belongs, such that $g(S),g(S')$  span the same $\Rr$- subspace of $\Rr^n$.    Then $S$ has finite index in $S'$.  \end{lem}
 
 \prf   
Let $s_1,\ldots,s_n \in S$ be such that the $g(s_i)$ form a basis $v_1,\ldots,v_k$ for the span $V$ of $f(S')$.  
 Now any point of $V$ is at bounded distance from a lattice point $\sum m_i v_i$; and $\sum m_i v_i$ is at bounded distance from $g(s_1^{m_1} \cdots s_n^{m_n})$, as $g$ is a  homogeneous  
 quasimorphism.  Thus we can find a ball $B$ in $V$ such that $f(S)+B = V$, and $B$ contains the error set of $g$.
 Let $x \in S'$.  Then for some $s \in S$ we have $g(s) \in g(x) + B$ and so  $g(x s \inv)  \in BB$.   Thus $g(x) \in g \inv(BB) S$.  Now $g \inv(BB)$ is covered by finitely many
 cosets of $S$.  Hence $G$ is covered by finitely many
 cosets of $SS=S$, i.e. $[S':S]<\infty$.  
 \eprf

\begin{prop} \label{minimalclass} For any commensurability class $\om$ of approximate subgroups of a group $G$, there exists a unique  
  commensurability class ${\la \om \ra} _{min}$ of subgroups   of $G$ such that    $\om$ belongs to a subgroup $S$ iff
  $S$ contains some element of ${\la \om \ra} _{min}$.    \end{prop}
  
    We will use several times the statements on intersections of approximate subgroups above.  
 In particular,  the class $\Om$ of subgroups containing (commensurably) an element of $\om$ is closed under finite intersections.   
  

%

%
%
  

We will   
use \thmref{ag2}  in the version of supplement (i).   We may assume $G=\wl$,
so we have   a connected Lie group $L$,  a rigid abelian subgroup $A \normal L$,  $A \cong \Rr^n$, 
a normal compact subset $K$ of   $A$, 
and a    quasi-homomorphism 
$f: G  \to L:K$
with $\om= \chi(f)$.  Let $\bar{f}$ be the homomorphism obtained as a composition of $f$ with the quotient map $L \to L/A$.

Let $H$ be an element of $\Om$ such that $L_H= cl( \bar{f}(H))$ has least possible dimension.   
The connected component 
$L_H^o$ of $L_H$ is a connected Lie group of the same dimension, and is open in $L_H$; so for a compact $U \subset L/A$,
$U \meet L_H$ is contained in finitely many cosets of $L_H^o$.  Thus if $X \in \om$ and $X \subset H$,
then $X$ is contained in finitely many translates of $H':=\bar{f} \inv (L_H^o)$.    So $H' \in \Om$; and $\bar{f}(H')$ is
dense in $L_H^o$.  Changing notation, and replacing $H$ by $H'$,  this shows that we can find $H$ with  $cl(\bar{f}(H))$ 
connected; we may change notation to    $cl(\bar{f}(H)) = L/A$.   By minimality of the dimension, for any $S \leq H$ with $S \in \Om$, we have $cl(\bar{f}(S)) = L/A$.  
 
$G$ is considered
   as a piecewise-definable group in a structure $M$ (where a piece is an element of $\om$).  We can take    all subgroups   $S$ of $G$ to be locally definable.   We take  a saturated extension $M^*$ of $M$,
   and write $S^*$ for $S(M^*)$.      
    By local compactness
   of $L$ and $L/A$,  $f$ and $\bar{f}$ extends to $G(M^*)$.   
         Let $\G =\ker(\bar{f})= f\inv(A)$.   So $\G$ is a  $\bigwedge$-definable group.  
         Note that if  $cl(\bar{f}(S)) = L/A$, we have $\bar{f}(S(M^*))=L/A$.
   
   

Since $\Om$ is closed under finite intersections, to find a minimal element (up to commensurability) it suffices to show
that there is no infinite descending chain of subgroups $H>H_1 > H_2> \cdots$ in  $\Om$,
and with $H_i/H_{i+1}$ infinite.  We show more:
   
 \claim{}  Let    $H > H_1 > \ldots >H_{n+1}$ be elements of $\Om$.  Then $[H_i: H_{i+1}] $ is finite for some $i$.
 
  \prf   
  Let $g$ be the homogeneous quasimorphism at bounded distance from $f | \G$.  Note that $\chi(g)$
  belongs to  $H_i^* \meet \G$ for each $i$ (a bounded subset of $A$ is still a bounded subset of $L$.)  
  The spans of $g(H_i^*)$ form a descending sequence of subspaces, so $g(H_i^* \meet \G)$ has the same span as $g(H_{i+1}^* \meet \G)$ for some $i$.
By  \lemref{minimalclass2lem}, $[\G \meet H_i^*: \G \meet H_{i+1}^*] $ is finite. 

So $\G \meet H_i^* = ( \G \meet H_{i+1}^*) F$ for some finite $F \subset H_i^* \G$.

Now let $x \in H_i^*$.  Let $y \in H_{i+1}^*$ be such that $\bar{f}(y) = \bar{f}(x)$.  Then $y \inv x   \in H_i ^* \meet \G$.
So $   y \inv x   \in H_{i+1}^* F$ and thus $x \in H_{i+1}^*F$.  This shows that $H_i^* \subset H_{i+1}^*F$.
So $H_i^*$ does not intersect $|F|+1$ distinct cosets of $H_{i+1}^*$.  Thus $H_i$ does not intersect $|F|+1$ distinct cosets of $H_{i+1}$; so $[H_i:H_{i+1}] \leq |F|$.   This proves the claim and hence the Proposition.

 \eprf

\ssec{The fundamental commensurability class}
\label{classofmorphism}

More generally,   if $f: G \to \lc:K$ is a quasi-homomorphism, with  $\lc$ locally compact and 
$K=K \inv$ compact, we consider pullbacks $Y=f \inv(W )$ of  compact neighborhoods $W$ of $K$ in $\lc$.
 We will see below that
   $Y Y \inv$ is an approximate subgroup, whose 
commensurability  class is uniquely determined.  We continue to use   the notation $\chi(f)$ for this class, and to refer to it as the fundamental class of $f$.  
We will write $X \sim \chi(f)$ to mean that the commensurability class of $X$ equals $\chi(f)$. 
Uniqueness again holds but now in the category of compact correspondences.

\begin{prop}\label{converse}  Let $H$ be a   locally compact topological group, $K$ a symmetric compact subset of $H$,
and assume \[f: \Lam \to \lc:K \]
is a quasi-homomorphism, with 
\[f(x y \inv) \in f(x) f(y) \inv {K} \] 
Then all sets $X=f \inv({U}{K})$, with   ${U}$   a compact neighborhood of the identity in  $\lc$,  are
 right- commensurable in $\Lam$.     Thus any symmetric set sandwiched between two such is an approximate subgroup.
 In particular   $XX \inv$ is one. 
\end{prop}

\prf    We begin by comparing translation in $\Lam$ and in $\lc$.

\claim{1}   For any ${V} \subset \lc$ and $a \in \lc$,    $f \inv (  {V}a) \subset   f\inv({V}{V}\inv {K}) c$ for some $c$.   
 
 If $f \inv ({V}a) = \emptyset$, there is nothing to prove.  Otherwise, let $c$ be such that $f(c) \in {V}a$.
If $f(x) \in {V}a$, then $f(x    c \inv) \in  ({V}a)  ({V}a) \inv {K} = {V}{V} \inv {K}$; so $xc\inv \in f\inv({V}{V} \inv {K})$ and 
$x \in  f\inv({V}{V} \inv {K}) c$.

\claim{2} Let $Z \subset \lc$ be   compact and let $U$ be a neighborhood of $1$ in $\lc$.  Then
 finitely many right cosets of $f \inv (U  {K})$ cover $f \inv (Z)$.   
 
 \prf Let $V$ be a neighborhood of $1$, with
 $VV\inv \subset U  $.    As $Z$ is compact,
 it is covered by finitely many cosets $Va$;
 so it suffices to show that $f \inv (Va)$ is covered
 by a right coset of $f \inv (V V \inv {K})$; this follows from (2) above.    \eprf

 Let $U$ be a compact neighborhood of $1$ in $\lc$, $X= f\inv(UK)$ and $Y = X X \inv$.      Note $Y$ is symmetric ($1 \in Y=Y \inv$).
 We have $Y \subset f \inv(UK(UK) \inv K)$; taking $Z=UK  (UK) \inv K$ in Claim 2, we see that $Y$ is covered
 by finitely many right translates of $X$.   
 Claim 2 also  implies the right-commensurability of any two such $X$, and hence of any set sandwiched $W$ between two such.
 Thus any such $W$ is an approximate subgroup; in particular this applies to $Y$.  
\eprf

\begin{rem} \begin{enumerate}
\item    \propref{converse} did not assume normality of $K$.   
\item  If we do assume that $K$ is normal and  that  $f(xy) \in Kf(x)f(y)$ (as is the usual assumption for quasimorphisms), 
we obtain   (with $x,y$ replaced by $x,y \inv x$) that $f(x) \in  f(y)f(y \inv x) K$,
and so changing sides and using $K=K \inv$, we have $f(y \inv x) \in f(y) \inv f(x) K$, and similarly $f(x y \inv) \in f(x) f(y) \inv K$.   
\item  
It is worth noting that although we began with an approximate homomorphism into a locally compact group in  the general sense of a {\em compact} error set $K$, the associated commensurability class in $\Lam$  consists of   approximate subgroups in
Tao's original sense of a {\em finite} set of translates.  
 \item   Responding to  a question of Kobi Peterzil during a talk, we consider the case of quasimorphisms, i.e.  take  
 $\lc=\Rr$.    we may take $K=[-1,1]$.   Modifying
 $f$ to be homogeneous, it becomes symmetric, $f(x \inv)=-f(x)  $ and $f(1)=0$.  
 Given $e>0$, let $k(e)$ be the smallest natural number  such that $X=f \inv(-1-e,1+e)$  is a $k$-approximate subgroup
 (\propref{converse}).  What kind of function is $k(e)$?

 Following the above proof, we have 
$XX \subset f \inv(-3 -2{e},3+2 {e})$;   this interval  is covered by  $m$ translates of $( -{e}/2, {e}/2)$, if $m> 4(1+  6/e)$;
the pullback of each such translate   is contained in a translate of $X$.    Hence $X$ is a $5$-approximate subgroup if $e>24$.

\begin{question} Is  there an absolute bound on the level of approximateness,   for small $e>0$?\end{question}

 In \exref{milnor}, with $\Lam$ the central extension of the homeomorphism group of the circle,  the  
 translation number $f$ commutes with translation by   $\Zz$;   for any $e>0$, if $X=  f \inv(-1-e,1+e)$, then
 $XX \subset f \inv(-3-2e,3+2e) \subset X \pm \{0,1,2,3\}$.     What about Brooks' examples,  \cite{kotschick}, Example 1?

 \begin{question}  What is the minimal $k$ such that the 
approximate kernel of any quasimorphism is commensurable to     a $k$-approximate subgroup?   \end{question}
 By the above discussion, the answer is one of the numbers 2,3,4,5.



 \end{enumerate}\end{rem}
 
 We have already proved the principal structure theorem for approximate subgroups, \thmref{grmain}; we will now refine it.  
See \defref{groupdefs} (6) for the definition of rigidity, and \secref{classofmorphism} for $\chi$.
 
\begin{thm} \label{ag2}  Let $\om$  be a commensurability class of approximate subgroups of a group $G$.  Then there exists  a subgroup $\wl \leq G$, a connected Lie group $\sl$,  a rigid abelian subgroup $A \normal \sl$,  a normal compact subset $K$ of   $A$, 
and a    quasi-homomorphism 
$f: \wl  \to \sl:K$
with $\om= \chi(f)$.   
\end{thm}

Some supplements:  \begin{enumerate}[label=(\roman*)]
 \item    One can demand that $A \cong \Rr^N$
for some $N \in \Nn$,  and that  $\la f(\wl) \ra A$ be dense in $\sl$.     
\item  One can add that $\sl$ has no compact normal subgroups. 
\item In  place of connectedness of $\sl$,  and as an alternative to (i), one can ask that   $\la f(\wl) \ra $ is dense in $\sl$,
and $A \cong \Rr^n \times \Zz^k$ for some $n,k$.
\item  The induced homomorphism $\wl \to \sl/A$ is definable in the sense of continuous logic.
 
\end{enumerate}

\medskip

Following the proof, we will give two ways to decompose the quasi-homomorphism  $f$ into  better studied  kinds.
In \thmref{ag2mt}, we will describe the situation in terms of a $\bigwedge$-definable group $H$ and $n$ definable quasimorphisms
on $H$, i.e. a definable quasi-homomorphism on $H$ into $\Rr^n$.   In \propref{ag2gt},  we will instead decompose
$f$ in terms of a rigid group extension $\widehat{\Lam}$ of $\wl$, and a group homomorphism $\hat{f}: \widehat{\Lam} \to \sl$.


%

\ssec{Refining the target.}  \label{mprov}  
The proof of   \thmref{ag2} will start with a quasi-homomorpism ${f}: \al \to \lc:K $ as in \thmref{grmain},
and gradually seek to simplify it, without changing the commensurability class $\chi({f})$  associated with ${f}$.    
Let us review in advance a few of the allowed modifications, along with an explanation of why the  pullback class $\chi(f)$    does not change.   Some of the steps parallel  \cite{FK}, where Fujiwara and Kapovich  proved that  any quasi-homomorphism with discrete target can be perturbed to one with central error.     Here the best we can aim for is a rigid error set.

\begin{enumerate}

\item  Replace $\lc$ by an open subgroup $\lc_1$ containing $K$; let ${\al}_1= {f} \inv(\lc_0)$, ${f}_1 = {f}|{\al}_1$; $K_1 = K \meet \lc_1$.   
Then ${\al}_1$ is a subgroup (if $x,y \in {\al}_1$ then ${f}(xy) \in {f}(x){f}(y) K \subset {\al}_1$.)
Since a compact open neighborhood of the identity in $\lc_1$ is also one in $\lc$, we have $\chi({f}|{\al}_1) = \chi({f})$.

\item 
Factor out a compact normal $N$ subgroup of $\lc$:  let ${f'}: \al \to \lc/N: K/N$ be the composition of ${f}$
with the quotient map $\pi:  \lc \to \lc/N$.    We have $({f'}) \inv (YN/N)  = {f} \inv (Y N) $ for any 
compact neighborhood  $Y$ of  $K$ in $\lc$, so $\chi({f})=\chi({f'})$.  

 \item Let $\lc_2$ be a  closed subgroup of   $\lc$, with $\lc_2 K^m = \lc$.   Let $K_2 = \lc_2 \meet K^{3m+1}$.  
 For any $x \in {\al}$ pick $c(x) \in K^m$ with $f(x)c(x) \in \lc_2$, $c(1)=1$, and define $g:{\al} \to \lc_2$ by $g(x)=f(x)c(x)$.  
 For any $x,y \in {\al}$ 
we have, for some $k=k(x,y) \in K$:
\[g(xy) = f(xy)c(xy) = f(x)f(y) k c(xy) = g(x) c(x) \inv g(y) c(y) \inv k c(xy) =  \]
\[=g(x) g(y) c(x)^{-g(y)} c(y) \inv k c(xy)  \in  g(x)g(y) K^{3m+1} \] so $g(xy) \in g(x)g(y) K_2$.  

For a compact neighborhood $U$ of $1$ in $\lc$,   $UK^m \meet \lc_2$ is a compact neighborhood of $1$ in $\lc_2$, 
and $f \inv(U) \subset g \inv(UK^m \meet \lc_2) \subset f \inv (U K^{2m})$.  Thus $\chi(f)=\chi(g)$.

\item  Let $H_0$ be an open subgroup of $\lc$ of finite index.  Let $\lc_3 = H_0 \la K \ra$.  Then $\lc_3 = H_0 K^m$
for some $m$.  Let ${\al}_1 = f \inv(\lc_3)$. 
  Then ${\al}_1$ is a finite index   subgroup of ${\al}$, and by (1) $\chi(f) = \chi(f|{\al}_1)$.
By the previous step (3), we can find $g: {\al}_1 \to H_0$ with $\chi(f) = \chi(g)$.

  \end{enumerate}

     \begin{rem}  To a model theorist,  the apparent   arbitrariness in the choice of $f_0$ in step  (3) is disturbing.  
     It will only be used via (4), so in a finite index setting, where only finitely many parameters need be chosen.  Note also
  in general
     that if   $\la K \ra$ is centerless, the correction function $c_0$, and hence $f_0$, are uniquely defined in the construction.   
 The general case could be made definable  using a relational presentation of the notion of a quasi-homomorphism, 
 namely the category $\cogr$ of Appendix \ref{categories}, 
 in effect allowing $f(x)$ to be a compact subset rather than an element; 
 and taking $k(x) = f(x)\inv H_2 \meet K^m$. 
\end{rem}

\prf[Proof of \thmref{ag2}]
Let $\Lam_0$  be an approximate subgroup in $\om$, and $\wl = \la \Lam_0 \ra$.  
By \thmref{grmain},  
 there exists $f: \wl \to \lc:K$ with $\lc$ locally compact,
$K$ normal, symmetric, compact,with $\Lam_0 \sim \chi(f)$.

By \propref{lc1cc}, there exist closed normal subgroups $N \leq A \leq H_1 \leq \lc$ with $\lc/H_1$ finite, 
$N$ compact,  with $A/N$  basic abelian and rigid in $H_1/N$,    such that $\la K \ra \meet H_1  \leq A$.
By \ref{mprov} (4), we may pass from $\lc$ to  the finite index subgroup $H_1$ (replacing $K$ by $K^m \meet H_1$ for appropriate $m$);  so we may assume $K \subset A$, $N \leq A \leq \lc$.  Since $N$ is compact, by \ref{mprov} (2) we may factor it out.
Now $A$ is  basic abelian and rigid in $\lc$, and $K \subset A$.  

  We have $A=DV$ with $D \cong \Zz^l$ central in $\lc$, and $V$ a rigid vector group.  Let $D_\Rr = D \tensor \Rr$, so that
  $(A,D) \cong (\Zz^l,\Rr^l)$.   Let $\bar{\lc} =  \lc \times_D D_\Rr  := (\lc \times D_\Rr) / D$, with $D$ embedded into
  $\lc \times D_\Rr$ by $d \mapsto (d,-d)$.   This diagonal  image of $D$ is a closed, discrete subgroup of $\lc \times D_\Rr$
  and so the quotient is still locally compact.  Now $\lc$ embeds into $\bar{\lc} $ as a closed subgroup. Let $\bar{A}=A \times_D D_\Rr$.  
If $U$ is a compact neighborhood of $1$ in $\bar{\lc} $, then $U \meet \lc$ is a compact neighborhood of $1$ in $\lc$.  
Hence we may replace $\lc,A,K$ by $\bar{\lc} ,\bar{A},K$; we have the same situation, and in addition $\bar{A} \cong \Rr^N$ for some $N$.   We are of course not claiming any density at this point.

  By Yamabe's theorem \cite{yamabe}, $\lc$ has an open subgroup $H'$ and a  compact subgroup $N'$, normal in $H'$,
 with $H'/N'$ a 
 Lie group.    
  As $V$ is connected and $H'$ is open,  $V \subset H'$; i.e. by the above paragraph, $A \subset H'$.  By 
   \ref{mprov} (1), we may pass from $\lc$ to $H'$; so we may assume $\lc=H'$.
 Factoring out the compact normal subgroup $N'$ using  \ref{mprov} (2)
 again, we have that $\sl$ is a Lie group. 

  Replacing $\lc$ by the closed subgroup
$\sl=cl(f(\wl)A)$, we may assume  $f(\wl) A $ is dense in $\sl$.     Now $\sl^0$ is an open subgroup of $\sl$; 
$\sl^0 \meet f(\wl)A$ is thus still dense in $\sl^0$: since $A$ is  connected, $A \subset \sl^0$; by  \ref{mprov} (1), we may replace $\sl$ by $\sl^0$, so that we may assume $\sl$ is connected.   Finally, by \ref{mprov}(2), we may, without changing $\chi(f)$, factor out the maximal normal compact subgroup of $\sl$, to obtain that $\sl$ has no nontrivial normal compact subgroups.
 
Now the continuity (iv) comes from   \corref{G00b}; the various modifications made
to the quasi-homomorphism of \thmref{grmain}, to wit restricting the domain to a  subgroup, factoring out a compact normal subgroup of the image, and modifying by a map into $A$,  will not effect this.

\eprf

\begin{rem} \label{ag2r} 
  If desired, say in order to connect with cohomological results phrased for unitary modules, 
   we can make the error set contained in a complex (unitary) rather than real rigid space.   
Take first the  semi-direct product of $\sl$ with the complexifization $A_{\Cc}$ (where $\sl$ acts on $A$ by conjugation, and 
using the derived action on $A_{\Cc}$),  and then factoring out the image of the  homomorphism $A \to \sl \times A_{\Cc}$,
$x \mapsto (-x,x)$,  which becomes a normal subgroup.    The three variants $\sl,\sl',\sl''$ do not detract from the canonicity of the construction, but give three canonical forms that may be convenient for different purposes.
\end{rem}


\ssec{Model-theoretic formulation}

Let us recall some definitions.  As usual, we   interpret  $\bigwedge$ and $\bigvee$-definable
sets in a sufficiently saturated model of the theory. 

\begin{enumerate}

  \item Let $\wl = \union H_n$ be a $\bigvee$-definable group.  
   A set $X \subset \wl$ is {\em locally $\bigwedge$-definable} if for each $N$, $H_n \meet X$ is $\bigwedge$-definable. 
 \item A subgroup $\G$ of $\wl$ is {\em locally $\bigwedge$-definable} if $H_n \meet \G$ is $\bigwedge$-definable, for all $n$.
 
   
        \item  A $\bigwedge$-definable approximate subgroup $W$ 
   is said to {\em determine a commensurability class $\omega$} if $W= \meet_n W_n$, with $W_1 \supset W_2 \supset \cdots $
   all in $\om$.   
   
\item Let $X$ be a    $\bigwedge$-definable set.   A function $\phi: X \to \Rr$   is said to be {\em continuous} (or: {\em definable in the sense of continuous logic}) if  it is bounded, and for any two disjoint compact sets $C,C' \subset \Rr$ there exists a definable $D$
with $\phi \inv(C) \subset D$ and $D \meet  \phi \inv C' = \emptyset$.  
\item Let $X$ be a locally   $\bigwedge$-definable subset of $\wl$.  Then $\phi: X \to \Rr$ is {\em continuous} if 
$\phi | (X \meet H_n)$ is continuous for each $n$.   
 \end{enumerate}
 
 %
  
%
%
%

\begin{thm} \label{ag2mt}  Let $G$ be a group,  $\om$    a commensurability class of approximate subgroups of a group $G$.
 Then there exists a $\bigvee$ -definable  
group $\wl$, 
a locally $\bigwedge$-definable subgroup $\G  \leq \wl$ of bounded index, and finitely many continuous
  quasimorphisms  $\phi_1,\ldots,\phi_n :  \G  \to \Rr$,  such that the 
   approximate kernel \[ \G_\phi : = \meet_{i=1}^n \phi_i \inv [-1,1] \]
   is a $\bigwedge$-definable approximate subgroup, and
  $\omega$ is the commensurability class determined by $\G_\phi$.
\end{thm} 

\prf 
 
 By \thmref{ag2},  there exists a  quasi-homomorphism $f: \wl  \to \sl:K$, where $K$ is a compact subset of   a rigid abelian subgroup $A \normal \sl$, $A \cong \Rr^n$,   $\om=\chi(f)$.
    Let $\bar{\sl} = \sl / A$,   let $\bar{f}:  \wl \to \bar{\sl}$ be the composed homomorphism $\wl \to \sl \to \bar{\sl}$;
    and let $\G$ be the kernel.  
   By \corref{G00b}, since $A$ contains $cl(\la \Sh \ra)$,  the induced homomorphism $\wl/ \G  \to \sl / A$ is continuous, and $\G$ is a 
 locally $\bigwedge$-definable group.  
 
 Let $\psi = f| \G$.  
 Identifying $A$ with $\Rr^n$ and  projecting
 to coordinates, we can  view $\psi$ as an $n$-tuple of quasimorphisms $\psi_1,\ldots,\psi_n$.  
  Let $\phi_i$ be the homogeneous
 quasimorphism at bounded distance from $\psi_i$; and let $\phi=(\phi_1,\ldots,\phi_n)$.  
 We have:
\begin{enumerate}
\item  $\phi \inv (-m,m)^n$ is an approximate subgroup, for large enough $m$ (by  \propref{converse}).
 \item   If $C \subset A$ is compact, then $\psi \inv(C) \subset f \inv(C)$
 so this set is contained in a definable subset of $\wl$.  (\thmref{grmain} (4)).
 \item   for some $m$,  the pullbacks under $\psi$ of two compact subsets of $\Rr^n$ whose Euclidean distance is at most $m$ are separated by  a definable set
\item $\psi \inv (-m,m)^n$  has bounded index in $\wl$ (even if we work in a large saturated $M$, if $(a_i: i < (2^{\aleph_0})^+$
are elements of $\wl$, then some $a_i,a_j$ with $i<j$ must have the same image in $\sl$, hence $a_i \inv a_j \in A$,
and $\psi(a_i) \in \psi(a_j) \Sh \subset (-m,m)^n$.) )
\end{enumerate} 
Properties  (2-4) of $\psi$ above are all also inherited by $\phi$.

 
Replacing $\phi$ by $\frac{1}{m} \phi$   we may assume $\phi \inv ((-1,1)^n)$ is an approximate subgroup, and the inverse images of two compacts that are $1$-separated (i.e. $[-1,1]^n C \meet C' = \emptyset$) are separated by a definable set.    
Along with the homogeneity of $\phi$, this implies that $\phi$ is  continuous:  if $C ,C'$ are disjoint compact subsets of $\Rr^n$,
they are $\e$-separated for some $\e >0$; let $N$ be an integer with $N \e >1$.  Note 
 $\phi \inv(NC), \phi \inv(NC')$ are contained in a definable subset of $\wl$.  As   $NC, NC'$ are $1$-separated,
  there exists a definable set $D$ with $\phi \inv(NC) \subset D , D \meet \phi \inv(NC')=\emptyset$.  Let
  $D^*=\{g \in \G:  g^n \in D\}$.  Then $\phi \inv(C) \subset D^*$, $D^* \meet \phi \inv (C') =\emptyset$.   Thus the definable sets
  separate pullbacks of an arbitrary pair of disjoint compacts, proving the continuity of $\phi$.
Since $\phi \inv [-1,1]^n$ has bounded index in $\wl$, any definable set containing  $\phi \inv [-1,1]^n$ must be commensurable with the definable approximations to $\wl$.
 
%
 
 \eprf

\begin{rems} \label{canchooseLie}
\begin{enumerate}
\item    In  case $n=0$, $\G$ is $\bigwedge \om$-definable.    Existence of such a $\G$ is assured in this case by Theorem 3.5 of  
  \cite{nqf}. 
 \item    We may choose $\wl,\G$ so that $\wl/\G$ is a connected
  Lie group with no nontrivial compact normal subgroup.  

 \end{enumerate}

\end{rems}

  \begin{cor}  Let $G$ be a group of bounded exponent.  Then any approximate subgroup of $G$ is commensurable to a subgroup.
  \end{cor}
  
 \prf Let notation be as in \thmref{ag2mt}; by \remref{canchooseLie}, we may take $\wl/\G$ to be a connected
  Lie group.  Since it has bounded exponent, it must be trivial.  Thus by the Theorem, $\omega$ is
   the commensurability class associated with a quasimorphism $\wl \to \Rr^n$; we can take the quasimorphism to be homogeneous.
   In this case, the bounded exponent again implies that the image is $0$.  So $\omega$ is the commensurability class of the group
   $\wl$.  \eprf

\ssec{Group-theoretic formulation} \label{ag2gt}

We   restate \thmref{ag2}  in more standard group-theoretic terms, somewhat in the spirit
  of the Bargmann-Wigner theorem 
  on projective representations of Lie groups.
  

  Let $G$ be a group, let $E_n$ be Euclidean space, $O_n$ the group of isometries of $E_n$ fixing $0$.  Assume
a homomorphism $\rho:G \to O_n$ is given; $\rho$ makes $E_n$ into a $G$-module $V$; we call a {\em rigid} $G$-module
to indicate that the image of $\rho$ preserves the Euclidean norm.

Let $C^n(G,\rho)=C^n(G,V)$ be the $\Rr$-space of functions from $G^n \to V$ view this as a complex with the usual differential, 
$ \left(d^{n+1}\varphi \right)(g_{1},\ldots ,g_{n+1})=g_{1}\varphi (g_{2},\dots ,g_{n+1})+\sum _{i=1}^{n}(-1)^{i}\varphi \left(g_{1},\ldots ,g_{i-1},g_{i}g_{i+1},\ldots ,g_{n+1}\right)+(-1)^{n+1}\varphi (g_{1},\ldots ,g_{n}) $.  

The kernel of $d^{n+1}$ is denoted
$Z^n$, the group of $n$-cocycles.  In particular any $\alpha \in Z^2(G,V)$ determines canonically
a group $G \ta V$ with universe $G \times V$, and multiplication:  $(g,e) (h,d) = (gh,  \rho(h) e +d+\alpha(g,h))$.  It fits into  an exact sequence
\[  0 \to E_n \to G \ta V \to G \to 1 \]


In case $\rho$ is trivial, this is a central extension of $G$ by $\Rr^n$; in general we may call it a rigid extension of $G$.

When $G$ is discrete, we give $ G \ta V$ the product topology of $G \times \Rr^n$.
  Let $C_{bdd}$ be the subcomplex   of the bounded chains,
  and $H^*_{bdd}$ or $H^*_b$ the cohomology.   

%
%
%
%
%

%

%
%
%


  \def\hG{\tilde{G}} 
%
 
\begin{thm} \label{ag2gr}  Let $\om$  be a commensurability class of approximate subgroups of a group $G$.  Then there exists  a subgroup $\wl \leq G$, a connected Lie group $\sl$ without compact normal subgroups,  
 a rigid $\wl$-module $V$, 
 a cocycle bounded $2$-cocycle $\alpha: \wl^2 \to V$, 
 and a continuous homomorphism $\phi:  \wl \ta V \to \sl$,  with $\om$
commensurable to $\{ g \in \wl:  \phi(g,0) \in U \}$, for $U$ a compact open neighborhood of $1$ in $L$.
   \end{thm}
   
 Above, continuity refers to the product topology with $\wl$ discrete; it amounts to  continuity on $1 \times V$.


  \prf  By \thmref{ag2}  there exists  a subgroup $\wl \leq G$, a connected Lie group $\sl$, 
a normal compact subset $K$ of a rigid abelian subgroup $A \leq \sl$,  $A \cong \Rr^N$,
and a   quasi-homomorphism 
\[f: \wl  \to \sl:K\] 
with $\om = \chi(f)$.       We will identify $A$ with $\Rr^N$, and write multiplication additively within $A$.
 Rigidity of $A$ means that there exists a homomorphism $\rho': \sl \to O_N$ with $x \inv a x = \rho'(x)a$,
 for all $x \in \sl, a \in A$.  
 For $g \in \wl$, let $\rho(g) = \rho'(f(g))$.  As $K \subset A$, $\phi(K) =1$ and hence $\rho:  \wl \to O_N$ is a homomorphism.
Let  $V$ be $A$ with this $\wl$-module structure.
 Define $\alpha: \wl^2 \to A=\Rr^N$ by $\alpha(x,y) = f(xy)-f(x)-f(y)$.  It is easy to check that $d^3 \alpha =0$ in $C^3(\wl,\rho)$.  
 Since $\alpha(x,y) \in K$, it is bounded, i.e.   $\alpha \in Z^2_{bdd}(\wl,\rho)$.

  Define $\phi:  \wl \ta V \to \lc$ by $\phi(g,a) =  f(g)a$.  It is easy to check that $\phi$ is a homomorphism, and that 
  $\chi(f) \sim \{ g \in \sl:  \phi(g,0) \in U \}$, for $U$ a compact open neighborhood of $1$ in $L$.  

\eprf

   Let $\beta$ be the bounded cohomology class represented
  by $\alpha$.  If $\beta=0$, then $\om$ is 
laminar: we have, above, $\alpha=  d^2 (b)$ for some bounded $b: \wl \to A$.  Define $f'(x) = f(x) b(x) \inv$.
Again it is easy to check that $f': \wl \to \sl$ is a homomorphism.  As $b(\wl)$ is a bounded subset of $A$, and hence
a precompact subset of $\sl$, we have $\chi(f)=\chi(f')$.  This shows that $\om$ is laminar.  

Similarly, if $\beta$ is represented by $\alpha$ and also by $\alpha'$, then the corresponding approximate subgroups are
commensurable.
%

 
%

 
An alternative description can be given in terms of cohomological data on $G$ and on the   quotient group $\bar{L} = \sl/A$. 
  In this case, the
data consists of a homomorphism $\bar{\phi} : \wl \to \bar{L}$, a finite-dimensional Euclidean space $V$,  a homomorphism
$\rho: \bar{L}  \to O_n = Aut(V)$,  a class $\beta \in H^2_b(G,V)$, agreeing in $H^2(G,V)$ with the $\bar{\phi}$-pullback
of the class describing the extension $L$ of $\bar{L}$.   More precisely:

\begin{cor}  \label{cohcor} Let $G$ be a group, $\omega$ a commensurability class of approximate subgroups of $G$.
Then $\omega=\chi(\bar{f},\bar{L},L, \beta)$ is determined by:
\begin{itemize} 
\item A connected Lie group $L$ with no compact normal subgroup; 
\item  A rigid abelian normal subgroup $A$ of $L$; let $\bar{L}  = L/A$, and let $M$ be the $\bar{L}$-module $A$.  
Let $\alpha \in H^2(\bar{L}, A)$ be the cohomology class corresponding to the extension $1 \to A \to L \to \bar{L} \to 1$,
 {\em viewed as   discrete groups}.
\item  A subgroup $H$ of $G$, to which $\om$ belongs. 
\item A homomorphism  $\bar{f}: H \to \bar{L}$.
\item  A class $\beta \in H^2_b(H,M)$, such that the image  of $\alpha$  in $H^2(H,M)$ coincides with the image
$(\bar{f}|H)^* \beta$ of  $\beta$ in $H^2(H,M)$.
\end{itemize}
 $\omega$ is laminar   if  $\beta=0$.  
 
In particular if $H^2_b(H,M) = 0$ for all $H$ such that $\om$ belongs to $H$, then $\om$ is laminar. 
\end{cor}
Except perhaps for the choice of $H$,  that is determined only up to finite index, 
 the data is canonically determined by $\omega$:  we begin with the locally compact pattern automorphism group;  $L$ is the unique connected Lie group with no compact normal subgroups; $A$ is the closed subgroup arising from the Lascar neighbour 
 relation $\Sh$; $\beta$ is defined as above.


  Conversely, given data $L,A, H, \beta,\bar{f}$, we can construct an approximate subgroup in $H$.   Let   $L^\delta$ be the discrete
  group underlying $L$, and represent it as $\bar{L} \times_{\alpha} A$.  Then $\alpha,\beta$ determine classes
  in $H^2(H,A)$, and by assumption we have  $\alpha - \beta = 0 \in H^2(H,A)$.   Let $\gamma: H \to A$ be a map
  with $d \gamma = \alpha - \beta$. Define $\phi(x) = (\bar{g}(x),\gamma(x))$.   The error set is still given by that of $\beta$, and 
  thus is contained in $A$ and bounded.
    
%


 Burger and Monod have proved vanishing theorems for  lattices in semisimple  groups $L$.  This will not immediately help us in
  understanding approximate lattices $\Lam$ of $L$, since  the criterion of  \corref{cohcor}   requires vanishing   of the bounded cohomology not of $L$ itself, but of the subgroups $G$ generated by 
   approximate subgroups.
  Prima facie, these groups $G$ seem perfectly likely to be free and hence have very large bounded cohomology.    We will   analyze
  $\Lam$ by other methods.  But after the fact, we will know that $G$ is in fact isomorphic to a lattice in a semisimple group
  (not in $L$ itself but in a product of $L$ with a `complementary' Lie group);  making the results of  \cite{burger-monod} and \cite{monod} relevant.  
   It will be too late for $G$ itself, but along with Gromov's mapping theorem, we will be able 
   to apply the criterion to soluble extensions of $G$.  
See \propref{othergroups}.

 With \propref{notsame}, we will see  a partial converse to    the statement that the vanishing of $\beta$ implies laminarity.  Let $\Lam$ be an approximate subgroup of $G$,
  let $L,A,\bar{L}$ be as above with $A \cong \Rr$ central in $L$, and   let $\beta  \in \ker( H^2_b(G,\Rr) \to H^2(G,\Rr))$
  be the associated class.  If $\Lam$ is laminar, then $\Lam$ belongs to a subgroup $H$ of $G$ such that
  the image of $\beta$ in $H^2_b(H,\Rr)$ vanishes.  Indeed we have $\beta = d \alpha$ for some quasimorphism 
  $\alpha: G \to \Rr$.  Let $\alpha'$ be a homogeneous
 quasimorphism with $\alpha-\alpha'$ bounded.   By \propref{notsame}, $d \alpha'$ is a homomorphism on 
 some subgroup $H$ representing the commensurability class of $\Lam$.  Passing to this subgroup $H$, we have 
 $d \alpha' =0$ so $d \alpha = d (\alpha - \alpha')$ and thus $\beta$ vanishes as a class in $H^2_b(H)$.

 \begin{rem} 
 It follows that any nonzero class in $H^2_b(G,\Rr)$ indicates a non-laminar approximate subgroup, if not in $G$ itself then 
 in a central extension:   the image of $\beta$ in $H^2(G,\Rr)$ determines a central extension $\widetilde{G}$ of $G$ by $\Rr$; 
 the image $\widetilde{\beta}$ of $\beta$ in $H^2_b(\widetilde{G}, \Rr)$ remains nontrivial by the  Gromov mapping theorem alluded to earlier;  but now the image of  $\widetilde{\beta}$ is trivial.  
 \end{rem}

\ssec{Approximate kernels of quasimorphisms are {non-laminar}.} \label{5.18}

The proposition below implies that an approximate kernel $E$ of a nontrivial quasimorphism is associated to no $\bigwedge$-definable group,
in any expansion of the language.      In other words, there is no   chain of approximate subgroups
$D_n$ commensurable to $E$,  with $D_{n+1}^2 \subset D_n$.      
In particular, the two classes of   approximate subgroups, arising from  homomorphisms into Lie group or  from (nontrivial) quasimorphisms into $\Rr$, have trivial intersection.

 \begin{prop}   \label{notsame}   
 Suppose $\om$ is a commensurability class of approximate subgroups, equal to $\chi(g)$ for some
  homogeneous quasimorphism $g: G \to \Rr^n$, and  also to $\chi(f)$ for some  homomorphism $f: G \to L$  into a connected
  locally compact group, with dense image.   Then 
   $g$ is a homomorphism. 
  \end{prop}
  
  \prf    
  Since $G = \union_d g \inv [-d,d]^n$, it suffices to prove for each $d$ that $g$ is a homomorphism
on the subgroup $<E>$ generated by $E:=g \inv [-d,d]^n$.   So we may  assume  $G=\la E \ra$.   
Factoring out a compact normal subgroup of $G$ does not change $\chi(f)$; so we may assume $L$ has 
h no compact normal subgroups.  
  Being homogeneous, $g$ is conjugation-invariant, so $E$ is normal, and thus $f(E)$ is normalized by $f(G)$ and so $cl(f(E))$ is normal in $L$.  But
 $f(E)$ is precompact, so $cl(f(E))$ is compact.    We can now apply  \propref{lc1cc} with $H=L$, $K=cl(f(E))$.  We must have
 $H_1=H=L$ and $N=1$.  Thus there exists a rigid $A \leq L$, with $\la cl(f(E)) \ra \subset A$.    But $\la  f(E)  \ra = f(G)$
 is dense in $L$; so $L=A$ is abelian, and   $f (\la [G,G] \ra) = 1$.  Since $\chi(f)$ and $E$ are commensurable,
 $\la [G,G] \ra$ is contained in finitely many translates of $E$, so $g( \la [G,G] \ra) $ is bounded. 
   By homogeneity, $g(x^n)=ng(x)$ for $x \in [G,G]$, so $\Zz g(x)$ is bounded and hence $g(x)=0$, i.e. $g$ vanishes on $  \la [G,G] \ra $.   
 Now let $x,y \in G$; then $x^ny^n (xy)^{-n} \in  \la [G,G] \ra $, so $g(x^ny^n (xy)^{-n}) =0$ and hence $g(x^n)+g(y^n)+g((xy)^{-n})$
 is bounded independently of $n$.  But $g(x^n)+g(y^n)+g((xy)^{-n}) = n(g(x)+g(y)-g(xy))$ so letting $n \to \infty$
 we have $g(x)+g(y)=g(xy)$, i.e. $g$ is a homomorphism.   
 \eprf
  
 If in \propref{notsame} we remove the assumption that $cl(f(G))$ is connected, we still obtain that $g$ is a homomorphism
 on some subgroup to which $\om$ belongs.  
 
%
%
%

\ssec{}  By way of contrast, we note a family of cases where an  approximate homomorphism $f$ with central error - hence in particular a quasimorphism - can be replaced with an actual homomorphism.  For instance this applies to a centerless finitely generated or linear   group, or a subgroup of a stable group.  Factoring out the center yields an actual homomorphism $\bar{f}$,  whose pullbacks of compact sets are not much larger in some sense: any 
any group topology on $G$ making the elements of $\chi(f)$  discrete, also must make the elements of $\chi{\bar{f}}$ discrete.  

See  \secref{complementarity} for the definition of discrete-on-compact.

\begin{prop}\label{centralcase}  Let $G,\lc$ be topological groups, $\tl$ a   subgroup of $G$.  
Let $K$ be a compact subset of the center $Z$ of $\lc$.    Assume:  

\begin{itemize}
\item   For  some finite set $F \subset \tl$, the centralizer $C_G(F)$   is finite.    
  \end{itemize}

Let  $f: \tl \to \lc:K$   a discrete-on-compact quasi-homomorphism.   Then the composition  $ \tl \to \lc/Z$ is a  
discrete-on-compact homomorphism.   \end{prop}

\prf  
 Let  $\pi: \lc \to \lcb:=\lc /Z$ be the quotient  map, and $\ff = \pi \circ f: \tl \to \lcb$.
    Let $T_h: \lc \to \lc$ be the map: $x \mapsto x \inv h x$.  Then $T_h$ factors through a map $t_h:\lcb \to \lc$.  By definition of the quotient topology, 
$t_h: \lcb \to \lc$ is continuous.  Hence it maps (pre)compact sets to (pre)compact sets.

Let $W$ be a subset of $\tl$ with $\ff(W)$ precompact.  It follows for any $g \in \tl$ that
 $t_{f(g)} (\ff (W)) $
is precompact.  Now $t_{f(g)}(\ff W) = T_{f(g)} ( f(W)) $ and
$ f(g^W)  \subset     t_{f(g)} ( f(W)) K^4 $ so  $f(g^W)$ is also precompact,
and hence $g^W$ is discrete in $G$.  

  Suppose $w_i \in W$ 
 and $w_i \to a \in G$.  Then for each $g \in F$, $w_i \inv g w_i \to  a \inv g a$.    Since $g^W $ is discrete in $G$, 
 it follows that $w_i \inv g w_i = a \inv g a$ for large $i$.  As this holds for all $g \in F$, it follows that all $w_i$ (for large enough $i$) lie in the same coset of $ C_G(F)$.  But $C_G(F)$ is finite.  Thus $W$ has no accumulation points in $G$, so it is discrete in $G$.
\eprf

\begin{question}  Can the `no chain' statement alluded to in the first paragraph of \secref{5.18} be made more effective?     For example, it follows from this statement
that there exists an approximate subgroup  $D$ commensurable to $E$, such that for no approximate subgroup $D'$
do we have $(D')^8 \subset D^4$.    How can such a $D$ be found or recognised?
\end{question}

\begin{example}   All known nontrivial definable approximate subgroups appear to live in theories with the strict order property.
 This is related to longstanding questions regarding simple theories, and to new ones concerning the existence property in theories of Shelah's class NSOP1.  It is interesting in this connection to examine a first example of an approximate subgroup arising from a quasimorphism.  We take Brooks' quasimorphism on the free group, see \cite{kotschick}, Example 1.   Let $F$ be the free group on two letters $x,y$, and $w$ a nontrivial cyclically reduced word; for definiteness take $w=xy$.    Let $A \subset F$ be the set of words
containing as many  nonoverlapping copies of $w$ as of $w \inv$;    $A$ is an approximate subgroup.   The group $x^{\Zz}=\la x \ra$ is contained in $A$.
Moreover $x^n y \in A$ for $n<0$ but not for $n>0$.  Hence the semigroup $x^{\Nn} = A \meet C_F(x) \meet A y \inv$ is definable
with parameters in $F$, and even definable in $A$ given with partial multiplication and conjugation by $x,y$.  Hence $(F,A,\cdot)$
has the order property.
\end{example}

 The following example is based on  the theory of generic types of \cite{pillay98},    extended to groups $G$  given as the limit of a directed family of commensurable subsets.   
This extension is explained below with respect to the existence of generics, and I expect that modifying the other points is routine.  

We need only the easier theory over a model, and the proof may apply to other theories with
a similar behaviour of generics; the chain condition on forking is not required, only the consequences for amalgamation. 
 The natural candidates would be theories with   Shelah's property NSOP1 (\cite{shelah-towards}, \cite{kr}).
 
\begin{example}  Every approximate subgroup defined in a simple theory is laminar.  \end{example}
\prf   Let  $G=\union_i G_i$ with $G_i$ symmetric and commensurable to each other, and
such that for any $i,j$, for some $k \geq i,j$ we hqve $G_iG_j \subset G_k$.     The commensurator of an approximate subgroup
admits such a presentation (\lemref{commlem}).    

We follow the definitions and proof of  \cite{pillay98}, over some base $M_0$; types and independence are over $M_0$ if not otherwise stated.  $p$ is {\em generic} if for any $a$, there exists $b \models p$ with $ab,a$ independent.   
    The immediate stabilizer $St_0(p,q)$ of one type $p$ into another $q$ 
is defined as the set of $a$ such that for some $b \models p$ and $c \models q$, $ab=c$ and $a,b$ and $b,c$ are independent.
We let $St(p)$ be the group generated by $St_0(p,p)$; it is $\bigwedge$-definable, and  generated in two steps,
indeed $St(p) \setminus St_0(p)$ is non-generic.   It follows that $St_0(p,q)$ agrees, up to a non-generic set, with a coset of $St(p)$; we denote this coset by $St(p,q)$. 
 We now follow \cite{pillay98}   to prove 
 existence of a generic type $p$ of $G$, lying
in any given definble set $G_i$.    Define $\Psi^i_\alpha(x)$
 as in the proof of  \cite{pillay98} Lemma 3.8, only starting with $G_i(x)$ at stage $\alpha=0$.  
    Since for any $\phi$ and $k$ the rank $D(\Phi,\phi*,k)$ is translation-invariant,
   and the rank of a finite union is the maximum of the ranks, for fixed $k_0$, the number  $n_0:=D(G_i,\phi*,k_0)$ independent of $i$.   By the same argument, inductively, $D(\Psi^i_\alpha(x) ) = D(\Psi^i_\alpha(x))$.  This allows  the proof of genericity
   of the type $p$ constructed for $G_i$ to go through. 
    
    Note that the immediate stabilizer of $p$ will lie in $G_i^2$,   and the stabilizer subgroup of $p$ will be contained in $G_i^4$.  
This is always the same group $C=G^{00}_0$, the smallest $M_0-\bigwedge$-definable  subgroup of $G$ of bounded index.
For two generic types $p,q$ within $C$, we have $St_0(p,q) \neq \emptyset$, $St(p,q) \subset C$, and $St(p,q)$ is a coset of $C$,
so $St(p,q)=C$.  Hence for any generic types $p,q,r$ in $C$ we can find $a,b,c \in C$ with
 $a  \models p$, $b \models q$, $c \models r$ and  and $ab=c$.   
 
 At this point an elementary but slighty less perspicacious argument shows that $pp \inv p p \inv = C$ and then using compactness, that any the image of any $G_i$ modulo $G^{00}$ has open interior, giving laminarity.   However we invoke instead the picture given by \thmref{ag2mt}.  If laminarity fails, then according 
 to this theorem, after perhaps replacing $G$ by the group $\widetilde{\Lam}$ generated by some definable set commensurable to $G_i$,  there exists a continuous quasimorphism $\alpha: G \to \Rr: [-1,1]$ with unbounded image.  
 Continuity implies that $\alpha(x)$ depends only on $tp(x)$.  
 Choosing 
 $p,q,r$ with $\alpha(r) > \alpha(p)+\alpha(q) +1$,   and elements $a,b,c$ as above, we see that the inequality contradicts
   the quasimorphism property.

 \eprf


     \pb

  \pb
\section{Approximate lattices in semisimple groups}

\label{discreteapprox}

 Let $G$ be a locally compact group.  
Recall the definition of an approximate lattice from the introduction, and  Appendix \ref{approxlattices}.

 
  
The definitions below parallel those in \cite{BjH}; we allow ourselves to define {\em Meyer sets} for all lattices (in our sense),
so that a  Meyer set in the sense of \cite{BjH} is what we call a {\em uniform Meyer set}.   (The requirement that $M,S$
be   commensurable is prima facie stronger than that in \cite{BjH}, but they are equivalent, see \lemref{commens}.)
 
 When we have a product $G \times H$ of two groups, we will usually denote the projections by $\pi_1,\pi_2$.
 
\begin{defn}[\cite{BjH}] \label{modelset} Let $G$ be a locally compact group.   \begin{itemize}
\item  A  subset $S \subset G$ is a {\em (uniform) model set}  if 
 there exists a locally compact group $H$, a (co-compact) lattice $\Gamma \subset G \times H$,
and a compact    $C \subset H$ with nonempty interior, such that  $\Gamma$ has dense image under  the projection $\pi_2$,  $\pi_1$ is injective on $\Gamma$, and  
\[S= \pi_1 ( \Gamma \meet (G \times C)) \]
$S$ is also called a {\em cut-and-project set} for $\G$.    
 \item A {\em (uniform) Meyer set} is  a  subset $M$ of a (uniform) model set $S$, commensurable to $S$.  
\end{itemize}
\end{defn}

 This generalizes a construction of lattices, fundamental to the definition of an arithmetic lattice 
 where $C$ is taken to be a compact
  {\em subgroup}; see  \cite{margulis} p.1.   Here it is replaced, in effect, by an approximate subgroup.   
    
 The injectivity and density assumptions in \defref{modelset} are inessential:   we can replace $H$ by the closure $H'$ of the image of $\pi_2(\Lam)$,   then factor out the closed normal subgroup
$\pi_2( \ker(\pi_1))$ from $H'$, without changing the cut-and-project sets.

Finding a 
   cut-and-project model $\Gamma \leq G \times H$ for an approximate sublattice $\Lambda$   is   equivalent to the existence of 
homomorphism $f:  \wl \to H$,  viewed as a partial map $G \to H$, with the  interesting property  of interchanging compact  sets in the image  with discrete   sets in the domain; the relation between them is that $\G$ is the graph of $f$.       

    The identity map on the ring $\Zz[1/2]$, as a partial map from the $2$-adic completion to the real one, or vice versa, is a suggestive example of this phenomenon.  

   \begin{problem}[Bj\"orklund and Hartnick, \cite{BjH} p. 2918] \label{bjhp1} In what locally compact  second countable groups $G$ is every cocompact approximate lattice a Meyer set? 
 \end{problem}
 
This  is equivalent to the apparently weaker statement that every approximate lattice
$\Lam$ be commensurable to a cut-and-project set:

\begin{lem} Let $G,H$ be locally compact groups, $\G$ a lattice in $G\times H$ projecting injectively to $G$,  $\Lam$ an approximate subgroup of $G$
commensurable to a cut-and-project set of $\G$.  Then there exists a locally compact $H'$ and lattice $\G' \leq G \times H'$,
such that $\Lam$ is a Meyer set of $\G'$, i.e. commensurable to and contained in a cut-and-project set of $\G'$.

   \end{lem}
\prf  Let $\wl$ be the projection of $\G$ to $G$, and $f: \wl \to H$ the homomorphism whose graph is $\G$.  
Let $\om$ be the commensurator of $\Lam$.  By \lemref{commsmooth}, there exists $f' : \om \to H'$ such that
$\chi(f')=\chi(f)= \om$.  Let $\G' \leq G \times H'$ be the graph of $f'$.  By \lemref{commlem}, $\Lam  \subset \comm(\om)$.  Since $\Lam \sim \chi(f')$, there exists
a compact $C \subset H'$ with $\Lam \subset (f') \inv(C)$.  Thus $\Lam$ is a Meyer set of $\G'$
\eprf
   
   Note that this may require a change of target $H$; even if $G,H$ are connected Lie groups, it may be necessary to augment $H$ with totally disconnected factors.  This will be clear in the arithmetic setting below, where $p$-adic factors appear.

      Existing results fall into the amenable class; see \cite{machado3}. 
We will complement this with a positive answer, without the restriction to cocompact approximate lattices, 
for  almost simple groups over local fields, and their finite products.

 

  \ssec{Abstractly semisimple groups}    
 
 \begin{defn} \label{ass} 
   A  locally compact group $G$ is {\em abstractly semisimple} if    \begin{enumerate}
 \item $G$ has no nontrivial normal abelian subroups.
 \item  Outside the center, $G$ has no discrete conjugacy classes.
 \end{enumerate}
 \end{defn} 

  \begin{rems} \label{assrems}  \begin{enumerate}
  
 \item    \label{assd} Abstract semisimplicity remains true upon passing to a finite index normal subgroup $H$ of $G$.
 \item   \label{products}  Abstract semisimplicity is preserved by finite products.
  \item  \label{adeless}  
The class is  also closed under restricted infinite products $\Pi'_i G_i$ with respect to compact subgroups $C_i$ (see \defref{groupdefs} (\ref{restricted}.)), provided each $C_i$
has trivial centralizer in $G_i$.

      \item  \label{agss}  A semisimple algebraic group over a local field is abstractly semisimple, at least upon moving to a finite index subgroup, and factoring out a finite center.    These operations will not alter the validity of the results we will prove.

  \item  \label{bt1} 
 Let $G$ be an abstractly semisimple locally compact group.   Then  for some compact normal subgroup $N$,  $G/N$ has a finite index subgroup of the form $G_1 \times G_2$ with $G_1$
  a totally disconnected abstractly semisimple locally compact group, and $G_2$ a connected semisimple Lie group. 
     \end{enumerate}
   \end{rems}
   
    \prf   
    (\ref{assd}) If $X$ is a discrete conjugacy class of $H$, then the union of $G$-conjugates of $X$ is a conjugacy class of $G$ and is still discrete (the union being finite), so trivial .  If  $A$ is a normal abelian subgroup of $H$, let $A_1,\ldots,A_n$ be the conjugate subgroups.  Then $\meet_{i=1}^n A_i$ is $G$-normal and abelian, so trivial.  We show by reverse induction on $k=n,\cdots,1$ that the
intersections of any $k$ of the groups $A_i$ is trivial.  Assume this true for $k+1$; we prove it for $k$.    If $A,B$
are each the intersection of $k$ of the $A_i$ and are distinct, then  $A \meet B = (1)$; so 
 the various $k$-intersections commute and so generate a commutative $G$-normal subgroup, hence  trivial.   For $k=1$
 this gives $A=(1)$.   
 
  (\ref{products})  The statement on discrete conjugacy classes is clear, since
a conjugacy class in the product is a product of conjugacy classes.   If $A$ is a normal abelian  subgroup of $G$,
the projection of $A$ to each factor is normal and abelian so trivial, hence $A$ is trivial.

(\ref{adeless})  The proof that there are no normal abelian subgroups is the same as in (\ref{products}).   
If $X$ is a discrete conjugacy class of $G$,
pick $x \in X$; then the centralizer $C_G(x)$ is an open subgroup of $G$;  thus for some finite $I_0 \subset I$,
 $1_{I_0} \times \Pi_{i \in I \m I_0} C_i$ is a subset of $C_G(x)$.    From this it follows that $x \in G_0:= \Pi_{i \in I_0} G_i \times 1_{I \m I_0}$.
  So the projection 
$\pi_0 X$ 
of $X$ to $G_0$ must already be discrete.  By (\ref{products}), $\pi_0 X = (1)$.  Since $x \in G_0$ we have $x=1$.  

(\ref{agss})  By factoring out the center and moving
 to a finite index subgroup,  we may take $G$ to be a finite product of simple 
 groups.  Then it is clear that $G$ has no nontrivial normal abelian subgroups, nor open ones.  If $a^G$ is a discrete conjugacy class,
 then $C_G(a)$ is open, and using the chain condition on centralizers, $C_G(a^G)=C_G(F)$
 for some finite set $F$ of conjugates of $a$; so $C_G(a^G)$ is an open normal subgroup; hence $C_G(a^G)= G$ so $a$ is 
central.

  (\ref{bt1}) If $G$ is connected, this  is Yamabe's theorem (with $G_1=1$).  Moreover in this case, $G$ has a unique maximal
  compact normal subgroup.   In general, factoring out the maximal compact $N$ of the connected component, $G^0$  becomes
 a semisimple Lie group.   It has have only finitely many   outer automorphisms, using the structure theory for semimple Lie groups, and Borel-Tits; thus up to finite index, $G/N$ decomposes as a product of $G^0$ with the totally disconnected group $C_G(G^0)$.
 
\eprf
  
As a consequence of (3), if  $G$ is a semisimple algebraic group over $\Qq$, then $G(\Aa_{\Qq})$  has a co-compact subgroup that is abstractly semisimple modulo a compact center.   For instance
   $SL_2(\Aa_{\Qq}) / Z$ is abstractly semisimple, where $Z$ is the center $\{\pm 1\}^{\{\infty,3,5,7,11\cdots\}}$.

Let $G$ be an abstractly semisimple group, and let 
$\Lam$ be a  discrete approximate subgroup.  
\thmref{ag2} explains    $\Lam$ in terms of a quasi-homomorphism into a Lie group; we would like to show that an ordinary homomorphism will already work.  In case the error set is central,
this is already carried out in  \propref{centralcase}, assuming in a strong sense that $\G$ is centerless: 
  the centrality of the error set is played against the semisimplicity of $G$, showing we can  factor out the centre of $L$ without harm.
 The general case, where the error set is contained in a rigid normal subgroup $A$, is more complicated; 
$A$, when factored out, may require taking with it a bigger   normal subgroup, that needs to be shown to be solvable.
 We postpone the technical details to \propref{ctp} at the end of the section, and prove the theorem assuming it.

\begin{defn}  Let $G$ be a topological group, and $\om$ a commensurability class of approximate subgroups.  $\om$ is
{\em strictly dense} if it does not belong\footnote{as defined in the beginning of \secref{approxsg}} to a closed subgroup of infinite index in $G$.    
\end{defn}
Note that lattices are never strictly dense; see also the Example below.  In view of this we  will also say that $\Lam$ is
{a {\em strictly approximate lattice}} if the commensurability class of $\Lam$ is strictly dense.  

In case   $\Lam \in \om$ and  $\la \Lam \ra \in {\la \om \ra}_{min}$ in the sense of \propref{minimalclass},  $\Lam$ is  a fortiori
strictly dense in $G$, of course.
 
\begin{example} \label{nondense}
   If $G=G_1 \times G_2$, $\Lam_1$ an approximate lattice in $G_1$,
$\Lam_2$ a lattice in $G_2$, then $\Lam_1 \times \Lam_2$ is an approximate lattice in $G$, that is not dense;
it   belongs to the closed subgroup $G_1 \times \Lam_2$.  
\end{example}
  We prove the first
{laminar}ity result assuming strict density.  Later, in  the classical setting at least,  we will recognize \exref{nondense}  as
 the only obstruction to strict density, so that the theorem applies to all approximate lattices.

\begin{thm} \label{discrete1} Let $G$ be an abstractly semisimple locally compact group.   Let  $\Lam$ a  discrete approximate
subgroup of $G$, with strictly dense  commensurability class.    Then there exist  a  connected semisimple real Lie group $\sl$ and a  discrete subgroup
$\Gamma \leq G \times \sl$ with  cut-and-project sets    commensurably covering   $\Lam$.  

\noindent{{\bf Supplements:}} 
\begin{itemize}
\item $\sl$ can be taken to be centerless and to have no nontrivial  compact normal subgroups.   
\item   We can take the projection  $\Gamma \to \sl$ to be injective with dense image, and $\Gamma \to G$ to
 be injective.
 \item 
 If $\Lam$ is an approximate lattice in $G$, then $\G$ is a lattice in $G \times \sl$, and the cut-and-project sets of $\G$ are commensurable to $\Lam$.
\end{itemize}
 \end{thm}

 \prf 
  By \thmref{ag2}, there  exists  a subgroup $\wl \leq \la \Lam \ra$, a Lie group ${L}$ with 
 $L/Z(L)$ connected, 
 a normal compact subset $K$ of a rigid abelian subgroup $A \leq {L}$, 
 and a   quasi-homomorphism 
 \[f:  \wl \to {L}:K\]
 with $f(\wl)A$   dense in $L$,
 such that $f\inv(C)$  commensurably covers $\Lam$
 for some compact $C \subset {L}$ with nonempty interior, and in fact is commensurable with $\Lam$.   
 We use the commensurability here to conclude that $f \inv(C)$ is discrete in $G$; so $f$ is discrete-on-compact.

By assumption, $cl(\wl)$ has finite index in $G$; so it is still abstractly semisimple (\remref{assrems} (\ref{assd})).  Thus we may assume 
$cl(\wl)=G$; and (i),(ii) of   \propref{ctp} hold (with $H=\wl$.)
%
  By \propref{ctp},   $L$ has a normal subgroup $R$ such that $L/R$ is semisimple, centerless, without nontrivial normal compact subgroups; and  the induced map $\bar{f}: \wl \to L/R$ is injective and discrete-on-compact.    
 Since $A \subset R$, $\bar{f}$
 is a homomorphism,
  The cut-and-project sets of $L/R$ can only be bigger than those of $L$.
 Let $\sl = L/R$, and let $\G$ be the graph of the composed map $\bar{f}: \wl \to L \to L/R=\sl$.  Since $\bar{f}$ is 
 discrete-on-compact,   $\G$  is discrete.   The projection $\Gamma \to G$ is injective by definition, and we saw that $\G \to L$ and hence $\G \to \sl$ have dense image.
 That the kernel of $\G  \to \sl$ is trivial follows from injectivity of $\bar{f}$; the kernel on the left is trivial since $\bar{f}$ is a function.
 
   Assume now that $\Lam$ is an approximate lattice.   Let $X$ be a cut-and-project set of $\Gamma$, commensurably covering $\Lam$.  
   As $\Gamma$ is discrete, it follows that $X$ intersects each compact of $G$ in a finite set, so that $X$ is discrete.   Now $X^4$ is also contained in a  cut-and-project set of $\Gamma$, so by the same argument $X^4$ is discrete.
  The commensurability of $\Lam$ and $X$  now follows from \lemref{commens}.   The fact that $\Gamma$ is a  
  lattice,   cocompact if $\Lam$ is, comes from \lemref{latticecoherence}.

  \eprf

\begin{rem} \label{discrete1r}  
  Assume $G$ has the centralizer chain condition (3) (as is  automatic in a linear group).   
Then in  \thmref{discrete1},
the strict density requirement can be replaced by the weaker conditions (1),(2).    In particular, it suffices to assume
that $<\Lam>$ is dense, and that $G$ has no closed normal subgroup $N$ with $G/N$ discrete.
\begin{enumerate}
\item    $\Lam$ is weakly Zariski dense.  This refers to the topology on $G$ whose basic closed sets are defined by equations $ax=xb$  or $b x \inv ax = x \inv a x b$, with $a,b \in G$.   Notably normalizers $N_G(C_G(X))$ of centralizer groups $C_G(X)$ are $\wz$-closed.
\item  (Weak irreducibility)  There is no closed subgroup $H$ of $G$ 
such that $\Lam \subset N_G(H)$ and $\Lam / H$ is  commensurably covered by an infinite  discrete subgroup of $N_G(H)/H$.  
\item \label{weaklin}    Any centralizer group $C_G(X)$ is the centralizer group $C_G(X_0)$ of a finite subset $X_0$ of $X$.
(Or at least, if $H$ is a closed subgroup and $C_G(x) \meet H$ is open in $H$, then the interesection of conjugates of 
$C_G(x) \meet H$ in $H$ is still open in $H$.)
\end{enumerate}\end{rem}
\prf  
Let $\wl$ be a subgroup as in the beginning of the proof of \thmref{discrete1}.  It suffices to check that $cl(\wl)$ is 
abstractly semisimple (and then proceed with the proof of \thmref{discrete1}.)  If $A$ is an abelian subgroup of $G$ normalized by $cl(\wl)$, then the centralizer group $C_G(C_G(A))$
is also normalized, hence by weak Zariski density it is normal in $G$, and thus by abstract semi-simplicity of $G$ it is trivial;
so $A=1$ also.   Now let $a \in cl(\wl)$ and suppose the conjugacy class of $a$ in $cl(\wl)$ is discrete.  Then 
$cl(\wl) \meet C_G(a)$ is an open subgroup of $cl(\wl)$.  By (3), $cl(\wl) \meet C_G(a^{cl(\wl)})$ 
is still an open subgroup $H$ of $cl(\wl)$.  Thus $cl(\wl) / H$ is discrete.  As $\Lam$ is commensurably covered by
$cl(\wl)$, by (2), it follows that $cl(\wl) / H$ is finite, i.e. a finite index subgroup of
$cl(\wl)$ centralizes $a$.  By weak Zariski density, a finite index subgroup of $G$  centralizes $a$.  So $a^G$ is finite and hence
trivial, and $a=1$.  
  \eprf

 \begin{cor}  \label{semisimple1}
 Let $G$ be a finite product of non-compact, centerless, simple algebraic groups over local fields.
   Let   $\Lam$ be an  approximate lattice of $G$.   Then there exists a decomposition $G=G_1 \times G_2$, $\Lam$ is commensurable to $\Lam_1 \times {\Lam}_2'$
where $\Lam_1 \leq \la \Lam \ra$ is an approximate lattice of $G_1$ and $\Lam_2'$ is a lattice of $G_2$,   and
there exists a semisimple real Lie group $\sl$ and an  approximate lattice 
$\Gamma \leq G_1 \times \sl$ with  cut-and-project sets   commensurable to $\Lam_1$. \bigskip

\noindent{\bf Supplements:} \begin{enumerate}
\item       If $\Lam$ is a uniform approximate lattice, then $\Gamma$ is a uniform lattice, as is $\Lam_2'$.  
\item     The projections $\Gamma \to G_1,  \Gamma \to \sl$ can be taken to be injective and a dense embedding,   respectively.      
\item  
The projection of $\Gamma$ to $G_1$ is a subgroup $\wl_1$.  One can choose $\G$ so that $\wl_1$ is minimal, in the sense
that no subgroup of infinite index commensurably contains $\Lam_1$; and   $\sl$ connected and without normal compact subgroups.  In this case, 
$\sl$ is uniquely determined,  and $\G,\wl_1$ are determined up to commensurability. 
 
\item    There exists no  nontrivial discrete  subset  of $G_1$   normalized by $\Lam_1$;
  in particular, no nontrivial subset normalized by $\Lam_1$ and commensurably covered by $\Lam_1$.  
  \end{enumerate}
    \end{cor}

  \prf 
   Let   $\wl, \sl$ and $f:  \wl \to \sl:K$ be as given by \thmref{ag2}, applied to $\Lam \subset \la \Lam \ra$.  We may at this point replace $\Lam$ by the commensurable $\Lam^2 \meet \wl$, and so assume $\Lam \leq \wl$.  
    If $cl(\wl)$ is itself abstractly semisimple,
we proceed as in  \thmref{discrete1},
 obtaining  $\sl$ and a discrete $\Gamma \leq G \times \sl$ with a cut-and-project set   $X$   commensurably covering $\Lam$. 
 We repeat the argument that $\Gamma$ is a lattice (in $G \times \sl$ and not only in $cl(\wl) \times \sl$):
     since $\Gamma$ is discrete, it follows that $X$ intersects each compact of $G$ in a finite set, so that $X$ is discrete.   Now $X^4$ is contained in a  cut-and-project set of $\Gamma$, so by the same argument $X^4$ is discrete.
  The commensurability of $\Lam$ and $X$  now follows from \lemref{commens}.   The fact that $\Gamma$ is a  
  lattice,   cocompact if $\Lam$ is, comes from \lemref{latticecoherence}. 
   The statements on injectivity and density are as in 
  \thmref{discrete1}.      Supplement (3)   follows from \propref{sensible}.    
 
  We thus need to handle the case where $cl(\wl)$  is not abstractly semisimple.  In any case, $\Lam$ and thus $cl(\wl)$
  are Zariski dense by \propref{boreldensity}.   Let $A$ be an  abelian subgroup of $G$ normalized by $\Lam$.  To show $A$ is trivial it suffices to show the triviality of every projection  to one of the simple algebraic factors of $G$.  Let $\pi_i:G \to G_i$ be such a projection;   then $\pi(A)$ is abelian and normalized by $\pi(\Lam)$;  hence the double centralizer of $\pi(A)$, a Zariski closed group,  is   normalized by $\pi(\Lam)$; so by Zariski density of the latter, it is trivial.

In a similar way one can show that a centralizer subgroup of $G$ normalized by $\Lam$ is normal.
Indeed let $C=C_G(Y)$; both $C$ and $N_G(C)$ are defined by coordinatewise conditions, so that if $\pi_i(C)$ is normal in $G_i$ for each $i$, 
then $C$ is normal in $G$.  Now $\pi_i(C)$ is a centralizer subgroup of $G_i$, hence Zariski closed; and thus so is
$N_{G_i}(\pi_i(C))$; since the latter contains the Zariski dense set $\pi_i(\Lam)$, it must equal $G_i$.

  The remaining issue is  condition (ii)  of abstract semisimplicity, forbidding discrete conjugacy classes.   
  Let $D$ be the set of all elements $a \in cl(\wl)$ with $a^{cl(\wl)}$ discrete in $G$; and let $G_1 =  C_G(D)$
    be the centralizer.   
 $G_1$ is clearly normalized by $cl(\wl)$ and hence by $\Lam$, so by the above, it is a  closed normal subgroup of $G$.   

  Let $G_2=C_G(G_1)$;   as $G$ is a product of simple groups, it is clear that $G_1$ is a product of some of them, and
  $G_2$ of the others.  We have $G \cong G_1 \times G_2$ naturally; $G_2=C_G(G_1) $ contains $D$ by definition.

 By linearity of the factors of the product group $G$,    $G_1=   \meet_{a \in D} C_G(a)$ is   a finite intersection, i.e. equals $C_G(D_0)$ for some finite $D_0 \subset D$.   
  For $a \in D_0$,  let $U$ be a neighborhood of $a$ with $U \meet a^{cl(\wl)}  =\{a\}$; then 
   $C_G(a) \meet {cl(\wl)} = \{g \in cl(\wl):  a^g \in U \}$, so it forms an open 
  subgroup of $cl(\wl)$; and hence $G_1  \meet cl(\wl)$ is open in $cl(\wl)$.   
 It follows that $cl(\wl) / (G_1  \meet cl(\wl)) \cong cl(\wl) G_1 / G_1$ is discrete as a topological space, and so also as
 a subgroup of $G/G_1 \cong G_2$.
  The   image $\pi_2(cl(\wl)) = \pi_2(\wl)$ is thus a discrete, closed subgroup of $G_2$.  
   A fortiori, $\pi_2(\Lam)$ is discrete in $G_2$.  
   By   \propref{productdec}, 
 $\Lam_1 = \Lam^2 \meet G_1$  is an approximate lattice in $G_1$;      $\Lam_2=\Lam^2 \meet G_2 $
   is an approximate lattice in $G_2$; and $\Lam_1 \times \Lam_2$ is commensurable with $\Lam$.       Since
   $\pi_2(\wl)$ is a discrete subgroup containing $\Lam_2$, it is commensurable with $\Lam_2$ by \lemref{commens};
   and so $\Lam_1 \times \pi_2(\wl)$ is also commensurable with $\Lam$; the second factor $\pi_2(\wl)$ is a lattice in $G_2$.
   Let $\Lam_2' = \pi_2(\wl)$. 
    
     Now  $G_1$ has no discrete conjugacy classes and so is abstractly semisimple.
     
 $f(\Lam_2)A$ is a normal subgroup, and contained in a compact set  (since $f(\Lam)$ is).   
 Since $\Lam_1 \times \Lam_2$ is commensurable with $\Lam$, $cl(f(\Lam_1)A) = cl (f((\Lam_1 \times \Lam_2)A)$
 is co-compact in  $\sl$.   Under the conditions of \thmref{discrete1}, that we have for $\sl$, this implies $cl(f(\Lam_1)A)=\sl$.
  
  We are now in the situation of the first paragraph with respect to $G_1$.    In particular, as shown there, (3) holds for $G_1$.
  Using induction on the number of simple factors of $G$, we obtain the main theorem and supplements (1-3).  
  
%
%
 
It suffices now to prove supplement (4) under the additional assumption (*) that no subgroup of infinite index
  in $\wl$ commensurably contains $\Lam$.   
  Suppose $S$ is a discrete  subset  of $G_1$, normalized by $\Lam_1$, and 
  $1 \neq a \in S$.    Let $\Lam' = \Lam_1 \times \Lam_2$.   Since $\Lam_1$ normalizes $S$, as does $\Lam_2$,
  so does $\Lam'$, and hence   the group $\wl'$ generated by $\Lam'$ normalizes $S$.  
  Now $S$ is closed in $G$, so $cl(\wl ')$ normalizes $S$ too.   By (*), $\wl ' $ has finite index in $\wl$.  Hence
  $cl(\wl')$ has finite index in $cl(\wl)$.   Thus    $a^{cl(\wl)}$.   In the notation within the proof, we have $a \in D$
  and hence  
  $G_1$ commutes with $a$.  But $a \in G_1$ so $a$ also commutes with $G_2$, and hence is central in $G_1 \times G_2=G$,
   contradicting our assumptions.   This proves (4) under (*), and hence as noted above in general.

   \eprf

 This solves Problem 1 of \cite{BjH} for semisimple algebraic groups over local fields, and their finite products.  
    
   \begin{center}{*}    \end{center}
 
We are still in debt of a technical lemma.  When applied to lattices in (real) Lie groups, it is a classical statement 
 going back to Zassenhaus and Auslander on the interaction of a lattice with the solvable radical.  For us
 the main import will be an ability to  convert approximate homomorphisms  with error in a basic abelian group, to actual homomorphisms into a quotient Lie group.

   \begin{prop} \label{ctp} \label{7.15} Let $L$ be a Lie group, with $L/Z(L)$  connected, containing a normal closed 
  abelian subgroup     $A \normal L$, and a compact normal symmetric $K \subset A$.    
    Assume $A=VD$,  where $V \cong \Rr^n$ is closed and normal in $L$, and $D$ is closed and central in $L$.
  Further assume that either $V$ is rigid, or $K=(1)$.   

   Let  $G$ be a second countable locally compact group,  with a subgroup ${H}$ 
      satisfying:
 \begin{enumerate}[label=(\roman*)]
 \item $G$ has no nontrivial    abelian subgroups normalized by ${H}$. 
 \item $H$ has no conjugacy classes that are discrete in $G$, other than (1).
 
     \end{enumerate} 
     
  Let $f: {H} \to L:K$ be a  discrete-on-compact quasi-homomorphism, with  $f({H})A$   dense in $L$.   
    
  Then      $f$ induces an injective, discrete-on-compact homomorphism $\ff: {H} \to L/R$, where $R$ is a closed normal solvable subgroup of $L$,   and $L/R$ is semisimple and centerless, and without nontrivial normal compact subgroups.
   \end{prop}
     
  \prf

  We first show:
   
 \claim{1}   $f \inv(A) = 1$.  
 
 \prf 
  If $A$ is rigid,    let $g \in f \inv(A)$, and let $g^{{H}}$ be the set of ${H}$-conjugates of
$g$.  Then 
$f(g^{{H}}) \subset  f(g)^{L} K^3$; 
since $A$ is rigid,   $f(g)^{L}$ and  hence also $f(g)^{L} K^3$ are precompact.
So $g^{{H}} \subset f \inv f(g^{{H}}) $ is discrete (in particular closed) in $G$.  By (ii) we have $g=1$, showing  that  $f \inv (A) =1$.  

Now assume instead of rigidity that $K=1$.  Then $f$ is a homomorphism.   Hence $f \inv(1)$ is a normal, discrete subset
of $H$, and it follows directly from (ii) that $f \inv(1)=1$.    Thus $f$ is an injective homomorphism; so 
 $f \inv(A)$ is an abelian subgroup of $G$.  It is  normalized by $H$, hence trivial by (i).  
 
  \eprf

It is only in the above proof of Claim 1 that we use the rigidity of all of $A$, rather than just $K$.

\claim{2}  There exists a closed, normal, solvable subgroup $I$ of $L$ containing $A$  such that
the induced homomorphism 
  $f_1: {H} \to L/I$  is discrete-on-compact.   
  
\prf
  Any Lie group has a compact neighborhood $U$ of $1$, satisfying (*) $[U,U] \subset U$;  if $\mathfrak{u}$ is any compact neighborhood of $0$   in the Lie algebra, then for small enough $t>0$, $exp(t \mathfrak{u})$ will do.   Moreover,
  we can arrange  topological nilpotence (**):  letting $U_1=U, U_{n+1} =[U,U_n]$, we have $\meet_n U_n = (1)$.
  
 It follows that every   locally compact group $G$ has   a neighborhood satisfying (*):   by Yamabe's theorem, $G$ has an open subgroup $S$ containing a compact normal subgroup $N$, with  $S/N$ a Lie group; pull back from $S/N$ any neighborhood with this property. 
 
 Fix a  a compact neighborhood $W$ of $1$ in $G$ with  $[W,W] \subset W$, and a compact neighborhood $U$ of $1$ in $L$ with  (*) and (**).  Writing $\bar{U_n}$ for the image of $U_n$ modulo $A$, we can also ensure that $\meet _n \bar{U_n} = {1}$.

 By assumption, $A= V  D$ is a product of a vector group $V$, normal
  in ${L}$, with  a closed  group $D$ central in  ${L}$.    
   
Fix a Euclidean inner product on $V \cong \Rr^n$; let $B_r(V)$ be the ball around $0$ of radius $r$; normalize it so that 
$K^3 \subset B_1(V) \times D$.     
   Also a compact $B_1(D) \subset D$, 
such that $K^3 \subset B_1(V) \times B_1(D)$.   Let $B_r (D) = B_1(D)^{  \lfloor r \rfloor   }$, and let
$A_r = B_r(V) B_r(D)$. 

The inner product   induces an operator norm on $End(V)$.   For $g,x \in {L}$,   we will write   $x{^g} $  for $g x g \inv$.
Let $ad: L \to End(A)^{op}$ be the endomomorphism mapping $g$ to the endomorphism $x \mapsto   x{^g}$ of $A$. 
 Write $||x||$ for the operator norm of $ad(x)-1$ when restricted to $V$.    This is a continuous function on ${L}$.

 We can shrink $U$   so as to satisfy:  
\begin{enumerate}   \setcounter{enumi}{1}
\item \hspace{20pt}  $||g||  < 1/10$ for $g \in U^2$.   
\end{enumerate}
      
In fact, $ad([x,y]) = [ad(x),ad(y)]$ so it suffices to pull back a  neighborhood $\Omega$ of  $1$ in $GL_n(V)$ 
with the analogous property, and such that $[\Omega,\Omega] \subseteq \Omega$ (so as not to lose this property for $U$.)

  \claim{2.1}      Let $n \geq 1, r \geq 2$.  We have  \[ [U A_r, U_n A_r ] \subset U_{n+1}  A_{r/4} . \]
  
 \prf  
  Consider elements 
  $\g_i = u_i a_i$, $u_i \in {U}, a_i \in A$.   Assume   $u_1 \in U_1,u_2 \in {U}$ and $a_i \in A_r$  $r \geq 2$.
 
 We compute 
\[  [u_1a_1,u_2a_2]  =  [u_1,u_2]  ( [u_2 \inv,a_1] [a_1,a_2]  [a_2,u_1 \inv] )^{u_2u_1} \]
since $[a_1,a_2]=1$,
\[  \ \ \ \ \ \ \ \ \ \ \ \ \ \ \ \  [\g_1,\g_2] =  [u_1,u_2] ad_{u_2u_1} ( ^{u_2 \inv -1} a_1  + ^{1-u_1 \inv} a_2)  \]
where multiplication inside the abelian group $A$ is written additively.   
   If $u_1,u_2 \in {U}_n$ then $[u_1,u_2] \in {U}_{n+1}$.   
 And if $a_1,a_2 \in B_r(V)+D$, then $ ^{u_2 \inv -1} a_1$ and
 $ ^{1-u_1 \inv} a_2$  lie in $A_{r/10}$ by (2), so $( ^{u_2 \inv -1} a_1  + ^{1-u_1 \inv} a_2) \in A_{r/5}$,
 and again by (2),  $ad( u_2  u_1 )   ( ^{u_2 \inv -1} a_1  + ^{1-u_1 \inv} a_2) \in A_{r/4} $ 
 \eprf

 (Actually since $D$ is central we have the stronger statement  $[U B_rD, U_n B_r D] \subset U_{n+1}  B_{r/4}$.)
 
   \claim{2.2}   For $r \geq 2$, the set $\Delta_r = U A_r \meet f(W \meet {H})$   is finite.

  \prf
  Since $f$  is discrete-on-compact, and $UA_r $ is compact,    $f \inv (UA_r)$ is discrete, so $W \meet f \inv (UA_r)$
  is finite.   Applying $f$, we see that $f(W \meet {H}) \meet UA_r$ is finite.   
 \eprf
 
  \claim{2.3}     For any $r \geq 2$,     $\Delta_r$  generates a nilpotent group modulo $A$.
  
  \prf    Let $v_1,v_2 \in  \Delta_r$.  So for  some $u_1,u_2 \in U$, $a_1,a_2 \in A_r$, $w_i \in W \meet {H}$,
  we have $v_i = u_ia_i=   f(w_i)$.  
 By the proof of Claim 2.1,  $[v_1,v_2] =[u_1a_1,u_2a_2] =  [u_1,u_2]  a$ with $a \in  A_{r/4}$.  
 On the other hand,  
$ f([w_1,w_2]) = [f(w_1),f(w_2)] k $ for some $k \in K^3 \subset A_1$; so 
$[v_1,v_2] = [f(w_1),f(w_2)] =  f([w_1,w_2)]) k \inv$.  Hence $[v_1,v_2] k \in f(W \meet {H})  \meet UA_{r/4+1}$.
Since $r \geq 2$  this shows that  $[v_1,v_2] k \in \Delta_r$.   

Let $\bdr = \Delta_r/A$ be the image of the finite set $\Delta_r$ in $L/A$.  We have just shown that $\bdr$
is closed under the commutator.    Let $F_1 = \bdr$ and $F_{n+1} =[F_1,F_{n}]$.  Then $F_1 \supset F_2 \supset \cdots$ and the sequence must stabilize at some $F_k$ with $[F_1,F_k]  = F_k$.  Since $F_k \subset {\bar{U}}_1$, and $[{\bar{U}}_1,{\bar{U}}_n] \subset {\bar{U}}_{n+1}$,
it follows inductively  that $F_k \subset {\bar{U}}_n$ for each $n$.  Thus $F_k \subset \meet_n {\bar{U}}_n = (1)$.  This proves
that $\bdr $ generates a nilpotent group.
 \eprf

 Write ${\G}$ for the graph of $f$, and (with slight abuse of notation)
 \[{\G A}  := \G (1 \times A) = \{(w, f(w) a): \  w \in {H} , a \in A \} \]
 While ${\G}$ may not be a subgroup of $G \times L$, since $f$ is an $A$-quasi homomorphism, 
  ${\G A}$   is a subgroup of $G \times L$; hence so is  the closure $cl( {\G A})$
Let $C$ be  the identity component:
 \[ C :=    cl(  {\G A}  )^0 \]
  As $A$ is normal, $C$  is normalized by ${\G}$.   
  Let $\pi_2: G \times L \to L$ be the second projection.

      By a theorem of Zassenhaus  (see \cite{dixon}),  all soluble subgroups $S$ of the linear group ${L}/ Z({L})$
      are soluble of some fixed degree;  hence the same holds for $L$.  Call this bound ${\bb}={\bb}(L)$.  
       \footnote{The theorem is easy to prove for characteristic zero linear groups, such as $L/Z(L)$:  by passing to Zariski closure  it suffices to consider Zariski closed subgroups $H$; the solvable degree of the connected  part $H^0$ is clearly bounded by the dimension;  while Jordan's theorem   presents the quotient $H/H^0$  as a normal abelian subgroup by a bounded subgroup.}

\claim{2.4}    ${\G A} \meet ( W \times U) $  generates  a ${\bb}$-step soluble\footnote{In fact we can show here, using connectedness of $C^0$, that $\pi_2(C)/A$ is nilpotent.} 
 subgroup of $G \times L $.

\prf      An element of ${\G A}$ has the form $(w, f(w) a)$ with $w \in {H} $.   Thus an element of 
   ${\G A} \meet ( W \times U) $   has the form $(w, f(w) a)$ with $w \in {H} \meet W$ and $f(w) a \in U$.
  Letting $\pi_2$ denote the second projection $G \times L \to L$, we have
  $\pi_2({\G A} \meet ( W \times U) ) \subset (U \meet f(W  \meet {H}) A)$.     For a finite $F \subset {\G A} \meet ( W \times U)$, 
  $\pi_2(F) \subset U \meet f(W \meet {H}) A_r$ for some $r \geq 2$.     Now 
\[ \pi_2 (F) / A  \subset (U \meet f(W \meet {H}) A_r ) / A = (UA_r \meet  f(W \meet {H}) ) / A = \bdr.\]  
So $\pi_2 (F) / A$ generates a nilpotent group, by Claim 2.3.    Hence $\pi_2(F)$ generates a soluble group.  By definition of ${\bb}$,
it is ${\bb}$-step soluble.

Let $z$ be an element of the
 ${\bb}$-th derived group of $\la F \ra$.  Write $z = (x,f(x)a)$ with $x \in {H}$.  
 Then $\pi_2(z) = f(x)a \in A$ so $f(x) \in A$ and hence $x=1$  by Claim 1.   Thus   $\la F \ra$ is ${\bb}$-step soluble.  
 
Since this holds for all finite $F \subset  {\G A} \meet ( W \times U)$, the group
 $\la   {\G A} \meet ( W \times U) \ra$ is  ${\bb}$-solvable.   
\eprf

It follows that $cl(  {\G A} \meet ( W \times U) )$ still generates a ${\bb}$-solvable subgroup  $S$ of $G \times L$.
   Now $ cl( {\G A}) \meet (W \times U) = cl (  {\G A} \meet (W \times U)) $ so $S$ is generated by an open subset of 
$cl({\G A})$; hence $S$ is an open subgroup of $cl({\G A})$.   In particular, $C =  cl( {\G A})^0  =S^0 \leq S$ and so $C$ is solvable.   
   
In particular the first projection $\pi_1(C)$ is  soluble.   Since $C$ is   normalized by ${\G A}$, $\pi_1(C)$
is normalized by ${H}$.  By assumption (i), $\pi_1(C) = (1)$.  Thus $C = (1) \times I$ for some $I$, which is clearly
a closed subgroup of $L$.  $I$ is normalized by $\pi_2(\G) Z(L)$, which is dense in $L$ by assumption.  So $I$ is normal
in $L$; and $I$ is soluble.

  \begin{center}{*}    \end{center}
   
It remains only to prove that the induced homomorphism $f_1: {H} \to L/I$ is discrete-on-compact.    
By Claim 1,  $f \inv(I)$ is soluble, and is normalized by ${H}$; hence  by assumption (i), $f \inv(I)=(1)$ and so  $f_1$ is
injective.   So we are in position to use   \lemref{doc} (3); we thus need to show, for one compact neighborhood $W' \subset G$, that  $f_1(W' \meet {H})$ is discrete in $L/I$.          

Let us write also ${\G} I$ for ${\G} (1 \times I)$.   
 Since $cl({\G A})$ contains $cl({\G})$ and $ (1 \times I)$, and is closed; and on the other hand $A \subset I$, we 
 have $cl({\G A}) = cl({\G}I)$.
 
 We saw that $cl({\G A})^0=  1 \times I$.   The group $cl({\G A})/(1 \times I)$ 
  is thus totally disconnected.  So it has a profinite open subgroup, that we can write as $P' / (1 \times I)$, with  $P'$   an open subgroup    of $cl({\G A})$.     
     The image of $P'$ in $L/I$, namely $\pi_2(P') / I$,   is also profinite, in particular totally disconnected;    but by the Von-Neumann-Cartan theorem, being a closed subgroup of $L/I$, it is itself a Lie group.  
   Thus $\pi_2(P')/I$ is finite; so $ \pi_2 \inv(I) \meet P'$ has finite index in $P'$; being closed, it is an open subgroup,
   hence open in   $cl({\G A})$ too.    Let $P=\{a: (a,1) \in P' \}$; then $P \times I$ is a clopen subgroup of $cl({\G} I)$.

   

   Let us now show that $f_1(W')$ is discrete.    Otherwise, some subsequence of distinct points of $f_1(W')$ has a limit point
      in $L/I$; and it follows that some sequence of distinct points from $(W')(W') \inv$ converges to $1$ in $L/I$.
    Refining further ,using compactness of $W'$ and of $W''=(W')(W') \inv$, we find 
      a sequence  $(a_i)$  of elements of ${H} \meet W''$,  converging to some $a \in W''$,    such that $f_1(a_i) \to 1$.
      Let $U'$ be a precompact open neighborhood of $1$ in $L$; then $U'  I / I$ is an open subset of $L/I$, containing $\bar{b}$;
      so $f(a_i) \in U' I$ for almost all $i$; $f(a_i) = u_i c_i$, with $c_i \in I$ and $u_i \in U'$.  We have $(a_i, u_i c_i) \in \G$,
      so $(a_i,u_i ) \in {\G} I$.   Since $U'$ is precompact, we may assume $u_i \to u \in cl(U')$.  Then $(a,u) \in cl({\G}I)$.
       Since $P \times I$ is open in $cl({\G} I)$,  almost all pairs $(a_i,u_i)$ lie 
   in a single coset of $P \times I$.  In particular, the $u_i$ lie eventually in a single coset
   of $I$; this implies that modulo $I$ they are finite in number.   This contradicts the assumption of distinctness of the 
   elements $f_1(a_i) = f(a_i) I = u_i I$, and 
   ends the proof of Claim 2.  

\eprf

 \begin{center}{*}    \end{center} 

We return to the proof of the Proposition.  
The composed map ${H} \to L \to L/A$ is a  group homomorphism, injective by Claim 1.  Since $I/A$ is solvable, $f \inv(I) $ is
 solvable.     By   (i)  we have:  
 $ f \inv(I)=1 $.    Let $A_1=I$  and let $f_1: {H} \to L/A_1$ be  the composition ${H} \to L \to L/A_1$.    So $f_1$ is injective.  By  Claim 2,  $f_1$ is discrete-on-compact.
         
 We have now achieved the main goal, of replacing the quasi-homomorphism $f$ by an actual homomorphism $f_1$,
 still  discrete-on-compact.   Moreover $f_1$ is injective.  We still wish to change it further so that the target Lie group becomes semisimple
 and centerless.  To this end we repeat the  procedure inductively,  constructing a sequence $A_1 \leq A_2 \leq \cdots $ 
 of closed normal subgroups of $L$, such that $f$ induces a  discrete-on-compact, injective homomorphism $f_k: G \to L/A_k$.

 Assume $A_k,f_k$ have been defined and satisfy these conditions.   We proceed as follows.
  In case $L/A_k$ has a nontrivial center,   let $A_{k+1}/A_k$ be the center of $L/A_k$.
 Then the composed homomorphism $f_{k+1}: {H} \to L/A_{k+1}$ is  discrete-on-compact by Claim (2),
 applied to $f_k$.

  The kernel of $f_{k+1}$ is an abelian normal subgroup of   ${H}$, hence trivial by assumption (i).  
 
 Assume then that $A_{k+1}/A_k$ has trivial center.
 If $L/A_k$ has a nontrivial normal compact subgroup, choose $A_{k+1}$ to be such a subgroup.  Then $f_{k+1} :{H} \to  L/A_{k+1}$ 
 remains discrete-on-compact trivially, since a pullback of a compact from $L/A_{k+1}$ to $L/A_k$ is compact.
 The kernel, being the pullback of a compact, is discrete in $G$, and normal in $H$, hence must by trivial by (ii).  (Note
 $f_k$ is an injective homomorphism at this point, i.e. $K=1$.)  
   
 Assume therefore that $L/A_k$   has no such normal subgroups. 
 
 If $L/A_k$ has a nontrivial normal closed connected abelian subgroup, let $A'/A_k$ be such a subgroup.  
 Then $A'/A_k$
 is a connected locally compact abelian group so it is isomorphic to $\Rr^n \times B$ with $B$ a maximal compact subgroup.
 If $B'$ is a $G$-conjugate of $B$, then $BB'$  is also compact, so $BB'=B$.  It follows that $B$ is normal in $G$,
so should have been chosen in the previous step, unless $B=1$.   Thus $B=1$ so $A'/A_k$ is a vector group.
In this case the hypotheses of the Proposition apply to $G,{H}$ with $f_k: {H} \to L/A_k$ 
and $A'/A_k$ in place of $f$ and $A$.
By Claim 2 (now with $K=1$), there exists a closed, solvable normal subgroup $A_{k+1}$ of $L$ containing $A'$,  
such that the induced map $f_{k+1}: {H} \to L/A_{k+1}$ is discrete-on-compact.         Again the kernel of $f_{k+1}$ is soluble and so trivial by assumption (i).  

Assume $L/A_k$ has no normal subgroups of the types considered so far.  If $L/A_k$ has  a nontrivial normal closed   abelian subgroup $A'$, the connected component of $A'$ is trivial.  Since $A'$ is a Lie group, it is discrete, and so central, returning to a previously considered case.  
 
 The remaining possibility is that $L/A_k$ has no nontrivial closed normal subgroups that are either abelian or compact.  When this
 occurs the construction terminates.
 
 Note that $\dim(L/A_{k}) \geq \dim(L/A_{k+1})$; equality holds only in case $A_{k+1}/A_k$ 
   is discrete; this  cannot recur twice in a row, since if $A_{k+2}/A_{k+1}$ and $A_{k+1}/A_k$ are both
 discrete then so is $A_{k+2}/A_k$, and it is central in $L/A_k$ by connectedness, so $A_{k+2}$ would have been chosen
 already at stage $k$. Thus $\dim(L/A_k)$   descends and strictly so at every other step, so the construction proceeds at most $2 \dim(L)$ times.  
 
If $A_k$ is the last stage, then $f_k: {H} \to L/A_k$ is a   discrete-on-compact homomorphism to a  centerless Lie group, with no normal abelian subgroups, hence semisimple, and with no compact normal subgroups.   
  \eprf

  When $\G$ is a discrete subgroup of a connected Lie group $G$, Claim 2 of
 \propref{ctp} is the  ``Zassenhaus Lemma''   of  \cite{auslander}.   
 Curiously, though Auslander explicitly invokes an effective version of (**),  
  his proof appears to make no use of it; he uses rather the exponential shrinking within $A$ due to (2).   
    In the  approximate situation, (2) cannot be used in this way since 
  upon iterated approximate commutation, the $A$-component remains bounded but need not 
  approach $0$.    We used instead, after some rearrangement, the convergence of the elements $u_i$ towards  $1$.

\pb
 
 \section{Arithmeticity, recognition, extensions}  \label{sec8}
 
 Let $K$ be a global field.
      Let $\Aa=\Aa_K$ be the adele ring of $K$.   If $S$  is a 
     finite set   of places of $K$, we may write
     $\Aa = {\Aa_S} \times \Aa'_S$, where ${\Aa_S}= \Pi_{v \in S} K_v$, and $\Aa'_S$ is the restricted direct product
     of the $K_v, v \notin S$.    Let $\pi_S: \Aa \to {\Aa_S}$ be the projection.

 \begin{defn}[Arithmetic approximate lattice]  \label{arithmodelset}   Let $K$ be a global field, $H$  
    a connected, almost simple algebraic $K$-subgroup  of $SL_n$.   
       Let $C$ be a compact open neighborhood of $1$ in $H(\Aa'_S)$, and      
       \[ \Lam = \pi_S( H(K) \meet ({\Aa_S} \times C) ). \]
         If $\psi: G \to H({\Aa_S})$ is a  topological group isomorphism,   
     we refer to $\psi \inv(\Lam)$ as an {\em arithmetic approximate lattice}.  
       
   An {\em arithmetic Meyer set}  (in $H({\Aa_S})$) is a subset of an arithmetic approximate lattice, commensurable to it. \end{defn}

 \begin{rem}  If one demands that $S$ contain all archimedean primes $v$ where
 $H(K_v)$ is not compact,  one may (without changing the commensurability class) choose $C$ to be a subgroup;
 in that case one obtains the definition of an arithmetic lattice.    (\cite{margulis} p.1)\end{rem}
 
 \begin{rems}  \label{8.2}  
  \begin{enumerate}
 \item If interested only in the approximate subgroup $\Lam$, one can take a finite set $S''= S \sqcup S'$;  let $\Oo_{K,S''}$ be the ring of $S''$-integers in $K$;  let $C$ be a compact neighborhood of $1$ in $H(\Aa_{S'})$
 and let $\Lam = \pi_S(H(\Oo_{K,S''}) \meet ({\Aa_S} \times C))$.    One can take $S'$ to consist just of archimedean primes, if desired.
 Since $K$ embeds into ${\Aa_S}$,  and also into $\Aa_{S'}$, we have a map  $f: H(\Oo_{K,S''})  \to H(\Aa_S')$,
 showing incidentally that $\Lam$ is laminar.
  
   \item  On the other hand,  consider a pair $(\Lam,\wl)$, where 
 $\wl$ is any subgroup of $G$ commensurating $\Lam$ and containing $\Lam$.  If $(G,\Lam)$ is arithmetic,
 laminarity is also shown by a homomorphism $f: \wl \to H(\Aa_K)$, with $   \Lam \in \chi(f)$,
 and $f(\wl) \leq H(K)$.  
  This is because
 $\wl \leq H(K)$ (\lemref{A18}), and by \defref{arithmodelset}, $\Lam = j \inv(C)$ where $j: \wl \to H(K)$ is the inclusion.
 
  For this, if $\wl$ is not finitely generated, a finite product would not suffice. 
  
  Note here that if $C$ is chosen to contain $H(\Oo_{K_v})$ whenever $v$ is non-archimedean, then $f(\wl)$
  contains $H(\Oo_K)$.
  
  Take the Weil restriction of scalars $H'$ of $H$ from $\Oo_K$ to  $\Zz$.  \footnote{$\Oo_K$ is isomorphic to  $\Zz^n$
  as an abelian group, and so a choice of basis gives a polynomial interpretation of the ring $\Oo_K$ over the ring $\Zz$, 
  with universe $\Zz^n$.  Any affine scheme $W \subset \Aa^N$ over $\Oo_K$ becomes interpreted as $W^* \subset \Aa^{Nm}$
  over $\Zz$, in such a way that $W(\Oo_k)$ can be identified with $W^*(\Zz)$.}   We obtain a homomorphism $h: \wl \to H'(\Aa_{\Qq})$
  with $ \Lam \in \chi(f) $, and   $H'(\Zz) \leq h(\wl) \leq H'(\Qq)$.   
  
  \end{enumerate}\end{rems}

  For simplicity we restrict to groups of adjoint type below.  This allows us to use the simpler definition of `arithmetic lattice'
  in \cite{margulis}, p.1; and it is convenient to have centerless groups.  
     For classification up to commensurability this is not a real issue;
 If $G$ is a group with finite center $Z$, $\pi: G \to G/Z$ the projection, then any $X \subset G$ is commensurable with
 $\pi \inv (\pi(X))$.

 \begin{thm} \label{arith1}  
     
Let $G$ be a finite product of    noncompact, adjoint, 
         simple  algebraic groups over local fields.    Then any  approximate lattice $\Lam_0$  in $G$ is commensurable with a
         product $\Lam$ of lattices and arithmetic approximate lattices of direct factors of $G$. \smallskip

\noindent{{\bf Supplements}} \begin{itemize}
\item We may choose $\Lam \subset \la \Lam_0 \ra $. 
\item          In case no lattices appear in the product,  $\Lam_0$ is an arithmetic Meyer set.
\end{itemize}
             \end{thm}
      
 \prf       
    By \corref{semisimple1}, $\Lam_0$ is commensurable to a  cut-and-project set of a lattice $\G \leq G \times \sl$.
  Here $\sl$ is a centerless semisimple Lie group, and can be viewed as the connected component (of finite index) in a semisimple
   adjoint algebraic group over $\Rr$. 
    
Passing to a finite 
index sublattice of $\G$ if necessary,  
we can express $\G$ as a product of irreducible lattices $\G_i \leq G_i \times \sl_i$, where $G=\Pi_i G_i$
and $\sl= \Pi_i \sl_i$;  in the  sense of \cite{margulis}, p. 1.   If $U_i$ is a compact open neighborhood of $1$ in $\sl_i$, 
the cut-and -project sets in $G$ for $\G$ and $\Pi_i U_i$ is the product of the corresponding sets for $\G_i$ and $U_i$.
Thus $\Lam_0$ is commensurable to the cut-and-project set of $\G$ with respect to $\Pi_i U_i$.  From this it is also 
clear   that $\Lam_0$ is  commensurable to $\Pi_i (\pi_i(\Lam_0))$; moreover evidently $\G \subset \Pi_i (\pi_i(\Lam_0))$.
Thus proving the result for each $\Lam_i$ will suffice, so 
 we may restrict attention to a single such factor; renaming, we take $\G$ to be an irreducible lattice
in $G \times \sl$.   

  In case $\sl $ is trivial, $\G$ is a lattice in $G$, commensurable with $\Lam_0$, and there is nothing further
to show.  Assume from now on that $\sl$ (and $G$) are nontrivial.  Then,
$G \times \sl$ has rank at least $2$, and so Theorem   1 of \cite{margulis}, p. 2 applies:  in the notation of \defref{arithmodelset}, we may identify $G$ with $H(K_S)$,    $\G$ with $H(K)$ as embedded diagonally in $H(K_S) \times H(\Aa'_S)$
and $\Lam$ with  $\pi_1( \G \meet (H(K_S) \times C) )$, with $C$ a compact open neighborhood of $1$ in $H(\Aa'_S)$.  
This yields the arithmeticity of $\G$.

 It remains to show  that   $\Lam_0$,  is  contained  in a cut-and-project set  $\Lam$ of $\G$.   
   We have in any case
 $\Lam_0 \subset \Lam F_0$ for some finite set $F_0$.  Note that one can take $F_0 \subset \Lam_0 \Lam$. 
 By \lemref{commlem},   $F_0 \subset \Lam_0 \Lam \subset   \comm(\Lam)$.
 

By Borel density (\secref{boreldensity}), $\pi_1(H(K) \meet (H(K_S) \times C))$   is Zariski dense in $H$.  (It is here
 that we use the assumption that $G=H(K_S)$ has no compact factors.)   By \lemref{A18}, 
  $\comm(\Lam) \leq H(K)$.   
So   $F_0 \subset   H(K)$.  
  Now we may view $F_0$ as a subset
  of $H(\Aa'_S)$ as well; and it is clear that $F_0C$ is compact and  $\Lam_0 \subset \pi( H(K) \meet (K_S \times F_0C) )$.
  
  \eprf

 \begin{cor}\label{char-p}   Let $G$ be a semisimple algebraic group over a local field of positive characteristic, or a finite product of such groups.  Then every approximate lattice in $G$ is commensurable with a lattice.  \end{cor}

 \prf  By \thmref{arith1}, any approximate lattice in $G$ is commensurable with a product of lattices, and arithmetic approximate subgroups.  In the arithmetic case, the field $K$ will also be of positive characteristic (say by Borel-Tits), and so 
  the complementary $H(\Aa_S')$ will be totally disconnected.  One can thus take the compact open neighborhood
  $C \subset H(\Aa_S')$  to be a subgroup; and in this case, the   cut-and-project set is a lattice. \eprf

 \ssec{Recognition}

  We will find that an abstractly semisimple  locally compact group, possessing a  strictly dense  
 approximate lattice, is itself of adelic origin.         The method of proof 
  involves two uses of  complementarity, relating two  locally compact groups by bouncing off an
   approximate subgroup of a finite dimensional  Lie group.   
 Margulis arithmeticity is invoked, this time also via the commensurator criterion.   Abstract semisimplicity is only used
 in ordre to invoke \thmref{arith1}; laminarity of the approximate lattice could be assumed instead.  
 
As mentioned in the introduction,  \thmref{arith2} and \remref{arith2b} essentially go over earlier results of 
  Bader, Furman and Sauer \cite{bfs}.   They use lattices in place of approximate lattices, and make the
semisimplicity assumption only on the archimedean factor, but assume the lattice embeds into that factor.  
Using \thmref{arith1}, it is possible deduce a version of \thmref{arith2} from \cite{bfs}. We retain it as an example of the flexibility of the approximate lattice framework.

%



  \begin{thm} \label{arith2}  Let $G$ be an abstractly semisimple  locally compact group.  
   Let ${\lambda}$ be a  strictly dense  commensurability class of approximate lattices in $G$.  Then 
   there exists a semisimple algebraic group $H$ over $\Qq$, and a subgroup $M$ of $H(\Qq)$ containing
   $H(\Zz)$, such that $G$ is compactly isogenous to the   closure  $T$ of $M$ in 
    $H(\Aa_{\Qq})$.   
   \end{thm}
       
      \prf   
      By  \thmref{discrete1}   there exists a connected semisimple real Lie group $\sl$ and a lattice
$\Gamma \leq G \times \sl$  such that $\pi_1: \G \to G$ is injective,
 $\pi_2: \G \to \sl$ is injective with dense image; and with $\wl= \pi_1(\Gamma)$ and
$f= \pi_2 \circ \pi_1 \inv$, we have ${\lambda}=  \chi(f)$.   Moreover $\sl$ is centerless has no nontrivial  compact normal subgroups.   

  By the strict density of ${\lambda}$, $\wl$ is dense in $G$.       
   Let  $\wo =f(\wl)$.  Then  $\G$ is also
  the graph of a map $f \inv: \wo \to \wl \subset G$.
 
 Now let $W$ be a 
  compact neighborhood of $1$ in $G$, and  let  $\Om= f(W \meet \wl)$.  By  \lemref{latticecoherence},
 applied on the right this time, 
 $\Om$  is an approximate lattice in $\sl$.  

 We now apply \thmref{arith1} to $(\sl,\Om)$.  We find a direct product decomposition $\sl=\Pi_{i=1}^k L_i$, and approximate
 lattices $\Om_i \leq L_i$ with $\Pi_i \Om_i$ commensurable with $\Om$, such that each $\Om_i$ is an irreducible lattice,
 or an arithmetic approximate lattice.  Let $\pi_i: \sl \to L_i$ be the projection.

Since $f$ has dense image, so does $\pi_i f$ for each $i$.  As $\Om_i$ is discrete,   $\pi_i \wo$ and $\Om_i$ 
cannot be  commensurable.   Hence in case $\Om_i$ is a lattice, while it may have rank one,  it is still arithmetic,
 by the commensurator criterion of Theorem 1 of \cite{margulis}  (p.2; note finite generation is automatic in characteristic zero, according to  point (b)).  
 
 The fundamental commensurability class of the locally compact group $G$ is preserved under $\wl$-conjugation;
hence the commensurability class of $\Omega$ is preserved under $\wo$-conjugation.     Thus $\wo \subset \comm(\Pi_i \Om_i)$.   Let $\wo_i = \pi_i(\wo)$.  It follows that $\wo_i  \leq \comm(\Om_i)$. 

By \remref{8.2} (2),   there exist    homomorphisms $h_i:  \wo_i  \to H_i(\Aa_{\Qq})$ with image between $H_i(\Zz)$
and $H_i(\Qq)$, and with
$\Om_i \in \chi(h_i)$.
 Let $H = \Pi_i H_i$ and $h= \Pi_i h_i$; then $\Om \in \chi(h)$, and $H(\Zz) \leq h(\wo) \leq H(\Qq)$.  Let 
 $M=h(\wo)$ and $T=cl(M)$.   Then $h: \wo \to T$ is a homomorphism with dense image and $\Om \in \chi(h)$.

But we already know a homomorphism with dense image from $\wo$ to a locally compact group and {canonical commensurability class $\Om$}
 namely $f \inv:  \wo \to G$.   By \propref{uniquelc}, $T$ and $G$ are compactly isogenous, by an isogeny respecting $\wo$.
 
\eprf

  \begin{center}{*}    \end{center}

 Let us refer to abstractly semisimple locally compact groups  $T$ of the form appearing in \thmref{arith2}
 as   {\em Margulis  groups}.   

   \begin{rem}  \label{arith2b} Up to compact isogeny,
 and up to subgroups of profinite index, any Margulis group is itself a restricted product of semisimple groups over $\Rr$ and the $\Qq_p$, with respect to the open subgroups $\Zz_p$.  \end{rem}
 
\prf  By the proof of \remref{assrems} (\ref{bt1}), after factoring out a compact normal subgroup, 
 we have $T= T_{na} \times T_{ar}$, where  $T_{ar}$ is a connected semisimple Lie group,  
 and $T_{na}$ is the closure of $\wl$ in $H(\Aa_{fin})$, where $A_{fin} =  \Pi'_p  H( \Qq_p)$ are the finite adeles.    It suffices to prove the statement for $T_{na}$.  
   
  $T_{na}$  contains the closure of $H(\Zz)$ in $H(\Aa_{fin})$, which in turn contains an open subgroup of $\Pi_{p} H(\Zz_p)$.
  Hence $T_{na}$ is open in $H(\Aa_{fin})$.  We may take $H$ to be centerless, and a product of groups $H_1,\ldots,H_m$
  that are almost $\Qq$-simple; rewrite the restricted product as $\Pi'_q H_q$, where $q=(p,j)$ and $H_q = H_j(\Qq_p)$.
  Let $\pi_q: \Aa \to H_q$ be the projection.
 Let $I$ be the set of $q$ 
 such that $\pi_q(T)$ is not compact.  
 Let $T'$ be the projection of $T$ to $\Pi'_{q \in I} H_q$,.  The kernel of $T \to T'$ is compact.  For $q=(p,j) \in I$, the projection $\pi_q(T)$ is open and unbounded in $H_j(\Qq_p)$.  Hence by the theorem of Tits  in the title of \cite{prasad}, Theorem T,
 $\pi_q(T)$ contains $H_j(\Qq_p)^+$.   Since the groups  $H_j(\Qq_p)^+$ are simple and non-isomorphic,
 the image of $T'$ in any finite  product $\Pi_{q \in I_0} H_q$ must contain $\Pi_{q \in I_0} H_q^+$.  But $T'$ is open, hence closed, and so it contains the full restricted product $\Pi'_{q \in I} H_q^+$.  
 \eprf

  This gives in particular an analogue for approximate lattices of some of the results of \cite{oh}.

As I learned from 
 Emmanuel Breuillard, many examples of abstractly semisimple groups not embeddable in a 
  restricted product of semisimple algebraic groups are known, including certain automorphism groups of trees and non-$p$-adic 
  analytic pro-$p$ groups.  Such groups thus do not contain a strictly approximate lattice.

\ssec{A non-laminar approximate lattice} \label{rough}

Let $\Lam=SL_2(\Zz[1/p])$, a  lattice in $G=SL_2(\Rr) \times SL_2(\Qq_p)$.
We will define a central extension $\hat{G}$ of $G$ by $\Qq_p$ and 
a certain extension of $\Lam$
by a laminar approximate lattice  in $\Qq_p$, yielding an approximate lattice in $\hat{G}$ that is not {laminar}.
 
The group $SL_2(\Rr)$ has a nontrivial central extension $\widetilde{SL_2(\Rr)}$, determined by a bounded $2$-cocycle $\beta$ taking values 
 in $\{-1,0,1\} \subset \Zz$.   (See e.g. \cite{HPP2} for an o-minimal approach;   see also the general theorem on characteristic classes in \cite{gromov} 1.3, p. 231.)    We will at different points view $\beta$ as taking values in $\Rr$, in $\Qq_p$ 
 and in $\Zz[1/p]$.

\claim{}  For any index subgroup $\Lam'$ of $\Lam$,  $\beta$ induces  a nontrivial class in $H^2_{bdd}(\Lam',\Zz[1/p])$.

\prf  We will prove the stronger claim, that $\beta$  induces  a nontrivial class in $H^2_{bdd}(\Lam',\Rr)$.  

Consider the exact sequence:
 \[ 1 \to \Zz \to   \widetilde{SL_2(\Rr)}  \to SL_2(\Rr) \to 1\]
 If we change coefficients from $\Zz$ to $\Rr$, the extension
 \[ 1 \to \Rr \to \Rr \times_{\Zz} \widetilde{SL_2(\Rr)}  \to SL_2(\Rr) \to 1\]
remains nonsplit: the commutator subgroup of $\widetilde{SL_2(\Rr)} $ is the entire group so 
$\widetilde{SL_2(\Rr)} $ admits no nontrivial
homomorphisms to $\Rr$.  Thus   $\beta$ represents a nonzero cohomology class in $H^2(SL_2(\Rr),\Rr)$;
and a fortiori in  the bounded cohomology group $H^2_{bdd}(SL_2(\Rr),\Rr)$.  

Define $\beta_1$ on $G=SL_2(\Rr) \times SL_2(\Qq_p)$ by $\beta_1(x,y)=\beta(x)$.  Then $\beta_1$ represents 
 a nonsplit extension 
\[ 1 \to \Rr \to   \Rr \times_{\Zz}  \widetilde{SL_2(\Rr)}  \times SL_2(\Qq_p) \to G \to 1 \]
 so $\beta_1$ has a nonzero class in 
 $H^2_{bdd}(G,\Rr)$.   By Theorem 1.4 of \cite{monod}, since $2<4=2 rk(G)$, 
 and $\Lam$ is an irreducible lattice, the restriction map
 \[H^2_{bdd} (G,\Rr) \to H^2_{bdd}(\Lam, \Rr) \]
 is an isomorphism.    Hence 
 the restriction of $\beta$ to $\Lam$ remains nonzero in  $H^2_{bdd}(\Lam,\Rr)$, and by the same token, in $H^2_{bdd}(\Lam',\Rr)$
 for any finite index subgroup $\Lam'$ of $\Lam$ (which is also an irreducible lattice.)
 \eprf

Let $\hat{\Lambda}$ be the central extension of $\Lam$ by $ \Zz[1/p]$ determined by $\beta$, 
\[1 \ot \Lam \ot \hat{\Lam} \ot \Zz[1/p] \ot 1 \]\
Concretely, we can let $\hat{\Lam} = \Lam \times  \Zz[1/p]$, with multiplication 
\[(\lam_1,a_1)(\lam_2,a_2) = (\lam_1\lam_2,a_1+a_2+\beta(\lam_1,\lam_2))\]

Let $ {\delta}(\lam,a)=-a$, where now $a \in \Zz[1/p]$ is viewed as a real number.  Then $ {\delta}: \hat{\Lam} \to \Rr :[-1,1]$ is a quasimorphism.    Let
$\Delta$ be the (symmetric) approximate kernel  $\Delta = \{x: -2< {\delta}(x)<2 \}$.   By \propref{converse},
$\Delta$ is an approximate subgroup  of $\hat{\Lam}$.  We have $\pi_1(\Delta) =\Lam$,
so if a subgroup $S$ of $\hat{\Lam}$ contains an approximate subgroup commensurable to $\Delta$,
then $\pi_1(S)$ has finite index in $\Lam$, and thus $Z S $ has finite index in $\hat{\Lam}$
 (where $Z$ is the central image of $\Zz[1/p]$ in  $\hat{\Lam}$).     Let $\Lam'$ be the image of $S$ in $\Lam$; 
 then $\Lam'$ has finite index, and we saw in the Claim that $\beta$ is not trivial in  $H^2_{bdd}(\Lam')$.
 
Now if we had $\delta=b+h$ on $S$, where $b$ is bounded and $h$ is a homomorphism, then
for $x,y \in \Lam'$, choosing $s,t$ with $(x,s),(y,t) \in S$, we would have:  
$\beta(x,y) = \delta( (x,s)(y,t)) -\delta(x,s)-\delta(y,t) = b((x,s)(y,t)) - b(x,s)-b(y,t)$, so that $\beta=d b = 0$ as a class
in $H^2_{bdd}(\Lam')$.   Hence 
  $ {\delta}$ is not at bounded distance from a homomorphism on $S$.   
  
  By \propref{notsame},
the commensurability class of the approximate subgroup $\Delta$ of $\hat{\Lam}$ is not {laminar}.

 Now view $\beta $ as a $2$-cocycle with values in  $\Qq_p$, and construct in the same way  a central extension
\[1 \to \Qq_p \to \widehat{SL_2(\Rr)_{\Qq_p}}  \to SL_2(\Rr) \to 1.\]  

Here $\widehat{SL_2(\Rr)_{\Qq_p}}$ is obtained from the universal covering group of $SL_2(\Rr)$ by taking a fiber product with $\Qq_p$
over $\Zz$.  This also describes the locally compact group topology on $\widehat{SL_2(\Rr)_{\Qq_p}}$; it has an open subgroup isomorphic
to the inverse limit of coverings of $SL_2(\Rr)$ of order $p^n$.   

Taking products with $SL_2(\Qq_p)$, we obtain 
\[1 \to \Qq_p \to \hat{G}  \to G \to 1.\]  

We can view the universe of $\hat{G}$ as being $G \times \Qq_p$, with multiplication twisted by $\beta$ as above.  
Then the product measure is a Haar measure; translation by an element of $\hat{G}$ involves a certain transvection
on $\Qq_p \times G$, and transvections are measure preserving.  Thus the product measure is a Haar measure.

 $\Delta$ is discrete in $\hat{G}$:   suppose $(\lambda_i,a_i) \in \Delta$ and they approach a limit point  $(\lam,a) \in \hat{G}$.
 Them $\lambda_i \to \lam \in G$; since $\Lam$ is discrete in $G$, we have $\lam_i = \lam \in \Lam$ for large $i$;
 we may arrange by translation that $\lam_i = \lam = 1$.  In this case we have $a_i \to a \in \Qq_p$.  
 But $\{a \in \Zz[1/p]:  |a|_\Rr < 2 \}$ is discrete in $\Qq_p$.  Thus $a_i=a$ for large $i$.

 $\Delta$ is a lattice in $\hat{G}$:  
  if $C_1 \subset  G$ is a Borel set of finite measure with $C_1 \Lam = G$, then $(C_1 \times \Zz_p) \hat{\Lam} = \hat{G}$,
and of course $\Zz_p \subset \Qq_p$ is compact and hence of finite measure.
 
 Thus 
  $\Delta$ is a non-laminar approximate subgroup of
$\hat{G}$.

\ssec{Approximate lattice extensions}  \label{othergroups}

Let $G$ be a connected linear algebraic group over a local field $k$ (the discussion extends naturally to a finite product of such groups.)  Then $G$ has soluble radical $R=R(G)$; it is a soluble Zariski closed normal subgroup, with $\bG:=G/R$   semisimple.     One can either choose $R$ connected, or such that  $G/R$ is centerless; the latter will be a little more convenient for us.   
We define $G_{ss} = G/R$.

Trying to reformulate   the Bj\"orklund-Hartnick problem, given the amenable and semisimple answers, we may first ask:

\begin{question} \label{othergroups-q}  Is $\Lam/R$ above  always an  approximate lattice in $G/R$? 
For connected $G$ at least, I am fairly sure, a  slight modification (and simplification) of  \lemref{ctp} gives a positive answer.  
  \end{question}
  
N.B.:   Machado has answered this positively at least over characteristic zero local fields.  

\bigskip

Let us consider the case where $\bl:=\Lam/R$ is an approximate lattice; then we know it is a product of
arithmetic approximate lattices, and lattices.  In case $\la \Lam \ra$ is dense in $G$ (as is natural to assume),
$\la \bl \ra$ will be dense in $G/R$, and then all the components of $\bl$ will be arithmetic.  By \propref{minimalclass}, $\la \bl \ra$ contains an element $\Gamma$ of ${\la \bar{\om} \ra } _{\min}$; and by pulling back such an element, we may
reduce to the case that $\la \Lam \ra \in {\la \bar{\om} \ra } _{\min}$.  In this case, $\la \Lam \ra$ is an $S$-arithmetic group
(of course it is not a lattice of $G/R$, but of some other product of semisimple groups.)    Now we may further ask:

\begin{question}  \label{fine?}Let $\bG$ be a product of semisimple groups over local fields, and $\bl$ an approximate lattice in $\bG$,
with $\la \bl \ra$ isomorphic to a given $S$-arithmetic group $\wl$.  When is it true that for all $G$  
with $G_{ss}\cong \bG$,
 every approximate lattice $\Lam$ of $G$ with $\Lam /R = \bl$  is {laminar}?  
 \end{question} 
 
 
\ssec{}  \label{othergr}   We sketch a possible argument that the answer to \qref{fine?} is:  just when: \\
($\diamond$)
 $H^2_b(<\bl>,V) =(0)$ for all rigid actions of $<\bl>$ on Euclidean space $V=\Rr^n$, induced from continuous  action of $G$.
 
 If $\diamond$ holds,   taking into account that $H^2_b(<\Lam>,V) = H^2_b(<\bl>, V^R)$
(by the mapping theorem \cite{gromov} p.40,  extended to nonconstant coefficients, \cite{noskov}, cf. \cite{monod}.),
the laminarity follows from \propref{cohcor} (using \propref{minimalclass} for $H$.)   
 
 Conversely   if $H^2_b(\bl,V) \neq (0)$, let $\beta$ be a nontrivial bounded $2$-cocycle, and use it to form an extension $\hat{\Lam}$ of $\bl$ by the module $V$.      By the  mapping theorem   we have 
  $H^2_b (\hat{\Lam}) = H^2_b(\Lam)$ canonically since $V$ is amenable.  Thus 
 pulled back   to   $\hat{\Lam}$,  $\beta $ still represents a nontrivial class in $H^2_b(\hat{\Lam})$;
 but now it maps to $0$ as a class in $H^2(\hat{\Lam})$.    Hence by the  the discussion below \ref{cohcor},
 we find a nonlaminar approximate subgroup in a central extension of $\hat{\Lam}$.   It remains to 
 embed  this in an appropriate algebraic group, as in \ref{rough}, so as to produce a negative answer for this instance of Question \ref{fine?}.
 
 The bounded cohomology of $S$-arithmetic groups, as $\bl$ must be here,  is  studied in depth in  \cite{burger-monod}, \cite{monod}, providing I believe both
 instances where $H^2_b(\bl,V) $ is zero for all $V$, and otherwise.
 

%

%
%
%
%

\

  \appendix

  \pb
  
  \section{Approximate lattices} \label{approxlattices}
  
  In this Appendix, $G,H$ denote   second countable, locally compact groups; $\mu$ will denote a  left invariant measure 
  on $G$.  
  
  We develop the basics of approximate lattices under a   slightly different (and more elementary) definition than given in  \cite{BjH}.
 In the case of strong or cocompact approximate lattices, these results appear in \cite{BjH} and \cite{machado3}.   
  
  \begin{defn}  An approximate subgroup $\Lam \leq G$ is an {\em approximate lattice} in $G$ if it is discrete, and the equivalent conditions (1)-(3) of \propref{approxlattices1} hold.  \end{defn}
   
  \begin{prop}\label{approxlattices1}  
     Let $\Lam$ be a discrete  approximate subgroup.    Then  the  following properties are equivalent:
  \begin{enumerate} 
  \item $\Lam D = G$ for some Borel set $D$ of finite measure.
  \item  There exists a finite bound on $\mu(C)$, over all Borel sets $C  \subset G$ with $\Lam^2 \meet C C \inv= (1)$;
  i.e. such that multiplication is 1-1 on $\Lam \times C$.
  \item  There is no Borel $C \subset G$ with $\mu(C) = \infty$, and $\Lam^2 \meet C C \inv= (1)$.
  \end{enumerate}
 If these conditions hold, then $G$ is unimodular, i.e. $\mu$ is also right invariant; and so the three dual conditions are also equivalent.         
  \end{prop}
  
  \prf  Note that as $G$ is second countable, the discrete subset $\Lam$ is countable.  
  
  Assume (1).  Let $C$ be a Borel set, with     $\Lam^2 \meet C^2 = (1)$.   
  We have $C \subset \union_{a \in \Lam} C_a$, where $C_a = (aD \meet C)$.  By the invariance assumption, $\mu(C_a) = \mu(a \inv C_a)$.
  The sets $a \inv C_a$
  are disjoint Borel subsets of $D$; so
  \[ \mu(C) \leq \sum_{a \in \Lam} \mu(C_a) = \sum_{a \in \Lam} \mu(a \inv C_a)  \leq  \mu(D).\]
  
  (2) trivially implies (3). 
  
 Now assume (3).   By \lemref{borel} (applied to $X=\Lam ^2 \m (1) $, $Y=G$, with opposite multiplication) there exists a Borel $M$ with $MM \inv \meet \Lam^2 = (1)$,  and  $G=   \Lam^2 M$.  By (3),   $M$ has finite volume.     
  But $\Lam^2 \subset   \Lam F$
 for some finite $F$, so $G= \Lam F M$, and $F M = \union_{f \in F} fM$ has finite volume.

    To prove unimodularity,  let $\rho(g)$ be the constant, given by the uniqueness of the Haar measure,
    such that for all Borel $C$ we have $\mu (Cg) = \rho(g) \mu(C)$.  
    Pick a symmetric compact open neighborhood $C$ of $1$, with $\Lam^2 \meet C^2 = (1)$; this is possible since
  $\Lam^2$ is discrete.  Then 
   for any $g \in G$,  we have $(Cg) (Cg) \inv= C C \inv$, so $\Lam^2 \meet (Cg) (Cg) \inv= (1)$.  By (2) 
   there is an absolute bound on $\mu(Cg)$ over all $g$.  Hence $\rho$ is bounded, and the image of $\rho$ is a bounded subgroup of $\Rr^{>0}$, i.e. $1$.
  
  
  \eprf
  
  Unimodularity was also proved in \cite{BjH} for strong approximate lattices, in the sense defined there.
  
The following statement on a Borel set of approximate representatives will be essential.    See \cite{kst} for a study of much deeper
 questions along this line.
 
 \begin{lem}\label{borel}     Let $G$ be a discrete or second countable locally compact group, $Y$ a Borel subset, and $X$ a discrete subset of $G$ with $1 \notin X$, and $X=X \inv$.   A subset $B \subset Y$ is called {\em $X$-free} if  $B   \inv B \meet X =\emptyset$, and {\em maximally $X$-free} if it is maximal with this property, i.e. $ Y \subset BX$.    Then there exists a Borel set $B \subset Y$ that is maximally $X$-free.  
 \end{lem}
 
 \prf \   \\
\claim{1}    Assume $Y$ is a countable union of Borel subsets $Y_n$, and for each $n$, any Borel $Y' \subset Y_n$ has a maximally $X$-free Borel set.   Then we can solve the problem inductively:  having defined $B_m \subset Y_m$, let $Y_n' = Y_n  \setminus \union_{m<n}   B_m X$, and let
  $B_n$ be a maximally $X$-free subset of $Y_n'$.   Then $B=\union B_n$ will be maximally $X$-free in $Y$.
  
 In view of this, since $G$ can be written as a countable union of compact neighborhoods, we may assume $Y$ is precompact.
 
For any Borel $Z \subset Y$, $b \inv Z \meet X$ is finite, since $X$ is discrete in $G$ and $Y$ precompact.   For a given
finite subset $H$ of $X$, let \[Z_H = \{b \in Z:  b \inv Z \meet X = H \}.\]  

\claim{2} For any Borel $Z \subset G$, and any finite $H \subset X$, we can find
a maximally $X$-free Borel subset of $Z_H$, or of any given Borel subset of $Z_H$.

We prove   Claim 2  by induction on $|H|$.  For $|H|=0$, $Z_H$ is already $X$-free, so any $Z \subset Z_H$ is maximally
$X$-free in $Z$.  
Assume $|H|=n$.  Let $Z'=Z_H$.  For $b \in Z'$, we have  $b \inv Z' \meet X \subset H$.   We can thus divide $Z'$
into finitely many sets, according to the value of $b \inv Z' \meet X$.  For any proper subset of $H'$, the induction hypothesis
applies to $Z'_{H'}$.  Thus it suffices to deal with $Z'_H$.  If $Z'_H= (Z_H)_H$ is empty, there is nothing to prove.
Otherwise, $H$ is a finite subgroup of $G$.   To see this, let $a \in (Z_H)_H$ and let $x, y \in H$.  
Then $ax \in Z_H$, so $axy \in Z$.  It follows that $xy \in a \inv Z \meet X$; since $a \in Z_H$, we have $xy \in H$.
So $Z_H$ is a union of $H$-orbits,   $y\inv z \in X$ iff $y,z$ lie in the the same $H$-orbit.  So given a Borel $W \subset Z_H$,
it  suffices to find a set of representatives for the $H$-orbit equivalence relation restricted to   $W$.  
(for instance we can fix a Borel linear ordering on $G$, and choose the least element of $W$ in each $H$-orbit.)   
This proves the Claim for $(Z_H)_{H'}$ for any $H' \subset H$, and hence by Claim 1, it is valid for $Z_H$.

Now $X$, being 
a discrete subset of a second-countable space, is countable; so there are only countably many possibilities for the finite set
$b \inv Y \meet X$.   We thus have a partition of $Y$ into countably many Borel sets $Y_H$.  By Claim  (2), any 
Borel subset of each of these
has a Borel maximally free set, and hence by Claim 1, so does $Y$.     \eprf

\begin{lem} \label{commens}  If an approximate lattice  is contained in a  set   $\Lam$ with $\Lam^4$ discrete, 
  they are commensurable, and $\Lam$ is an approximate lattice too.
 \end{lem}
\prf   Let $X$ be an approximate lattice contained in $\Lam$.  
%
Let $U$ be a symmetric compact neighborhood of $1$, small enough so that $U^2 \meet \Lam^4 = (1)$.
Let $\{c_i: i \in I\}$ be a maximal collection of elements of $\Lambda$, such that the sets $Xc_i$   are disjoint.
If $x,x' \in X, i,i' \in I$, $u,u' \in U$  and $xc_iu = x' c_{i'}u'$, then $u' u \inv \in \Lam^4$ so $u=u'$; then by the disjointness
of the sets $Xc_i$   we have $i=i'$; so finally $x=x'$.  In particular, the sets $c_iU$ are disjoint; and 
  multiplication is 1-1 on $X \times (\union_{i \in I} c_iU)$.   Hence, using \propref{approxlattices1} (3), 
  \[  |I| \,  \vol(U) =  \sum_{i \in I} \mu(c_i U) = 
  \mu(\union_{i \in I} c_i U) < \infty.\]    
So $I$ is finite.  From maximality of $I$, it follows by the usual  argument
 that $\Lam \subset \union_{i \in I } X \inv X c_i $:   if  $d \in \Lam$, then $X c_i \meet Xd \neq \emptyset $
for some $i$, so $d \in X\inv X c_i$.  As we assumed $X$ is an approximate subgroup, $\Lam$ and $X$ are commensurable. 
Using condition  \propref{approxlattices1} (2) it is clear that $\Lam$ is an approximate lattice.

\eprf

\ssec{Complementarity}  \label{complementarity}

 Let  $f \vdots G \to H$ be a partial map  between topological spaces, i.e. a map $f: X \to H$ where $X \subset G$.
We will say that $f$ is {\em  discrete-on-compact} if the pullback of any (pre)compact set in $H$ is discrete in $G$.  
 
 $f$ is {\em precompact-to-discrete} if  whenever $C \subset G$ is compact, then $f(C \meet \wl)$ is discrete in $H$.

\begin{lem}\label{doc}  Let $G,H$ be   locally compact groups, $\wl$ a subgroup of $G$, and  $f: \wl \to H$  a group  homomorphism.

\begin{enumerate}
\item   If $f$ is discrete-on-compact, and $X \subset H$ is compact, then $f \inv(X)$ is  uniformly discrete.

\item   If $f$ is discrete-on-compact,  it is   precompact-to-discrete.

\item   Assume  $G$ is second countable, and $f$  has finite kernel.   Then the converse to (2) holds:  
if $f$ is precompact-to-discrete, then $f$ is discrete-on-compact.   In fact if $f(U \meet \wl )$ is discrete in $H$ for one    neighborhood  $U$ of $1$ in $ G$
then $f$  is discrete-on-compact. 
\end{enumerate}\end{lem}

\prf (1)  Clear since $f \inv (X X \inv)$ is discrete.

(2)  In a locally compact space, to prove discreteness of $f(C \meet \wl)$, it suffices to show that $f(C \meet \wl) \meet C'$ is finite for every
compact $C'$; this is clear since $f \inv(C')$ is discrete in $G$, so $C \meet f \inv(C')$ is finite.

(3)  Reversing (2), given a compact $C' \subset H$, we must show that $f \inv(C') $ is discrete;  it suffices
to show that no 1-1 sequence $(c_i: i \in \Nn)$ of elements of $f \inv(C')$ converges to any $c \in G$.   Indeed if $c_i \to c$
then  $\{c,c_i : i \in \Nn\}$ is
compact, so  $\{c_i: i \in \Nn \}$ is precompact, and 
$\{f(c_i): i \in \Nn \}$    is discrete in $H$;    hence intersects $C'$ in a finite set; since $f$ is finite-to-one,   $\{c_i: i \in \Nn \}$ is finite,
  in contradiction to the assumption that the $c_i$ are distinct.   

To prove the stronger version,  we may  assume
that $\wl$ is dense in $G$.     Let $U$ be a nonempty  open set with   $f(U \meet \wl)$   discrete; we will show that
$f(D \meet \wl)$ is discrete for any  compact subset $D$  of $G$.  For each point $d \in D$,
let $\g(d) \in \wl \meet   U d \inv$.    Then $D \subset \union_{d \in D}  \g(d) \inv U$; so a finite number of translates $\g_i U$ 
cover $D$, with $\g_i \in \wl$.  Hence $f(D \meet \wl) \subset \union_{i} f(\g_i) f(U \meet \wl)$; so $F(D \meet \wl)$ is a finite union of discrete subsets of $H$ and hence discrete in $H$.

\eprf

 \ssec{Discrete-on-compact homomorphisms and Meyer sets}   \label{d-on-c}
 
 Let $G,H$ be   locally compact groups; let $\pi_1: G \times H \to G, \pi_2: G \times H \to H$ be the projections.
 Given a discrete subgroup $\G \leq G \times H$  such that $\pi_1$ is injective on $\G$, let $\tl = \pi_1(\G)$ and ${f}=\pi_2 \circ (\pi_1|\G) \inv: \Lam \to H$.  Then ${f}$ is discrete-on-compact.    Indeed if $C \subset H$ is compact,
 and $U$ any compact   set in $G$, then ${f} \inv(C) \meet U $ is finite since $(U \times C) \meet \G$ is finite;
 hence ${f} \inv(C)$ is discrete in $G$.  
  Conversely, if ${f}: \tl \to H$ is discrete-on-compact, then the graph of ${f}$ is a discrete subgroup of $ G \times H$.

 The uniform part of the  proposition below  is Prop. 2.13 of \cite{BjH}.   We do not  assume in advance that $H$ is unimodular,
 so that if $G$ has an approximate lattice and $f$ exists, then $H$ is automatically unimodular.  
   
 \begin{prop}\label{latticecoherence} Let $G,H$ be second countable locally compact groups, $\wl $ be a subgroup of $G$,  $f: \wl \to H$ a discrete-on-compact homomorphism with dense image.  Let $\G$ be the graph of $f$, viewed as a subset of $G \times H$.  
 Let $\Lam = f \inv(C)$, where $C $ is a compact, symmetric subset of $H$ with nonempty interior.   Then:
 \begin{enumerate}
 \item  $\G$ is a   lattice iff $\Lam$ is an approximate lattice.  
 \item  $\G$ is a uniform lattice iff $\Lam$ is a  uniform  approximate lattice.  
  \end{enumerate}
 \end{prop}

 \prf    Assume $\Lam$ is an approximate  lattice.   Let $U$ be a Borel  subset of $G$ of finite volume,  with $ \Lam U = G$.
We will show that $ \G (U \times C^2) = (G \times H)$.   As $U \times C^2$ has finite volume (and is compact if $U$ is), 
it follows that $\G$ is a lattice, and is uniform if $U$ is compact.
 
 To see that  $(U \times C^2)\G = (G \times H)$,  let $(g,h) \in G \times H$.   By density, $f(a) \in hC$ for some $a \in \tl$.   So $f(a)C \subset hC^2$.
  Now $a \Lam \meet gU\inv \neq \emptyset $ since $ g \in a  G =    a \Lam U$.
 Let $b \in a \Lam \meet gU \inv$.  Then $g \in bU$; and $f(b) \in f(a)  f(\Lam) \subset hC^2$, so $h \in  f(b)C^2$.
 Thus $(g,h) \in (b,f(b)) (U \times C^2) \subset \G (U \times C^2)$. 
 
 Now assume $\Lam$ is not a lattice.  By   \propref{approxlattices1} (3), there exists a
 Borel $U \subset G$ with $\mu(U) = \infty$, and $\Lam^2 \meet U U \inv= (1)$.  In other words
 $f \inv(C)^2 \meet U U \inv = (1)$, so $(U \times C)(U \times C) \inv \meet \Gamma =  (1)$.   Thus $\Gamma$ is not a lattice.  
 
   In the uniform case, it remains to show that if  $\G$ is a cocompact lattice,
 then there exists a compact set of the form  $U \times C$, with $(U \times C) \G = (G \times H)$; 
this is clear since any compact subset of $G\times H$ is contained in a compact rectangle.   
\eprf

Let us call a closed, symmetric neighborhood of $1$ in a totally disconnected locally compact group $H$ {\em sensible} if it is contained in some compact open subgroup; more generally for locally compact $H$, call $C$ sensible if the image $C/H^0$
in $H/H^0$ is sensible.  In particular if $H$ is connected, then any symmetric compact neighborhood of $1$ is sensible.

\begin{prop}\label{sensible}  Let $G,H$ be locally compact, $\Gamma \leq G \times H$ a lattice projecting densely to $H$ and 1-1 to $G$.  Let $\wl= \pi_1(\G)$, $f = \pi_2 \circ \pi_1 \inv:  \wl \to H$.  
Let $C$ be a compact symmetric  subset of $H$ with nonempty interior, $\Lam=f\inv(C)$, and $X \subset \wl $ a   symmetric set  commensurable with $\Lam:=f \inv(C)$.    Then $X \subset f \inv(C')$ for some compact $C' \subset H$.  

 In case $C$ is sensible,  $\la X \ra$ contains a subgroup of $\la \Lam \ra$ of finite index.   Thus  the commensurability class of $\Lam$ belongs to $\la \Lam \ra$ but not to any infinite index subgroup.

   In particular
 for any two sensible $C,C'$, not only the approximate subgroups $f \inv(C), f \inv(C')$ but also the groups   generated by them are commensurable.
 \end{prop}
 \prf     We have $X \subset \union_{i=1}^m \Lam a_i$, for some $a_1,\ldots,a_m \in \wl$; let
$C' = \union _{i=1}^m C f(a_i)$.   Clearly  $X \subset f \inv(C')$.   At this point it is also clear that 
$X$ is an approximate lattice in $G$.    ($X^2 \subset f \inv((C')^2)$ and the later is commensurable with $f \inv(C)$, so $X$
is an approximate subgroup.  On the other hand $\Lam \subset \union_{i=1}^l Xb_i$ so if $\Lam M = G$ with 
$M$ Borel of finite measure, then $X (\{b_1,\ldots,b_l\} M) =G$ and $(\{b_1,\ldots,b_l\} M)$ has finite measure.)

We have $\Lam \subset \union_{i=1}^l Xb_i$.   Since $C \subset \union_{i=1}^l cl(fX) b_i$, 
it follows from Baire category that $cl(f(X))$ has nonempty interior. 
 Thus $cl(f(X))^2$ contains a neighborhood of $1$ in $H$.  

Now assume $H$ is connected.   Then $H=\union_n cl(f(X))^n $.  So $f(\la X \ra)$ is dense in $H$.
By \propref{latticecoherence} (1),   the graph $\G'$ of $f | \la X \ra$ is a lattice in $G \times H$.   
By the (easier lattice case of) \lemref{commens},  
 $\G'$ has   finite index in $\G$.   Projecting to $G$ we see that $\la X \ra$ has finite index in $\wl$.
 
 Finally, assume more generally that $C$ is sensible.   We may replace $X$ by the commensurable set $X^2 \meet \Lam^2$
 (\ref{icc}); so we may assume that $X \subset \la \Lam \ra$.    By assumption, $C/H^0$ is contained in an compact open subgroup $O/H^0$.   For statements about $\Lam=f \inv(C)$,  and subsets of $\la \Lam \ra$,  we may as well replace $H$ by $O$; so 
 we may assume $H/H^0$ is compact.  In this case every  open subgroup of $H$ has finite index in $H$.  In particular
 $H':=cl(f(\la X \ra))$ has finite index in $H$.  As above,   the graph $\G'$ of $f | \la X \ra$ is a lattice in $G \times H'$,
 hence also a lattice in $G \times H$; and we conclude as in the previous paragraph.   \eprf

 \begin{rem} It can also be shown that   if  $\G \leq G \times H$ is a   lattice, with dense projection to $H$, then 
there exists  a Borel $M \subset G$ of finite volume and a {\em compact} $C \subset H$    with $(M \times C) \G = G \times H$.
(The proof starts similarly to \lemref{borel}, then invokes 4.2, 4.5 of \cite{kst}.)  
Thus  the lack of uniformity in an irreducible lattice of a product can be relegated to any  one of the factors.  

This gives an alternative proof of  \propref{latticecoherence} (2).    

The same fact on representative sets is  moreover true for approximate lattices. 
 \end{rem}   
 
\ssec{Borel-Wang density}  

Borel density is addressed in \cite{BjHS} for strong and for uniform approximate lattices.   We will give a direct proof for  approximate lattices in the sense we use here.  

If $G=\Pi_i G_i$ is  a product or restricted product of algebraic groups over possibly different fields (viewed with this structure, disregarding possible isomorphisms between the fields), we will say that a subset $X \subset G$ is {\em Zariski dense} if
the projection to each $G_i$ is Zariski dense.


\begin{prop}  \label{boreldensity}   Let $G$ be an algebraically connected 
 semisimple algebraic group  with no anisotropic factors over a local field.    Then every approximate lattice in $G$ is
  Zariski dense.

  Likewise for finite or restricted products of such groups over possibly different local fields; an approximate lattice is not commensurably contained in 
a product of Zariski closed subgroups, other than the entire ambient group. \end{prop}

\prf  Let $\Lam$ be an approximate lattice in $G$, and let $C$ be the Zariski closure.   We have $1 \in \Lam= \Lam \inv$ so $ 1 \in C = C \inv$.  We have
$\Lam^2 \subset \Lam F$ for some finite set $F$, so $C^2 \subset CF$.   Let $C_1$ be an irreducible component of $C$ of maximal dimension; then $C_1^2 \subset C_2 F$ for some component $C_2$ of $C$, and it follows that $C_1^2 = C_2 f$ for some $f \in F$;
from this   it follows easily that  $C$ is contained in a finite union of cosets $Hg_i$ of a Zariski closed subgroup $H$ of  $G$,
$g_i \in N_G(H)$, and that $\Lam$ is commensurable to $H$.  

  If $\Lam$ is a closed approximate subgroup of a locally compact group $G$, say $\Lam$ has {\em property S} if
  for every symmetric neighborhood $U$ of the identity in G and every element $g \in G$, there is a positive integer $n$ such that
   $g^n \in U \Lam^2 U$.    Equivalently, the sets $g^n U \Lam$ are not pairwise disjoint.  
   For closed subgroups, this agrees with the definition in \cite{wang}.     
        
Now  if $\Lam $ is an approximate lattice, (S) holds.  To see this take $U$ symmetric, with $U^2 \meet \Lam^2 = (1)$.
If the  sets  $g^nU$
    are not  pairwise disjoint,   the    condition holds trivially (using only $1 \in \Lam$).   Thus we may assume they are
    pairwise disjoint, so $\union g^n U$ is a set of infinite measure.
By \propref{approxlattices1} (3),  $(\union g^n U)\inv (\union g^n U)  \meet \Lam^2 \neq (1)$; so for some $n \geq m$,
$ (g^n U)\inv  g^m U  \meet \Lam^2 \neq (1)$.   We have $n>m$ since   $U^2 \meet \Lam^2 = (1)$.
Thus $U g^{n-m} U \meet \Lam ^2 \neq (1)$  and $g^{n-m} \in U \Lam^2 U$, proving (S).   

 Of course this implies 
that (S) holds for any subgroup $H$ containing $\Lam$.  
According to the Borel-Wang density theorem for groups, \cite{wang} Corollary (1.4), if $H$ is Zariski closed this implies that $H= G$.

In case $G$ is a finite or restricted product of such groups $G_i$, in a  product $\Pi_i G_i$, since $\Lam$ has property (S), it is 
clear that each projection $\pi_i(S)$ has property (S) too; it follows by definition that $\Lam$ is Zariski dense.  
  \eprf
%
\ssec{Product decomposition} 

The hypotheses of the following Proposition hold notably (using \propref{boreldensity}) when  $G_1,G_2$ are centerless, semisimple groups with no anisotropic factors over a local field.

%

\begin{prop}\label{productdec}  Let $G_1,G_2$ be  centerless locally compact groups, both satisfying the chain condition for centralizers,  and such that an approximate lattice is not contained in any proper centralizer subgroup.  
  Let $\Lam$ be an approximate lattice in $G=G_1 \times G_2$.
Let $\pi_i: G \to G_i$ be the projection, and let $\Lam_i = G_i \meet \Lam^2$ (we identify $G_1$ with $G_1 \times 1$, notationally.)
Then the following conditions are equivalent:  \begin{enumerate}
\item $\pi_2(\Lam)$ is discrete in $G_2$.
\item  $\Lam_1$ is a approximate lattice in $G_1$.
\item  $C_{G_1}(\Lam_1) = 1$.
\item[(1')] $\pi_1(\Lam)$ is discrete in $G_1$.    
\item  $\Lam$ is commensurable with a product of an approximate lattice in $G_1$ and an  approximate lattice in $G_2$.    \end{enumerate}
\end{prop}
\prf  
($1 \to 2$)    Certainly $\Lam_1$ is discrete in $G_1$, since $\Lam$ is discrete.    It is an approximate subgroup of $G_1$
by \secref{icc}.   To show that it is an 
   approximate lattice, it remains to show that conditions (2) or (3) of \lemref{approxlattices1} hold, i.e. that $\Lam_1$ has finite covolume.
   We will check this at the same time for $\pi_2(\Lam)$.  Let $\mu_1,\mu_2$ be  Haar measures on $G_1,G_2$, and let $\mu=\mu_1 \times \mu_2$   be the product Haar measure on $G$.  

 Let $E_i$ be any Borel subset of $G_i$, with $E_1^2 \meet \Lam_1 =(1)$, $E_2^2 \meet (\pi_2 \Lam)^2 = (1)$.
 Then $(E_1E_2)^2 \meet \Lam^2 =(1)$:   if $e_i \in E_i^2$ and $e_1e_2 \in \Lam^2$, then $e_2 \in E_2^2 \meet (\pi_2 \Lam)^2 $; by assumption $e_2  = 1$, so $e_1 \in \Lam^2 \meet G_1 = \Lam_1$, hence  $e_1 = 1$.   
   Thus  $E_1E_2$ has finite (indeed bounded) measure.  As $\mu(E_1E_2)=
   \mu_1(E_1)\mu_2(E_2)$, by fixing some small open neighborhood $E_2$ of $1$, we see that there is an absolute bound
   on $\mu_1(E_1)$ and in any case it cannot be infinite.  Hence $\Lam_1$ is an approximate lattice in $G_1$. 

($2 \to 3$)   Assume (2).  By assumption,  if $a \in G_1$ and $C_{G_1}(a) \neq G_1$, then $C_{G_1}(a)$
 cannot contain $\Lam_1$; so $a \notin  C_{G_1}(\Lam_1)$.   Thus $C_{G_1}(\Lam_1)$ consists only of central elements of $G_1$, hence equals $(1)$.  
 
($3 \to 1'$)   By the chain condition on centralizers in $G_1$ and (3), there exists a finite $F \subset \Lam_1$ with $C_{G_1}(F) = (1)$.  We now use
the argument of \lemref{centralcase}.  Let $a_i \in \pi_1(\Lam)$, and suppose $a_i \to a \in G_1$, and let $f \in F$.  We have
$a_i f a_i \inv \in  \Lam^4 $, and $a_i f a_i \inv \to af a \inv$.  
   Now $\Lam^4$ is discrete in $G$, so $a_i f a_i \inv = af a \inv$ for large $i$.    Thus $a \inv a_i \in Z_G(f)$ for each $f \in F$.
   So $a =a_i$ for large $i$.   This proves the discreteness of $\pi_1(\Lam)$.
   
 By repeating the circle, it is now clear that (1-3) are equivalent, and also equivalent to their duals (1'-3') replacing $G_1,G_2$.
  Thus $\Lam_1,\Lam_2$ are both approximate lattices,  hence so is $\Lam_1 \times \Lam_2$.  Clearly 
  $\Lam_1 \times \Lam_2 \subset \Lam^2$.  Since $\Lam^8$ is discrete, by \lemref{commens}, $\Lam_1\times \Lam_2$ is
   commensurable with $\Lam$.  (As is  $\pi_1 \Lam \times \pi_2 \Lam$, incidentally.)

  This gives (4).  The implication ($4 \to 1$) is clear. 
\eprf

By  \propref{boreldensity},  the assumption that proper centralizer subgroups not contain a lattice is true for groups $G_1,G_2$ satisfying the hypothesis there.  
 
In case only $G_1$ satisfies the chain condition on centralizers, we still have the implications $1 \to 2 \to 3 \to 1' \to 2' \to 3'$,
and $2+2'  = 4$.  (So $2 \to 4$.)  

The Proposition justifies a definition.   Let  $G$ be a connected, centerless, semisimple groups with no anisotropic factors over a local field.

\begin{defn}  \label{irrapp}
An approximate lattice $\Lam$ of $G$ is {\em irreducible} if there is no closed normal subgroup $N$ of $G$ 
such that the image of $\Lam$ in $G/N$ is discrete. \end{defn}

\begin{cor} \label{irrappcor} Let $G$ be a finite product of noncompact simple algebraic groups over local fields.     Then every approximate lattice of $G$ is commensurable to   a product $  \Pi_{i=1}^k \Lam_i$,
where $G=\Pi_{i=1}^k G_i$ and $\Lam_i$ is an irreducible approximate lattice in $G_i$.
\end{cor}

\ssec{Commensurator in algebraic groups}

\begin{lem} \label{commensurator2}  Let $L$ be a field, $K$ a  perfect subfield,  ${H}$    a centerless,  connected algebraic
group over    $K$.    Let $g \in H(L)$ and assume $g \inv H(K) g  \,   \meet H(K)$ is Zariski dense in $H(K)$.  Then
$g \in H(K)$.  \end{lem}
\prf   Let  $D:=g \inv H(K) g  \,   \meet H(K)$.  As $D$ is Zariski dense, $C_{H}(D) = (1)$;  hence $C_H(F)=(1)$ for some finite set $ \{d_1,\ldots,d_k\}  \subset D$;
we have $gd_i g \inv = b_i $ with $d_i,b_i \in {H}(K)$; this uniquely determines $g$ in ${H}(L^{alg})$, so $g$ is definable in
$ACF_K$ and hence lies in ${H}(K)$
\eprf

 \begin{lem} \label{A18}  Let $H$ be a connected algebraic group over a perfect field $K$.  Let $\Lam$ be
 an approximate lattice, Zariski dense in $H(K)$.  Let $L$ be a field extension, and let $C = \comm(\Lam)$
 be the commensurator of $\Lam$ in $H(L)$.  
    Then  $C \leq H(K)$. \end{lem}
  \prf As noted in \secref{icc},
 if $A,B$ are commensurable approximate subgroups, they are both commensurable to $AA \meet BB$, hence
 using connectedness of $H$,  $AA \meet BB$ is Zariski dense in $H$.   Let $g \in \comm(\Lam)$.  
 Let $A=\Lam$, and $B=g \inv \Lam g$; since
 $A \subset H(K)$ we have $A^2 \subset H(K)$, and similarly $B^2 \subset g \inv H(K) g$; so if  
  $g \in H(L)$ commensurates $\Lam$, then $H(K) \meet g H(K) g \inv$   is Zariski dense in $H$.  
  By  \lemref{commensurator2}, this
   implies that  $g \in H(K)$.   \eprf

\begin{question} \label{unimodularityq} If $ G_1 \leq G_2$ are locally compact groups, and $\Lam$ is an approximate lattice of $G_2$ contained
in $G_1$, is $G_1$ unimodular?  Equivalently, is $\Lam$   
 an approximate lattice of $G_1$?   The case $G_1 \normal G_2$ is easy  (in this case, carring a finite Haar measure, $G_2/G_1$ is compact).   But the general case seems to require a little care (if not a surprise).     \end{question}

  \begin{question}  Formulate an adelic version of \corref{semisimple1}.\end{question}  
     
\begin{question}\label{bjhquestion} Clarify the relation to the Bj\"orklund -  Hartnick definitions of approximate lattice.  

 N.B.:  Simon Machado has shown that a strong approximate lattice is also a finite covolume approximate lattice; and
 in turn a finite covolume approximate lattice is  an approximate lattice in the sense of \cite{BjH}.   It follows from \cite{BjH} and the results here
 that in semisimple groups, the   notions of strong and  finite covolume approximate lattice coincode; but the general case remains open.  
\end{question}

 \begin{question}  Let $\Lam$ be  an approximate lattice in  a non-amenable locally compact group.  
 Could $\Lam$ be definably amenable in a reasonable language, or even for that matter outright amenable in the sense of \cite{machado3}?  For 
 co-compact $\Lam$,  \cite{BjH} answer negatively for metric amenability.  
\end{question}

 \pb
 
\section{Categories of approximate spaces.} \label{catcomcor} \label{categories}
 
 We will use no category theory in this paper, at least not explicitly.   
Nevertheless it seems useful to describe a number of categories that  clarify the difference between  the automorphism group of the core, and the two Lascar groups.

\ssec{The categories $\cosp$ and  $\appsp$}  

Let $X,Y$ denote  locally compact, Hausdorff topological spaces.  A {\em proper correspondence} is a closed subspace
 $S \leq X \times Y$ such  that the projections $\pr_{S,X}: S \to X$, $\pr_{S,Y}: S \to Y$ are proper maps.  If $S \subset X \times Y$
 and $T \subset Y \times Z$ are proper correspondences, we let $S^{tr} = \{(y,x): (x,y) \in S\}$,  $T \circ S = \{(x,z) \in X \times Z: (\exists y \in Y) (x,y) \in S \wedge (y,z) \in Z \}$.      In case $X=Y=Z$ and $S=T$, we write this as $S^{\circ 2}$, and similarly define higher {\em iterates}.

 The objects of  $\cosp$ are pairs $(X,A_X)$ with  $X$ a locally compact space, and $A_X \subset X^2$ a 
 symmetric, reflexive  proper correspondence.     A {\em morphism} $(X,A_X) \to (Y,A_Y)$ is a 
  proper correspondence $S \leq X \times Y$ such that $\pr_{S,X}=X$ and  
 $S \circ  A_X   \circ S^{tr} \subset A_Y$.     
  
   \ssec{The category $\appsp$}     In $\appsp$, the objects are the same as in $\cosp$;
 the second condition on a morphism is weakened to: $S \circ  A_X   \circ S^{tr} \subset A_Y^{\circ n}$ for some $n$;
 and two morphisms  $S,T: (X,A_X) \to (Y,A_Y)$ are {\em identified} if  $S \subset A_Y^m \circ T $ for some $m$, and vice versa.

The map $X \mapsto (X, Id_X)$ extends to a fully faithful functor of the usual category of locally compact topological spaces, into  $\appsp$.    The more general objects of $\appsp$, of the form $(X,S)$, can be viewed as the space
$X$ where the points $x \in X$ are only known with some `experimental error'; if   $x$ is the result of a measurement, 
the actual point may be any $x'$ with $(x,x') \in S$.     The traditional approach in model theory, in this situation, was to factor
out the equivalence relation generated by $S$,  losing track of the topology; or the smallest closed equivalence relation containing $S$, retaining a locally compact space but usually a much smaller one.  Viewing $(X,S)$ as an object in $\appsp$ is an alternative way of `factoring out $S$';   it preserves, in general, much more information.     Only in case $S^{\circ m}$ is an equivalence relation  $E$ for some $m$, does the object $(X,S) $ of $\appsp$ reduce to the quotient space $X/E$.  
 
In some cases   we can do better, obtaining an object of $\cosp$  canonically.


 \ssec{The categories $\cogr$ and  $\appg$}  Here the objects are pairs $(G,K)$ where $X$ is a locally compact topological group
 and $K$ is a compact, normal, symmetric subset of $G$.    A morphism $(G,K) \to (H,L)$ is a   subset $M \leq G \times H$,
 with $(1,1) \in M$, $M=M \inv$ and $M^2 \subset  M (K \times L)$.  
 
 We have a natural functor $\cogr \to \cosp$, mapping $(G,K)$ to $(G,A_K)$ where $A_K = \{(x,y) \in G^2:  xy \inv \in K\}$.
 
 The category $\appg$  is obtained from this in a similar way to $\appsp$.
 
 An example is the object $(T,K)$ of $\appg$, where $T = \Rr^2/ \Lam$ is a torus, $\Lam \cong \Zz^2$ a torus in $\Rr^2$,
 and $K = [0,1] \times \{0\}$ is a one-dimensional interval.   Typically, $K$ generates a dense subgroup of $T$, 
 and so $K/\la K \ra$ is an indiscrete space, and $K / cl(\la K \ra)$ is a point.  On the other hand $(T,K)$ retains considerable information, in fact $T$ may be recovered from it.

\pb

\section{Group-theoretic terminology}  

 We   collect here a few   group-theoretic definitions and statements that we will need.
  
\noindent
\begin{defn}\label{groupdefs}
Let $G$ be a topological group.
\begin{enumerate}

\item  We denote the connected component of the identity in $G$ by $G^0$.  $[x,y] $ is the commutator $xyx\inv y \inv$.
 For $X,Y \subset G$, the  commutator {\em set} $[X,Y]$ is defined as $\{[x,y]:  x \in X, y \in Y \}$. 
 The group generated by $X$ is denoted $\la X \ra$.  (The commutator subgroup is thus $\la [G,G] \ra$.)
 $C_G(A)$ denotes the centralizer of $A$ in $G$, $C_G(A):= \{g \in G: (\forall a \in A) (ga=ag) \}$. 
 $Z(G):=C_G(G)$.    A {\em centralizer subgroup} is one of the form $C_G(X)$ for some $X$; equivalently, a subgroup $H$ with $C_G(C_G(H))=H$.     We say that a subset of a group is {\em normal} if it is conjugation invariant.


\item   We will   denote the image of a subset $X$ of  $G$ modulo a normal subgroup $N$ by $X/N \subset G/N$.

\item   The {\em commensurator} $\comm(\Lam)$ of a symmetric set $\Lam \subset G$ is the group of elements $g \in G$ with $g\Lam g \inv$ commensurable to $X$.  When $\Lam $ is an approximate lattice, the commensurator always contains
the group $\wl$ generated by $\Lam$.

\item  A {\em vector group} is a topological group isomorphic to $\Rr^n$ for some $n$.  A {\em Lie group}
is a group with a real analytic structure.   (We will consider $p$-adic groups too,  as algebraic groups with their natural $p$-adic topology, but will not use the term Lie groups for them.)  

\item
By a {\em basic abelian group} we mean a compactly generated locally compact Abelian group with no nonzero compact subgroups.  Such a group has the form $A=\Zz^l \times \Rr^n$, $l,n \in \Nn$.  

\item   A normal subgroup $N$ of  $G$ will be called   {\em rigid} if it is closed in $G$,
and can be written $N=DE$ with $D \cong \Zz^n$, $D$ central in $G$, and $E$ a Euclidean space on which $G$ acts by rigid motions; thus the action of $G$ by conjugation factors through $O_n=Aut(E_n)$.

 \item  A precompact subset of $G$ will be called   {\em rigid} if it is contained in some rigid normal subgroup.
 
 \item  $G$ is {\em compactly generated} if some compact $X \subset G$ generates $G$ as a pure group.
 Given such a compact generating set $X$,  let $f: Free(X) \to G$ be the natural map from  the free   group generated by $X$ to $G$;  then $G$ is {\em compactly presented} (via $X$) if   $\ker f $ is generated (algebraically) by words in $X$ of some bounded length $r$.   

\item $cl(S)$ dentoes the closure of $S$, when $S$ is a subset of a given topological space.

 \item We will say that a set $W \subset G$ is discrete in $G$ if it has no accumulation points in  $G$.
 $W$ is  {\em uniformly discrete} if $1$ is not an accumulation point of $W \inv W$.
 
 \item  A subset $W \subset G$ is {\em precompact} if the closure of $W$ is compact.

\item \label{bounded} Let a group $G$ act on a metric space $X$ by isometries.
We say that a subset $W$ of $G$ is {\em bounded} (with respect to this action) if for some $a \in X$, $Wa :=\{w(a): w \in W \}$ is
contained in a  ball of finite radius.  Equivalently, for {\em all} $a \in G$, $Wa$ has finite radius.   (Since $d(a,b)=d(ga,gb)$, if 
$d(a,ga)\leq \beta$  we have
$d(b,gb) \leq d(b,a)+d(a,ga)+ d(ga,gb)  \leq \beta+2 d(a,b)$.)

  \item  \label{quasimorphism-def}  A quasimorphism  is an quasi-homomorphism $f: G \to \Rr:K$ for some compact $K \subset \Rr$.  It is 
   {\em homogeneous} if $f(x^n)=nf(x)$.   Any quasimorphism is uniquely a sum of a bounded function into $\Rr$ and a homogeneous quasimorphism (\cite{kotschick}.)  Note that   a homogeneous quasimorphism $f$ is symmetric, i.e. $f(x \inv)=-f(x)$, since $f(x)+f(x^{-1})$ is also homogeneous, and bounded.  We also have $f(y \inv x y)=f(x)$, since $|f(y \inv x^ny) - f(x^n)| $ is bounded.

\item \label{restricted}   Let $G_i$ be a topological group, for $i \in I$.
  Let $C_i$ be an open subgroup of $G_i$.  The restricted product $\Pi' G_i$ is defined to be the subgroup $G$ of $\Pi_i G_i$ 
consisting of elements $(a_i)_{i \in I}$ with $a_i \in C_i$  for all but finitely many $i \in I$.    If each $G_i$ is locally compact
and almost every $C_i$ is compact, then $G$ is locally compact.  

\end{enumerate}

\end{defn}

\begin{lem}\label{rigid}   Let $A$ be a basic abelian group, and $H$ a group of automorphisms of $A$.  
Let $D$ be a dense subset of $A$, and assume $Hd$ is  precompact, for each $d \in D$.
    Then $H$ preserves an innder product on $A^0$; and one can write $A$ as a direct sum  $A = A^0 \oplus B$,   with $B \cong \Zz^k$ fixed pointwise by a finite index subgroup $H_1$ of $H$.   
\end{lem}

\prf   We have $A^0 \cong \Rr^n$.    Since $D$ is contained in no proper subspace of $A^0$, one can find a basis  $c_1,\ldots,c_n \in D$ for $A^0=\Rr^n$.  
 So $H c_1 ,\cdots,H c_n $ are precompact.   Let $G$ be the subgroup of $GL_n(\Rr) = End(A^0)$
preserving $cl(c_i^H)$ for each $i$.    Then it is clear that $G$ is a compact subgroup of $GL_n(\Rr)$;  hence $G$
preserves an inner product on $\Rr^n$ (namely the  integral  over $G$ of the $G$-conjugates of the standard inner product).  
Thus $H$ too preserves this inner product; and so every $H$-orbit on $A^0$ is precompact.

Next consider the action on $A/A^0 \cong \Zz^k$.   Since $A/A^0$ is discrete, the image of $D$ is all of $A/A^0$.    
Let $e_1,\ldots,e_k$ be a $\Zz$-basis for $A/A^0$.  Then the  orbit $E_i$ of $e_i$ is finite;
we have restriction maps $H \to Sym(E_i)$; the common kernel of these maps is a 
  finite index subgroup $H_1$ of $H$,  fixing each $e_i$ and hence fixing $A/A^0$.   Lift $e_i$ to a point $d_i$ of $D$.  
  
 \claim{}  Some point $b_i$ of $d_i+ A^0$  is fixed by $H_1$.
 \prf    Since $H_1$ fixes $A/A^0$, it  leaves invariant the $A^0$-coset $d_i+A^0$.   Also, $H_1 d_i$ is precompact.   The
  closed convex hull of $H_1 d_i$ is a  compact, convex, $H_1$-invariant subset of $d_i+A^0$. 
 Since the Euclidean distance is preserved, so is Lebesgue measure, and so also the center of mass  of this convex set; it is  a fixed
 point for the action of $H_1$. \eprf
  
  Let $B= \la b_1,\ldots,b_k \ra$.  Since the image of $B$ modulo $A^0$ is $\Zz^k$, we also have $B \cong \Zz^k$.  The decomposition we promised is achieved.

 \eprf

Discrete groups whose  every conjugacy class is finite  (called FC-groups) 
were studied by I. Schur, R. Baer, and B.H. Neumann, see \cite{neumann}.
This was later extended to  locally compact groups, and precompact conjugacy classes, by Ushakov 
and 
 Tits \cite{tits}, with emphasis on the connected case, and by
 \cite{GM}, \cite{wu-yu} in general.    Roughly speaking, such groups are extensions of abelian groups by compact ones.   
 We require a version describing   a single precompact conjugacy class.   
  Thanks to Yves de Cornulier for explaining the relevant   literature; previously, the deduction of this local version was a little more elementary but considerably more elaborate.

%
%
\begin{prop} \label{lc1cc}   Let $H$ be a  locally compact topological group. 
Let $K $ be a normal precompact subset of $H$. 
Then there exist closed normal subgroups $N \leq A \leq H_1 \leq H$ with $H/H_1$ finite, 
$N$ compact,  with $A/N$  basic abelian and rigid in $H_1/N$,    such that $\la K \ra \meet H_1  \leq A$.
Moreover, $N,A \leq cl(\la K \ra)$. 
   \end{prop} 
   
\prf
   
Let $G=cl(<K>)$ be the closed subgroup of $H$ generated by $K$.  Then $G$ is normal in $H$.   
   
By \cite{co-ha} 2.C.6, a locally compact  group  is compactly generated iff every open subgroup $U$ generates it,  once it is    augmented by some finite set.  
This condition clearly holds for $G$, since the compact set $cl(K)$ generates a dense subgroup;   
it suffices to add to $U$  a finite set $F$ such that $cl(K)  \subset FH$. 
Thus $G$ is compactly generated.  

Let $B(G)$ be the union of all precompact conjugacy classes of $G$.   Then $B(G)$ contains $<K>$, so it is dense.  
By \cite{cornulier}, Theorem 2.8, B(G) is closed; hence must equal $G$.  

 Now \cite{GM}, Theorem 3.20, applies; by this theorem,  $G$  has a characteristic compact normal subgroup $P$, with $G/P$ basic abelian.   Now $H$ acts on $G/P$ by conjugation.  Let $D$ be the image of $<K>$ in $G/P$; any $d \in D$ has   pre-compact orbit 
 under the action of $H$.  By \lemref{rigid}, a subgroup $H_1$ of finite index in $H$ acts rigidly on $G/P$.   Let $A=G \meet H_1$,
 $N=P \meet H_1$; then all the statements are clear.  \eprf

 \end{document}